\documentclass[12pt,leqno]{article}

\usepackage{amsthm}

\newcommand{\bm}[1]{{\mbox {\boldmath $#1$}}}

\usepackage{amsmath}
\usepackage{amscd}
\usepackage{amsmath}
\usepackage{amssymb}
\usepackage{latexsym}
\makeatletter

\@addtoreset{equation}{section}
\makeatother

\newtheorem{thm}{Theorem}[section]
\newtheorem{pr}[thm]{Proposition}
\newtheorem{df}[thm]{Definition}
\newtheorem{lm}[thm]{Lemma}
\newtheorem{cor}[thm]{Corollary}

\newtheorem{rmk}[thm]{Remark}
\newtheorem{ex}[thm]{Example}
\newtheorem{df-lm}[thm]{Definition-Lemma}

\newcommand{\Spec}{{\rm Spec}}

\newcommand{\Ql}{\overline{\mathbb{Q}_\ell}}
\newcommand{\rk}{{\rm rk}}
\newcommand{\sw}{{\rm Sw}}
\newcommand{\dt}{{\rm dimtot}}
\newcommand{\Proj}{\mathbb{P}}
\newcommand{\Zl}{\overline{\mathbb{Z}_\ell}}
\newcommand{\F}{\mathbb{F}}
\newcommand{\cyc}{{\rm cyc}}
\newcommand{\Z}{\mathbb{Z}}

\input xy \xyoption{all} \CompileMatrices

\begin{document}

\title{Characteristic epsilon cycles of $\ell$-adic sheaves on varieties}
\author{Daichi Takeuchi\thanks{Department of Mathematics,
Institute of Science Tokyo,
2-12-1, Ookayama, Meguro-ku, Tokyo, 152-8551, Japan. 
Email: 
\texttt{daichi.takeuchi4@gmail.com}}}
\date{}
\maketitle

\begin{abstract}
Let $X$ be a smooth variety over a finite field $k$. 
Let $\ell$ be a rational prime number invertible in $k$. 
For an $\ell$-adic sheaf $\mathcal{F}$ on $X$, 
we construct a cycle supported on the singular support of $\mathcal{F}$ 
whose coefficients are $\ell$-adic numbers modulo roots of unity. 
It is a refinement of the characteristic cycle $CC(\mathcal{F})$, in the sense that 
it satisfies a Milnor-type formula for local epsilon factors. 
After establishing fundamental results on the cycles, we prove a 
product formula of global epsilon factors modulo roots of unity. 
We also give a generalization of 
the results to varieties over the perfections of finitely generated fields. 
\end{abstract}

\section{Introduction}\label{intro}
Let $k$ be a finite field and let $\ell$ be a 
prime number invertible in $k$. 
Let $X$ be a smooth projective variety over $k$, and 
$\mathcal{F}$ be an object of the triangulated category $D^b_c(X,\Zl)$. In this article, 
 we use the term ``$\ell$-adic sheaf'' 
  to refer to  an object of this category. 

For an object ${\cal F}\in D^b_c(X,\Zl)$, the $L$-function $L(X,\mathcal{F};t)$ is defined as  
the infinite product 
\begin{equation*}
L(X,\mathcal{F};t)=\prod_x\frac{1}{{\rm det}
(1-{\rm Frob}_xt^{\deg(x/k)},\mathcal{F}_{\bar{x}}\otimes_{\Zl}^L\Ql)}\in\Ql[[t]], 
\end{equation*}
where $x$ runs through closed points of $X$, $\deg(x/k)$ denotes the degree of the extension of the residue field at $x$ over $k$, ${\rm Frob}_x$ is the geometric Frobenius at $x$, and $\mathcal{F}_{\bar{x}}$ is the stalk of $\cal F$ at a geometric point $\bar{x}$ lying above $x$. 
 Using his theory of \'etale cohomology, Grothendieck showed that this function is indeed a polynomial (rather than merely a formal power series) and satisfies the  functional equation 
\begin{equation*}
L(X,\mathcal{F};t)=\varepsilon(X,\mathcal{F})t^{-\chi(X_{\bar{k}},\mathcal{F})}
L(X,\mathbb{D}_X\mathcal{F};t^{-1}), 
\end{equation*}
where $\mathbb{D}_X$ denotes the Verdier duality functor. Let $H^i:=H^i(X_{\bar{k}},{\cal F}\otimes_{\Zl}\Ql)$, where $\bar{k}$ is an algebraic closure of $k$. Then the invariants 
$\chi(X_{\bar{k}},\mathcal{F})$ and $\varepsilon(X,{\cal F})$ can be described as 
\begin{equation*}
\chi(X_{\bar{k}},\mathcal{F})=\sum_i(-1)^i
\dim H^i, \quad \varepsilon(X,{\cal F})=\prod_i\det(-{\rm Frob}_k,H^i)^{(-1)^{i+1}}. 
\end{equation*}
 It is a long-standing problem in arithmetic geometry to express these invariants in terms of some invariants which can be computed locally on $(X,{\cal F})$. For the Euler characteristic $\chi(X_{\bar{k}},\mathcal{F})$, Beilinson and T. Saito made significant  progress on this problem.  Saito \cite{Sai17} 
 constructed the characteristic cycle $CC({\cal F})$, a cycle supported on the singular support $SS({\cal F})$ defined by Beilinson \cite{Bei}, and they proved 
the formula 
\begin{equation}\label{indexf}
\chi(X_{\bar{k}},\mathcal{F})=(CC({\cal F}),T^\ast_XX)_{T^\ast X},
\end{equation}
which expresses $\chi(X_{\bar{k}},\mathcal{F})$ as an intersection number (we will return to this formula after Theorem~$1.3$ below). When $X$ is a curve, this  specializes to the Grothendieck--Ogg--Shafarevich formula. 

The main subject of this article is the term $\varepsilon(X,{\cal F})$, called {\it the global epsilon factor}. When $X$ is a smooth  projective curve, this problem of describing this invariant was settled by Laumon \cite[(3.2.1.1)]{Lau}. He showed that the global epsilon factor decomposes as the product of the local epsilon factors of $\cal F$ defined at the  closed points of $X$; 
 this is commonly known as the product formula of global epsilon factors. 
  The purpose of this article is to attempt to generalize this formula for $\varepsilon(X,{\cal F})$ to higher-dimensional varieties,  following the 
strategy of Beilinson-Saito employed in the study of 
$\chi(X_{\bar k},\mathcal F).$

Let $X$ be a smooth variety over $k$. Let $\cal F$ be an $\ell$-adic sheaf on $X$ (i.e., an object of $D^b_c(X,\Zl)$). The singular support $SS({\cal F})$ of $\cal F$, defined by Beilinson \cite{Bei}, is  a closed subset of the cotangent bundle $T^\ast X$ of $X$.  Strictly speaking, $SS$ is defined in loc.~cit. for \'etale  sheaves with finite coefficients.  For an $\ell$-adic sheaf $\cal F$, we define its singular support by considering ${\cal F}\otimes^L_{\Zl}\overline{\mathbb{F}_\ell}$ and using the invariance of the singular support under extension of scalars. See Definition \ref{SSCCadic} for the precise construction. 

In this article, we introduce a cycle ${\cal E}({\cal F})$, called {\it the epsilon cycle of ${\cal F}$}, which is a cycle supported on $SS({\cal F})$ with coefficients in {\it $\ell$-adic numbers modulo roots of unity}. We also prove various formulae for ${\cal E}({\cal F})$; among them, we obtain a higher dimensional generalization of Laumon's product formula, up to roots of unity. 

Let us explain the content of this article more precisely. Set 
\[
\Theta:=\Zl^\times\otimes_{\mathbb{Z}}\mathbb{Q}, 
\]
which we regard as  the group of $\ell$-adic numbers modulo roots of unity. Let $SS({\cal F})=\cup_a C_a$ be the decomposition into irreducible components. As observed by Beilinson in \cite{Bei}, each $C_a$ has the same dimension as $X$. Therefore $[C_a]$ may be viewed as elements of $Z_n(T^\ast X)$, the group of $n$-cycles on $T^\ast X$, where $n=\dim X$. 

The epsilon cycle ${\cal E}({\cal F})$ is then defined as an element of $\Theta\otimes Z_n(T^\ast X)$, and 
takes the form 
\begin{equation*}
{\cal E}({\cal F})=\sum_a\xi_a\otimes[C_a]
\end{equation*}
with $\xi_a\in \Theta$. 
In the following, we describe several fundamental properties of ${\cal E}({\cal F})$ established in this article.

First, we explain that ${\cal E}({\cal F})$ satisfies a Milnor-type formula for local epsilon factors. As in the case of characteristic cycles, this property uniquely characterizes ${\cal E}({\cal F})$, and we in fact take this as the definition of the epsilon cycle. To state this result, we  recall a basic property of $SS({\cal F})$. 

Let $f\colon X\to\mathbb{A}^1_k$ be a function on $X$, and let $df$ denote the differential of $f$, regarded as a section $X\to T^\ast X$.  Beilinson shows in \cite{Bei} that $f$ is universally locally acyclic relatively to $\cal F$ provided that $df$ does not meet $SS({\cal F})\subset T^\ast X$.

Let $f\colon X\to\mathbb{A}^1_k$ be a function, and suppose  that there is a closed point $x\in X$ such that $df$ meets $SS({\cal F})$ only at $x$. Following \cite{Sai17}, we call such a point  {\it an isolated $SS({\cal F})$-characteristic point of $f$}. Let $x$ be such a point. Then  $f$ is universally locally acyclic relatively to $\cal F$ away from $x$ (as recalled above). The vanishing cycles complex of $(f,{\cal F})$ is thus supported at $x$; we denote it by $R\Phi_f({\cal F})_x$. This is a bounded complex of $\ell$-adic representations of the local field of $\mathbb{A}^1_k$ at $f(x)$. 

In the Milnor-type formula given below, we consider its local epsilon factor 
\[\varepsilon_0(\mathbb{A}^1_{k,(x)},
R\Phi_f(\mathcal{F})_x,dt)\in\Zl^\times,\]
 where $t$ is the standard coordinate of the affine line. This notation for the local epsilon factor is the one used just after \cite[(3.1.5.6)]{Lau}. 

The following is the Milnor formula satisfied by ${\cal E}({\cal F})$. 
\begin{thm}(Theorem \ref{epcygenmil})\label{ayoo}
Let $X$ be a smooth variety over a finite field $k$. 
Let $\mathcal{F}$  be an element of $D^b_c(X,\Zl)$. Write 
$SS({\cal F})=\bigcup_aC_a$ for the  decomposition into irreducible components of the singular support. 
 Then there exists a unique cycle 
\begin{equation*}
{\cal E}({\cal F})=\sum_a\xi_a\otimes[C_a]\qquad \xi_a\in\Theta=\Zl^\times\otimes_{\mathbb{Z}}\mathbb{Q},
\end{equation*}
called the epsilon cycle, 
which satisfies the following property.  For any diagram of $k$-schemes 
\begin{equation*}
X\xleftarrow{j}U\xrightarrow{f}\mathbb{A}^1_k
\end{equation*}
with $j$ \'etale and any isolated ($SS({\cal F})\times_XU$)-characteristic point $u\in U$ of $f$, we have an equality 
\begin{equation*}
\varepsilon_0(\mathbb{A}^1_{k,(u)},
R\Phi_f(\mathcal{F})_u,dt)^{-1}=
(\mathcal{E}(\mathcal{F}),df)_u^{\deg(u/k)}:=
\prod_a\xi_a^{\deg(u/k)\cdot(C_a,df)_{u}}
\end{equation*}
in $\Theta$. Here $(C_a,df)_{u}$ denotes the intersection number of $j^\ast C_a:=C_a\times_XU$ with $df$ at $u$. 
\end{thm}

In particular, epsilon cycles can be computed \'etale locally on $X$. 

Let us make a technical comment on the definition of local epsilon factors. To define local epsilon factors, one has to choose and fix a non-trivial additive character of $k$. However, if another character is chosen, the resulting local epsilon factor  differs from the original one only by a root of unity. Thus the two are equal when considered as elements of $\Theta$. For this reason,  we omit the choice of an additive character in the statement of the theorem and in what follows. 

In \cite{Sai17}, Saito defines the characteristic cycle $CC({\cal F})$, which is an integral  combination of the components $C_a$: 
\begin{equation*}
CC({\cal F})=\sum_an_a[C_a],
\end{equation*}
where $n_a\in \Z$. By the compatibility of local epsilon factors with unramified twist and by the characterizations of ${\cal E}({\cal F})$ and $CC({\cal F})$ in terms of Milnor-type formulae, we obtain the following relation  between them in $\Theta\otimes Z_n(T^\ast X)$:
\begin{equation}\label{SSCCrel}
\mathcal{E}(\mathcal{F}(1))-\mathcal{E}(\mathcal{F})=\sum_aq^{-n_a}\otimes[C_a]
\end{equation}
where $q$ denotes the cardinality of $k$ and ${\cal F}(1)$ is the Tate twist (see Lemma \ref{lm}.1 for a more general statement). Since  $CC({\cal F})$ can be recovered from the right-hand side of (\ref{SSCCrel}), we may regard the theory of epsilon cycles as a refinement of the theory of characteristic cycles. 

We briefly review the proof of Theorem \ref{ayoo}. 
We make use of the machinery developed in \cite{Sai17} for the construction of 
characteristic cycles.  
Roughly speaking, this machinery shows that if an invariant varies ``continuously'' in any family of isolated characteristic points, then there 
exists a cycle whose intersection numbers recover the invariant;  
see Proposition~\ref{flcy} for the precise statement. 

To apply this machinery in our setting, we need 
a continuit result for local epsilon factors. In another paper \cite{CLE}, we prove that the local epsilon factor indeed varies ``continuously'': more precisely, for any family of isolated characteristic points parametrized by a $k$-scheme $S$ of finite type, there exists a continuous character 
\[\pi_1(S)^{ab}\to \Zl^\times
\]
whose values at geometric Frobenius elements equal  the corresponding 
local epsilon factors. For the precise statement, see \cite[Theorem 4.9.2]{CLE} or Theorem \ref{contiep}.  
By the finiteness theorem of Katz--Lang 
\cite{KL}, such a character gives a flat function (Definition \ref{flatfunc}) after taking modulo 
roots of unity, which ensures the existence of the cycle 
$\mathcal{E}(\mathcal{F})$ with the desired property. 

After establishing basic properties of epsilon cycles, we prove 
a pull-back formula for properly transversal morphism, using an argument 
similar to Beilinson's treatment of $CC({\cal F})$ in \cite[Section 7]{Sai17}.

\begin{thm}(Theorem \ref{prtr})
Let $k$ be a finite field with $q$ elements. 
Let $h\colon W\to X$ be a morphism of smooth $k$-schemes. 
Let $\mathcal{F}\in D^b_c(X,\Zl)$. 
Assume that $h$   is properly $SS(\mathcal{F})$-transversal 
(Definition \ref{trans}.2). Then we have an equality 
\begin{equation*}
\mathcal{E}(h^\ast\mathcal{F})=h^!(\mathcal{E}(\mathcal{F})(\frac{\dim X-\dim W}{2}))
\end{equation*}
in $\Theta\otimes Z_m(T^\ast W)$ ($m={\rm dim}W$). Here, for $r\in\mathbb{Q}$, we set 
 \[{\cal E}(\mathcal{F})(r):={\cal E}(\mathcal{F})+\sum_aq^{-rn_a}\otimes[C_a], 
 \]
  where the integers $n_a$ are determined by  \[CC({\cal F})=\sum_an_a[C_a]. \]
  Note that $q^{-rn_a}$ is well-defined in $\Theta$. 
  
  The symbol $h^!$ denotes the pull-back map  defined for cycles supported on the singular support; see Definition \ref{!plbk} for the precise definition. 
\end{thm}

Finally, we state and prove a product formula for global epsilon factors. 
\begin{thm}(Theorem \ref{thmEP})
Let $X$ be a smooth projective variety of dimension $n$  over a finite field $k$. 
Let $\mathcal{F}\in D^b_c(X,\Zl)$, and let ${\cal E}({\cal F})=\sum_a\xi_a\otimes[C_a]$ be  its epsilon cycle. Then 
\begin{equation*}
\prod_i\det({\rm Frob}_k,H^i(X_{\bar{k}},{\cal F}\otimes_{\Zl}\Ql))^{(-1)^i}
=(\mathcal{E}({\cal F}),
T^{\ast}_{X}X)_{T^{\ast}X}:=\prod_a\xi_a^{(C_a,T^\ast_XX)_{T^\ast X}}
\end{equation*}
as elements of $\Zl^\times\otimes_\mathbb{Z}\mathbb{Q}$. Here $(C_a,T^\ast_XX)_{T^\ast X}$ denotes the degree of the intersection 
\[[C_a]\cap[T^\ast_XX]\in{\rm CH}_0(T^\ast_XX)={\rm CH}_0(X). \]
\end{thm}

When $X$ is a curve, the above formula is equivalent to Laumon's product formula modulo roots of unity, which follows from the description of ${\cal E}({\cal F})$ in the curve case given in Example \ref{pfexample}. 

Note that we have the relation
\begin{equation*}
\prod_i\det({\rm Frob}_k,H^i(X_{\bar{k}},{\cal F}(1)\otimes_{\Zl}\Ql))^{(-1)^i}=q^{-\chi(X_{\bar{k}},{\cal F})}\cdot\prod_i\det({\rm Frob}_k,H^i(X_{\bar{k}},{\cal F}\otimes_{\Zl}\Ql))^{(-1)^i}. 
\end{equation*}
Therefore, Theorem~$1.3$ is compatible with (\ref{indexf}) and (\ref{SSCCrel}). 

As a consequence of Theorem $1.3$, we obtain a formula expressing the $p$-adic 
valuations of global epsilon factors as  products of those of 
local epsilon factors (Example \ref{p-adicval}). Fix a field embedding 
\[\iota\colon\Ql\to\overline{\mathbb{Q}_p}\] where $p$ is the characteristic of $k$. Then Theorem $1.3$ yields the formula 
\begin{equation*}
|\iota(\det({\rm Frob}_{k},R\Gamma(X_{\bar{k}},{\cal F})))|_p=\prod_a
|\iota(\xi_a)|_p^{\deg(C_a,T^\ast_XX)_{T^\ast X}}, 
\end{equation*}
where $|-|_p$ denotes the $p$-adic valuation on $\overline{\mathbb{Q}_p}$. 

In \cite{UYZ}, N. Umezaki, E. Yang, and Y. Zhao prove the twist formula of global epsilon 
factors \cite[Theorem 5.23.]{UYZ}. A weaker 
version modulo roots of unity can also be  deduced from 
the theorem above together with Lemma \ref{lm}.1. 

More generally, we construct the theory of epsilon cycles when the base field $k$ is the perfection of a finitely generated field over its prime field (the characteristic may be $0$). In this general setting, the coefficient group  $\Theta=\Zl^\times\otimes_{\mathbb{Z}}\mathbb{Q}$ is replaced by 
\[{\rm Hom}(G_k,\Zl^\times)\otimes_{\mathbb{Z}}\mathbb{Q}\]
 where $G_k$ is the absolute Galois group and ${\rm Hom}$ denotes continuous group homomorphisms. 
When $k$ is finite, we identify ${\rm Hom}(G_k,\Zl^\times)$ with $\Zl^\times$ by sending $\chi\mapsto \chi({\rm Frob}_k)$; in this way, the results described above are regarded as  special cases of the theorems obtained for arbitrary $k$.  

To formulate the Milnor-type formula when $k$ is the perfection of a finitely generated field, we need a theory of local epsilon factor in this general setting. When $k$ is perfect of positive characteristic, we use the theory of local epsilon factors over general perfect fields independently developed by Yasuda (\cite{Y3}, \cite{Y1}) and Guignard (\cite{geomep}). Yasuda defines local epsilon factors for representations with torsion coefficients, whereas Guignard defines them for $\Ql$-representations. Since we  work mainly with
 $\Zl$-representaions, which can be obtained as inverse limits of representations with torsion coefficients, we primarily follow Yasuda's results. 
 However,  Guignard's theory could also be used, because the two constructions give the same local epsilon factor for a $\Zl$-representation $V$: 
 \begin{equation*}
 \varepsilon_0(T,V,\omega)=\varepsilon_0(T,V\otimes_{\Zl}\Ql,\omega)
 \end{equation*}
 where $\varepsilon_0$ on the left is difined via  Yasuda's theory and the one on the right via   Guignard's. 
 
When the characteristic is $0$, we define local epsilon factors for 
 representations $V$ with unramified determinant using Jacobi sum characters constructed in \cite{Jacob}. For general $V$, 
we take a direct sum of copies of $V$ so that the determinant 
becomes unramified; thus in characteristic $0$ we can construct local epsilon factors only modulo 
roots of unity.

In \cite{geomephigh}, Guignard provides another method for computing global epsilon factors in higher dimensions, different in nature from the theory of 
epsilon cycles.  It would be interesting to clarify the relation between his results and ours.

We explain the construction of this paper. In Section $2$, we recall the theory of relative singular support given by Hu--Yang \cite{relsing}, which is a relative version of Beilinson's theory. We also recall the machinery that will be used to establish the existence of epsilon cycles. 
In Section $3$, we collect and prove necessary  results on local epsilon factors for the construction of epsilon cycles. Section $4$ is devoted to the construction of epsilon cycles and the proof of their basic properties. In Section $5$, we introduce the notion of epsilon classes, which is an analogue of characteristic classes defined in \cite{Sai17}. Using the epsilon classes together with the formalism of Radon transforms, we prove a product formula for $\det(R\Gamma(X_{\bar{k}},{\cal F}))$ in Theorem \ref{thmEP}. 
We also give an axiomatic description of 
epsilon cycles (Theorem \ref{epcygen}).

We give notation which we use throughout this paper. 
\begin{itemize}
\item We denote by $G_k$ the absolute Galois group of a field $k$. 
\item We denote by $\chi_\cyc\colon G_k\to\mathbb{Z}_\ell^\times$ 
the $\ell$-adic cyclotomic character. 
\item 
For a finite separable extension $k'/k$ of fields, we denote by ${\rm tr}_{k'/k}\colon G_k^{ab}
\to G_{k'}^{ab}$ the transfer morphism induced  
by the inclusion $G_{k'}\hookrightarrow G_k$. The determinant character of the 
induced representation ${\rm Ind}_{G_{k'}}^{G_k}
\mathbb{Q}_\ell$ of the trivial representation is denoted by $\delta_{k'/k}$. 
\item For a scheme $X$ and its point $x$, $k(x)$ is the 
residue field of $X$ at $x$. 
\item For a finite extension $x'/x$ of the spectra of fields, we denote by 
$\deg(x'/x)$ the degree of the extension. When $x={\rm Spec}(k)$ and 
$x'=\Spec(k')$, we also denote it by $\deg(k'/k)$. 
\item Let $x$ be a geometric point on a scheme $X$. We denote the strict henselization of $X$ at $x$ by $X_{(x)}$. On the other hand, we denote the henselization at a point $x\in X$ by $X_{(x)}$. 
More generally, for a finite separable extension $y$ of $x\in X$, we denote the henselization of $X$ at $y$ by $X_{(y)}$. 
\item For a scheme $X$ of finite type over $S$, we say that $X$ is of 
relative dimension $n$ if all fibers of $X\to S$ are equidimensional and of 
dimension $n$. 
\item We fix an algebraic closure $\Ql$ of $\mathbb{Q}_\ell$. Let $\mu$ be the 
group of roots of unity in $\Ql$. For a finite extension $E/\mathbb{Q}_\ell$, 
the ring of integers of $E$ is denoted by $\mathcal{O}_E$. 
\item For the $\ell$-adic formalism of sheaves of a noetherian topos $T$, we refer to 
\cite{adic}, which we review in the appendix. 
The derived category of constructible 
complexes of $\mathcal{O}_E$-sheaves on $T$ is denoted 
by  $D^b_c(T,\mathcal{O}_E)$. 
\end{itemize}

\tableofcontents

\section{Relative Singular Supports and Characteristic Cycles (\cite{Bei}, \cite{relsing}, \cite{Sai17})}
In this section, we recall the theory of singular supports over general bases and 
characteristic cycles. 

Let $S$ be a noetherian scheme. 
Let $\Lambda$ be a finite local ring  whose characteristic is invertible in $S$. $\Lambda_0$ will denote its residue field. 
Let $X$ be an $S$-scheme of finite type. 
A complex $K$ of \'etale sheaves of $\Lambda$-modules on $X$ is said to be {\it constructible} if each cohomology sheaf 
${\cal H}^i(K)$ is constructible and is zero except for finitely many $i$. We write $D_{\rm ctf}(X,\Lambda)$ for  the full subcategory of $D(X,\Lambda)$ consisting of constructible 
complexes of finite tor-dimension. 
\subsection{Relative singular support}
The singular supports are firstly defined by Beilinson in \cite{Bei} for smooth varieties over a field. Later Hu--Yang \cite{relsing} generalize the results to smooth schemes over a general noetherian scheme. In what follows, we briefly recall their results; 
for detailed discussions,  we refer to 
\cite{relsing}. 

Let $X$ be a smooth scheme over $S$. We write $T^\ast (X/S)$ for the 
cotangent bundle of $X$ relative to $S$. 
 $T^\ast_X(X/S)$ denotes the $0$-section of $T^\ast (X/S)$. When $S$ is the spectrum of a 
field, we often simply write  $T^\ast X$, $T^\ast_XX$ for them. 
For a 
morphism $x\to X$ where $x$ is the spectrum of a field,  $T^\ast_x(X/S)$ denotes  the base change $T^\ast(X/S)\times_Xx$. 
We say that a closed subset $C$ of 
$T^\ast (X/S)$ is conical if $C$ is stable under the action of $\mathbb{G}_{\rm m}$ on $T^\ast(X/S)$. 

First we recall the notions of $C$-transversality. 
\begin{df}\label{trans}
Let $X$ be a smooth scheme over $S$ and $C$ be a closed conical subset of 
$T^\ast (X/S)$. 
\begin{enumerate}
\item(\cite[2.4]{relsing}) We say that an $S$-morphism $h\colon W\to X$ from a smooth $S$-scheme $W$ 
is {\rm $C$-transversal} if, 
for every geometric point $w\to W$, 
non-zero elements 
of $C\times_Xw\subset T^\ast_{h(w)}(X/S):=T^\ast(X/S)\times_Xw$ map to non-zero elements of $T^\ast_w(W/S)$ via $dh_w\colon T^\ast_{h(w)}(X/S)\to T^\ast_w(W/S)$. 
\item Assume that $X$ and $C$ is of relative dimension $n$. Let $W$ be a 
smooth scheme over $S$ of relative dimension $m$. 
We say that an $S$-morphism 
$h\colon W\to X$ is {\rm properly $C$-transversal} if $h$ is $C$-transversal in the sense of $1.$ and 
$W\times_XC$ is of relative dimension $m$. 
\item(\cite[2.4]{relsing}) We say that an $S$-morphism $f\colon X\to Y$ to a smooth $S$-scheme $Y$ is 
{\rm $C$-transversal} if, for every geometric point $x\to X$, no non-zero elements 
of $T^\ast_{f(x)}(Y/S):=T^\ast(Y/S)\times_Yx$ map into $C\times_Xx$ via $df_x\colon T^\ast_{f(x)}(Y/S)\to T^\ast_{x}(X/S)$. 
\end{enumerate}
\end{df}
\begin{lm}(cf. \cite[1.2]{Bei})\label{ctrfin}
Let $h\colon W\to X$ be a morphism of smooth $S$-schemes and $C$ be a closed 
conical subset of $T^\ast (X/S)$. If $h$ is $C$-transversal, then the map 
$dh\colon C\times_X W\to T^\ast (W/S)$ is finite. 
\end{lm}
\proof{Replacing $W$ with an open cover, we may assume that $W$ is affine. Then, the same argument given in Lemma (ii) in \cite[1.2]{Bei} also works in our setting. 
\qed
}
\begin{df}\label{trtr}
Let $X$ and $C$ be as in Definition \ref{trans}. Let $W$ and $Y$ be 
smooth $S$-schemes. 
\begin{enumerate}
\item(\cite[2.4]{relsing}) Let $h\colon W\to X$ be a morphism of $S$-schemes. When $h$ is 
$C$-transversal, we define $h^\circ C$ to be the image of 
$dh\colon C\times_XW\to T^\ast (W/S)$. 
This is a closed conical subset of $T^\ast(W/S)$ by Lemma \ref{ctrfin}. 
\item(\cite[2.4]{relsing}) Let $f\colon X\to Y$ be a morphism of $S$-schemes. Assume that 
$f$ is proper on the base of $C$, i.e., on $C\cap T^\ast_X(X/S)$. We define 
$f_\circ C$ to be ${\rm pr}_1(df^{-1}(C))$ where $df\colon T^\ast (Y/S)\times_YX\to T^\ast (X/S)$ and ${\rm pr}_1\colon T^\ast (Y/S)\times_YX\to T^\ast 
(Y/S)$. This is a closed conical subset of $T^\ast(Y/S)$. 
\item Let $(h,f)$ be a pair of $S$-morphisms 
\begin{equation*}
X\xleftarrow{h}W\xrightarrow{f}Y. 
\end{equation*}
We say that $(h,f)$ is {\rm $C$-transversal} if 
$h$ is $C$-transversal and $f$ is $h^\circ C$-transversal (where $h^\circ C$ is defined in $1.$). 
\end{enumerate}
\end{df}

\begin{df}
Let $X$ and $C$ be as in Definition \ref{trans}. Let $K\in D_{\rm ctf}(X,\Lambda)$. 
We say that $K$ is micro-supported on $C$ if, for any $C$-transversal pair $(h,f)$ as in 
Definition \ref{trtr}.3, $f$ is locally acyclic relatively to 
$h^\ast K$. 
\end{df}

Let $\Lambda\to\Lambda'$ be a local ring homomorphism between finite local rings. 
The following lemma shows that local acyclicity for $K\in D_{\rm ctf}(X,\Lambda)$ is equivalent to that for $K\otimes^L_\Lambda\Lambda'$. 
\begin{lm}\label{locacyadic}
Let $f\colon X\to Y$ be a morphism of schemes of finite type. Suppose that the characteristic of $\Lambda$ is invertible in $Y$. 
Let $K\in D_{\rm ctf}(X,\Lambda)$. Then, $f$ is (resp. universally) locally acyclic relatively to 
$K$ if and only if so is $f$ relatively to $K\otimes^L_\Lambda\Lambda'$. 
\end{lm}
\proof{
Let $\Lambda_0$ be the residue field of $\Lambda$. 
Let $x$ be a geometric point of $X$ and $y$ be a geometric point of $Y$ 
which is a generalization of $f(x)$. Since the functor 
$\Gamma(X_{(x)}\times_{Y_{(f(x))}}y,-)$ is of finite cohomological dimension, 
we have $R\Gamma(X_{(x)}\times_{Y_{(f(x))}}y,K)\otimes^L_\Lambda
\Lambda_0\cong R\Gamma(X_{(x)}\times_{Y_{(f(x))}}y,K\otimes^L_\Lambda
\Lambda_0)$. On the other hand, 
$\Lambda$ has a filtration of finite length whose graded quotients are finite free $\Lambda_0$-modules. This consideration implies that local acyclicity for $K$ is equivalent to that for $K\otimes_\Lambda^L\Lambda_0$. Therefore, we may assume that $\Lambda$ and $\Lambda'$ are finite fields. Then, the assertion is clear since $\Lambda'$ is a 
finite free $\Lambda$-module. 
\qed
}

\begin{lm}\label{lmtr}
Let $X$ and $C$ be as in Definition \ref{trans}. 
Let $K\in D_{\rm ctf}(X,\Lambda)$. Assume that $K$ is micro-supported on 
$C$. 
\begin{enumerate}
\item(\cite[Lemma 4.7(ii)]{relsing},\cite[Lemma 2.1.(ii)]{Bei})
Suppose that $(h,f)$ is a $C$-transversal pair (Definition \ref{trtr}.3). Then, 
$f$ is universally locally acyclic relatively to $h^\ast K$. 
\item Let $h\colon X'\to X$ be an $S$-morphism where $X'$ is a smooth $S$-scheme. 
If $h$ is $C$-transversal, then $h^\ast K$ is micro-supported on 
$h^\circ C$. 
\end{enumerate}
\end{lm}
\proof{
$1.$ The case where $\Lambda$ is a field is treated in \cite[Lemma 4.7(ii)]{relsing}. The general case then follows from this and Lemma \ref{locacyadic}. 

2. Let $X'\xleftarrow{h'}W\xrightarrow{f}Y$ be an 
$h^\circ C$-transversal pair. Then, the pair $(h\circ h',f)$ is 
$C$-transversal. Thus, $f$ is locally acyclic relatively to $h'^\ast 
h^\ast K$. 
\qed
}

We recall the theorem on the existence of relative singular supports. 
\begin{thm}\label{relSS}
Let $X$ be a smooth $S$-scheme of finite type. Let $K$ be a complex in 
$D_{\rm ctf}(X,\Lambda)$. After replacing $S$ with a  Zariski open dense subscheme, 
the following hold. 
\begin{enumerate}
\item(\cite[Theorem 5.2(2)]{relsing}) There exists a smallest closed conical subset $C$ of $T^\ast(X/S)$ on which 
$K$ is micro-supported. In this case, the smallest $C$ is called {\rm the relative singular support} and denoted 
by $SS(K,X/S)$. 
\item(\cite[Theorem 5.3]{relsing}) For a morphism $s\to S$ from the spectrum of a field, we have
\begin{equation*}
SS(K|_{X_s})=SS(K,X/S)\times_Ss. 
\end{equation*}
Here $SS(K|_{X_s})=SS(K|_{X_s},X_s/s)$ is the singular support, defined in \cite{Bei}, in the case over a field. 
\item Let $\Lambda\to\Lambda'$ be a local morphism of finite local rings. 
Then, $SS(K,X/S)$ exists if and only if $SS(K\otimes^L_\Lambda\Lambda',X/S)$ exists. If these conditions are satisfied, then we have
\begin{equation*}
SS(K,X/S)=SS(K\otimes^L_\Lambda\Lambda',X/S). 
\end{equation*}
\end{enumerate}
\end{thm}
\proof{ When $\Lambda$ is a field, the assertions $1.$ and $2.$ are proved in loc.~cit. Then, the general case follows from this case and Lemma \ref{locacyadic}. 
The assertion $3.$ also follows from Lemma \ref{locacyadic}.  
\qed
}
\begin{rmk}\label{rmkss}
If the relative singular support $SS(K,X/S)$ exists, then 
the structure morphism $X\to S$ is universally locally acyclic relatively 
to $K$, since $X\xleftarrow{{\rm id}}X\to S$ is $SS(K,X/S)$-transversal, as is explained in \cite[Proposition 4.5]{relsing}. 
If further $X$ is projective over $S$, then the existence of $SS(K,X/S)$ is 
equivalent to the universal local acyclicity of $X\to S$ relative to $K$ 
\cite[Theorem 5.2]{relsing}. 
\end{rmk}
We give some examples of singular supports. 
\begin{pr}\label{example}
Suppose that $S={\rm Spec}(k)$ is the spectrum of a field $k$. We use the abbreviation $SS(-)$ for $SS(-,X/S)$. 
\begin{enumerate}
\item(\cite[Lemma 4.3.2]{Sai17}) Let $X$ be a smooth curve over $k$. Let $K\in D_{\rm ctf}(X,\Lambda)$. 
Then we have 
\begin{equation*}
SS(K)\subset T^\ast_XX\cup\bigcup_xT^\ast_xX,
\end{equation*}
where $x$ runs through the closed points around which $K$ is not locally constant. 
The equality holds if and only if the generic fiber of $K$ is not acyclic. 
\item(\cite[Theorem 2.2.3.]{Saiex}) Let $X_1$ and $X_2$ be smooth schemes over $k$. Let $K_i\in D_{\rm ctf}(X_i,
\Lambda)$ for each $i=1,2$. Then, we have 
\begin{equation*}
SS(K_1\boxtimes^LK_2)=SS(K_1)\times_k SS(K_2) \subset 
T^\ast X_1\times_k T^\ast X_2\cong T^\ast(X_1\times_k X_2).
\end{equation*}
\end{enumerate}
\end{pr}
\proof{
The assertions are proved in the references given above. The assertion $2.$ is proved by applying the Thom--Sebastiani theorem in the \'etale setting \cite{Ill}. 
\qed}

\subsection{Flat function and characteristic cycle}

In this subsection, we recall key ingredients for the construction of characteristic cycles given in \cite{Sai17} as they play 
an essential role also in constructing epsilon cycles. 

 Let $k$ be a perfect 
field. We generalize some notions and results
on isolated $C$-characteristic points given in \cite{Sai17}, where  $k$ is assumed algebraically closed, to a non-algebraically closed $k$. 
For a scheme $X$, we write $|X|$ for the set of closed points of $X$. 
 We fix an abelian group $A$. 
\begin{df}(cf. \cite[Definition 5.5]{Sai17})\label{flatfunc}
Let $Z$ be a scheme locally of finite type over $k$. Let 
$\varphi\colon |Z|\to A$ be a function. 
\begin{enumerate}
\item For a morphism of schemes of finite type $f\colon Z'\to Z$, define a function $f^\ast\varphi\colon 
|Z'|\to A$ by $f^\ast\varphi(z'):=\deg(z'/f(z'))\varphi(f(z'))$. 
If no confusions occur, we simply write $\varphi|_{Z'}=f^\ast\varphi$. 
\item The function $\varphi$ is said to be {\rm constant} if 
there exists a function $\psi\colon|{\rm Spec}(k)|\to A$  with $\psi|_Z=\varphi$. 
\item Let $g\colon Z\to S$ be a quasi-finite morphism of schemes locally of finite type 
over $k$. We say that $\varphi$ is {\rm flat} over $S$ if the following condition holds:

For every closed point $z\in Z$, There exists a commutative diagram 
\begin{equation}\label{flat0}
\xymatrix{
U\ar@{^{(}-_>}[r]\ar[rd]_{\tilde{g}}&V\times_SZ\ar[r]\ar[d]&Z\ar[d]^g\\
&V\ar[r]&S
}
\end{equation}
of $k$-schemes with the following properties:
\begin{enumerate}
\item $V\to S$ is \'etale and there exists a closed point $v\in V$ whose image 
in $S$ is equal to $g(z)$. The map $v\to g(z)$ induces  an isomorphism on the residue fields. 
\item U is an open neighborhood of $(v,z)$ in $V\times_SZ$. 
\item U is finite over $V$. The fiber of $\tilde{g}$ over $v$ only consists of 
$(v,z)$. 
\item The function $\tilde{g}_\ast(\varphi|_U)\colon |V|\to A$ defined by 
\begin{equation*}
\tilde{g}_\ast(\varphi|_U)(x)=\sum_{y\in\tilde{g}^{-1}(x)}\varphi|_U(y)
\end{equation*}
 is 
constant in the sense of 2. 
\end{enumerate}
\end{enumerate}
\end{df}

Let $X$ be a smooth $k$-scheme and $C$ be a closed conical subset of $T^\ast X$. Let $f\colon X\to Y$ be a 
$k$-morphism to a smooth curve $Y$ over $k$. Let $x\in X$ be a closed point. 
\begin{df}\label{isolchar}
Let the notation be as above. 
\begin{enumerate}
 \item We say that $x$ is {\rm an at most isolated $C$-characteristic point} 
 of $f$ if 
there exists an open neighborhood $U$ of $x$ such that 
the restriction $f|_{U\setminus \{x\}}$ is $C$-transversal. 
\item Suppose that $X$ is purely of dimension $n$ and that every irreducible 
component $C_a$ of $C$ is of dimension $n$. Let $\alpha=
\sum_a\beta_a\otimes [C_a]$ be a cycle with A-coefficient that is supported on $C$. Namely, $\alpha$ is a linear combination of the irreducible components of 
$C$ with $\beta_a\in A$. 
Assume that $x\in X$ is an at most isolated $C$-characteristic point of $f$. We define the intersection number 
$(\alpha,df)_{T^\ast X,x}\in A$, or simply written as $(\alpha,df)_x$, by setting 
\begin{equation*}
(\alpha,df)_{T^\ast X,x}=\sum_a(C_a,df)_{T^\ast X,x}\cdot \beta_a,
\end{equation*}
where $(C_a,df)_{T^\ast X,x}$ is the intersection number, supported on the fiber 
of $x$, of $C_a$ and the section 
$f^\ast\omega$ of $T^\ast X$ defined by the pull-back of a basis $\omega$ of 
$T^\ast Y$ on a neighborhood of $f(x)\in Y$. 
\end{enumerate}
\end{df}

\begin{df}\label{fcnchar}
Let $X$ and $C$ be as above. 
\begin{enumerate}
\item $\varphi$ is said to be {\rm an $A$-valued function on isolated $C$-characteristic points} if, for every diagram 
\begin{equation}\label{char}
\xymatrix{
U\ar[r]^f\ar[d]_j&Y\\X
}
\end{equation}
of $k$-schemes and a closed point $u\in |U|$ such that $Y$ is a smooth curve over $k$, $U$ is \'etale 
over $X$, and $u$ is an at most isolated $C$-characteristic point 
of $f$, an element $\varphi(f,u)\in A$ is given. Further this assignment  
should satisfy the following 
conditions: 
\begin{enumerate}
\item $\varphi(f,u)$ is $0$ when $u$ is not an isolated $C$-characteristic point. 
\item For every commutative diagram 
\begin{equation}\label{char2}
\xymatrix{
U'\ar[r]^{f'}\ar[d]&Y'\ar[d]\\
U\ar[r]^f\ar[d]&Y\\X
}
\end{equation}
of $k$-schemes such that the vertical arrows are \'etale and $Y$, $Y'$ are smooth curves over $k$ 
and  an at most isolated $C\times_XU'$-characteristic point $u'\in U'$ 
of $f'$, we have $\varphi(f',u')=\deg(u'/u)\cdot\varphi(f,u)$ 
where $u$ is the image of $u'$ by $U'\to U$. 
\end{enumerate}
\item Let $\varphi$ be an $A$-valued function on isolated $C$-characteristic 
points. $\varphi$ is said to be {\rm flat} if, for every commutative diagram 
\begin{equation}\label{c}
\xymatrix{
Z\ar@{^{(}-_>}[r]&U\ar[rr]^f\ar[rd]\ar[d]_{{\rm pr}_1}&&Y\ar[ld]\\
&X&S&
}
\end{equation}
of $k$-schemes such that 
\begin{itemize}
\item S is a smooth scheme over $k$, 
\item $Y\to S$ is a relative smooth curve, 
\item the map $U\to X\times_kS$ is \'etale, 
\item $Z$ is a closed subscheme of $U$ quasi-finite over $S$, and 
\item the pair $({\rm pr}_1 ,f)$ is $C$-transversal 
(in the sense of Definition \ref{trtr}.3) outside $Z$, 
\end{itemize}
the function $\varphi_f\colon |Z|\to A$ defined by $\varphi_f(z)=\varphi(f_s,z)$, 
where $s$ is the image of $z$ by $Z\to S$ and $f_s\colon U_s\to Y_s$ is the base 
change of $f$ by $s\to S$,
 is flat over $S$ in the sense of Definition \ref{flatfunc}.3. 
\end{enumerate}
\end{df}

\begin{pr}(\cite[Proposition 5.8]{Sai17})\label{flcy}
Assume that $A$ is uniquely divisible (i.e., the canonical map $A\to A\otimes_\mathbb{Z}\mathbb{Q}$ is an isomorphism). 
Let $X$ be a smooth scheme purely of dimension $n$ over $k$. Let $C$ 
be a closed conical subset of $T^\ast X$. Assume that every irreducible 
component $C_a$ of $C$ is of dimension $n$. Let $\varphi$ be an 
$A$-valued function on isolated $C$-characteristic points. 
The following conditions are equivalent. 
\begin{enumerate}
\item $\varphi$ is flat.
\item There exists a cycle $\alpha=\sum_{a}\beta_a\otimes [C_a]\in 
A\otimes_\mathbb{Z}Z_n(T^\ast X)$ (where $Z_n$ denotes the group of $n$-cycles) with A-coefficient and supported on 
$C$ such that 
\begin{equation}\label{b}
\varphi(f,u)=\deg(u/k)(j^\ast\alpha,df)_{T^\ast U,u}
\end{equation}
holds for every diagram (\ref{char}) and every at most isolated $C$-characteristic point $u\in U$ of $f$. 
\end{enumerate}
Further, if these conditions hold, the cycle 
$\alpha$ in $2$ is unique. 
\end{pr}
\proof{
Since the proof is completely similar as \cite[Proposition 5.8]{Sai17} and we only use the implication 
$1\Rightarrow 2$ in the sequel, we sketch the proof of $1\Rightarrow 2$

First we consider the case when $k$ is algebraically closed. By the 
similar argument in \cite[Proposition 5.8]{Sai17}, we find a unique cycle 
$\alpha_X\in
A\otimes Z_n(T^\ast X)$ satisfying $\varphi(f,u)=
(\alpha_X,df)_{T^\ast X,u}$ for every diagram (\ref{char}) with 
$U\to X$ being an open immersion. Let $j\colon W\to X$ be an \'etale morphism. 
Restricting $\varphi$ to $W$, we have an $A$-valued function on 
isolated $j^\ast C$-characteristic points. Since this is also flat, we find a 
cycle $\alpha_W\in A\otimes Z_n(T^\ast W)$ satisfying 
$\varphi(f,u)=(\alpha_W,df)_{T^\ast W,u}$ for every diagram (\ref{char}) 
replaced $X$ by $W$. We need to show the equality $\alpha_W=j^\ast
\alpha_X$, which is a consequence of \cite[Proposition 5.8.2]{Sai17}. 

Next we consider the general case. Take an algebraic closure $\bar{k}$ of 
$k$. We put the symbol 
 $\bar{k}$ to mean the base change by $k\to
\bar{k}$. 
From $\varphi$, we define an $A$-valued function $\varphi_{\bar{k}}$ on 
isolated $C_{\bar{k}}$-characteristic points as follows. 
Let 
\begin{equation}\label{a}
\xymatrix{
U\ar[r]^f\ar[d]&Y\\X_{\bar{k}}
}
\end{equation}
be a diagram as in (\ref{char}), and $u\in U$ be an 
at most isolated $C_{\bar{k}}$-characteristic point of $f$. 
We assume that $U$ and $Y$ are quasi-compact. 
Take a finite subextension 
$k'/k$ in $\bar{k}$ such that there exists a diagram of 
$k'$-schemes 
\begin{equation*}
\xymatrix{
U'\ar[r]^{f'}\ar[d]&Y'\\X_{k'}
}
\end{equation*}
whose base change by $k'\to\bar{k}$ is isomorphic to (\ref{a}). Let $u'\in
U'$ be the image of $u$. Put $\varphi_{\bar{k}}(f,u):=\frac{1}{\deg(u'/k)}
\varphi(f',u')$, where $\varphi$ in  the right hand side is defined for  the diagram 
\begin{equation*}
\xymatrix{
U'\ar[r]^{f'}\ar[d]&Y'\\X_{k'}\ar[d]&\\X.
}
\end{equation*}
This is independent of the choice of $(k',f')$ and defines an 
$A$-valued function on isolated $C_{\bar{k}}$-characteristic points. Since 
$\varphi_{\bar{k}}$ is flat, we find a cycle $\alpha_{\bar{k}}\in A\otimes_{\mathbb{Z}}Z_n(T^\ast X_{\bar{k}})$ satisfying 
(\ref{b}). From the construction of $\varphi_{\bar{k}}$, $\alpha_{\bar{k}}$ 
is contained in the Galois-fixed part $(A\otimes_{\mathbb{Z}}Z_n(T^\ast X_{\bar{k}}))^{
{\rm Gal}(\bar{k}/k)}$. By \'etale descent, we get a cycle $\alpha$ that satisfies the condition. 
\qed
}

Let $X\xleftarrow{j} U\xrightarrow{f} Y$ be as in (\ref{char}). 
Take a closed point $u\in U$. 
Let $\Lambda$ be a finite local ring as before and let $K\in D_{\rm ctf}(X,\Lambda)$. 
Suppose that $u$ is an at most isolated $j^\ast SS(K)$-characteristic point of $f$. Let $U_{(u)}$ and $Y_{(f(u))}$ denote the henselizations of $U$ and $Y$ at  $u$ and $f(u)$ respectively. Let $Y_{(u)}$ be the 
unramified extension of $Y_{(f(u))}$ that corresponds to the separable 
extension $k(u)/k(f(u))$ of residue fields. 
\begin{df}\label{vancyc}
Let the notation be as above. Let $f_{(u)}\colon U_{(u)}\to Y_{(u)}$ be the induced 
map. 
We write $R\Phi_f(K)_u$ for the vanishing cycles complex of 
$K|_{U_{(u)}}$ with respect to the map $f_{(u)}$ and call it {\rm the vanishing cycles complex of $K$ supported at $u$}. This is a constructible complex 
on $U_{(u)}\times_{Y_{(u)}}\bar{u}$ supported at $\bar{u}$ endowed with an equivariant 
action of the absolute Galois group of $\eta_u$ where 
$\eta_u$ is the generic point of $Y_{(u)}$ and 
 $\bar{u}$ is a geometric point over $u$ (\cite{SGA7-2}). Since it is supported at one point, we usually regard $R\Phi_f(K)_u$ as a complex of representations of the absolute Galois group of $\eta_u$ on $\Lambda$-modules, or as 
 an object of $D_{\rm ctf}(\eta_u,\Lambda)$. 
\end{df}
For an object $M\in D_{\rm ctf}(\eta_u,\Lambda)$, define the total dimension 
$\dt M$ of $M$ to be $\dt M:=\rk M+\sw M$. 

Finally we recall the definition of characteristic cycles.  
\begin{thm}(\cite[Theorems 5.9, 5.18]{Sai17})\label{CC}
Let $X$ be a smooth scheme over $k$ and $K\in D_{\rm ctf}(X,\Lambda)$. Let $C$ be a closed conical subset of $T^\ast X$ on which 
$K$ is micro-supported. Assume that each irreducible component of $X$ and 
that of $C$ is of dimension $n$. 
\begin{enumerate}
\item 
There exists a cycle 
$CC(K)$ in $\mathbb{Q}\otimes Z_n(T^\ast X)$, supported on $C$, satisfying  
the following property: 

For every diagram as (\ref{char}) and an at most isolated $C$-characteristic 
point $u\in U$ of $f$, we have
\begin{equation*}
-\dt R\Phi_f(K)_u=(CC(K),df)_u.
\end{equation*}
Moreover, $CC(K)$ is unique and independent of a  choice of $C$ on 
which $K$ is micro-supported. $CC(K)$ has coefficients in $\mathbb{Z}$. 
\item Let $\Lambda\to\Lambda'$ be a local ring homomorphism between finite local rings. Then, we have 
\begin{equation*}
CC(K)=CC(K\otimes_\Lambda^L\Lambda'). 
\end{equation*}
\end{enumerate}
\end{thm}
\proof{
$1$. The existence and the uniqueness are direct consequences of Proposition \ref{flcy}, once one knows that the 
$\mathbb{Q}$-valued function on isolated $C$-characteristic points defined 
by $\varphi(f,u)=-\deg(u/k)\cdot\dt R\Phi_f(K)_u$ is flat. 
This flatness is proved in \cite[Proposition 2.16]{Sai17}. 
The integrality is proved in \cite[Theorem 5.18]{Sai17}. 

$2$.  Let $X\xleftarrow{j}U\xrightarrow{f}Y$ be as (\ref{char}). The assertion 
follows from $R\Phi_f(K)\otimes^L_\Lambda\Lambda'\cong 
R\Phi_f(K\otimes^L_\Lambda\Lambda')$ and the equality 
$\dt M=\dt( M\otimes^L_\Lambda\Lambda')$ for $M\in D_{\rm ctf}(\eta_u,\Lambda)$. 
\qed
}

Let $\ell$ be a prime number which is invertible in schemes under consideration. Take an algebraic closure 
$\Ql$ of $\mathbb{Q}_\ell$. 
We give definitions of the local acyclicity, 
singular supports, and characteristic cycles for $\Zl$-sheaves. 
For the $\ell$-adic formalism, we follow Ekedahl's method, which we recall in Section \ref{adical}. In the sequel, 
$E$ is a finite subextension in $\Ql/\mathbb{Q}_\ell$ and 
${\cal O}_E$ denotes its ring of integers. 
\begin{df}\label{LAadic}
Let $f\colon X\to Y$ be a morphism of noetherian schemes of finite type.  Let $\cal O$ 
be either of $\mathcal{O}_E$ or $\Zl$. 
Let $\mathcal{F}$ be an element in $D^b_c(X,{\cal O})$. 
\begin{enumerate}
\item When ${\cal O}={\cal O}_E$, we say that 
$f$ is (resp. universally) locally acyclic relatively to 
$\mathcal{F}$ if, for some (hence all) $n\geq0$, 
$f$ is (resp. universally) locally acyclic relatively to 
$\mathcal{F}\otimes^L_{\cal O}{\cal O}/\ell^{n+1}$. 
\item 
When ${\cal O}=\Zl$, we take a finite subextension $E$ and ${\cal F}_{E}\in D^b_c(X,{\cal O}_E)$ with ${\cal F}_{E}\otimes_{{\cal O}_E}\Zl\cong {\cal F}$. 
 we say that 
$f$ is (resp. universally) locally acyclic relatively to 
$\mathcal{F}$ if so is  
$f$ relatively to 
$\mathcal{F}_{E}$ in the sense of $1$. 
\end{enumerate}
{\rm The independences of auxiliary choices made above follow from Lemma \ref{locacyadic}. }

\end{df}

\begin{df}\label{SSCCadic}
Let $X$ be a smooth scheme of finite type over a noetherian scheme $S$. 
Let ${\cal O}$ be either of $\mathcal{O}_E$ or $\Zl$. 
Let $\mathcal{F}\in D^b_c(X,{\cal O})$. 
\begin{enumerate}
\item If ${\cal O}=\mathcal{O}_E$, then 
we define 
$SS(\mathcal{F},X/S):=SS(\mathcal{F}\otimes^L_{\mathcal{O}_E}
\mathcal{O}_E/\ell^{n+1},X/S)$ and, when $S$ is the spectrum of a perfect field, 
$CC(\mathcal{F}):=CC(\mathcal{F}\otimes^L_{\mathcal{O}_E}\mathcal{O}_E/\ell^{n+1})$ for some (hence all) $n\geq0$.  

\item When ${\cal O}=\Zl$, take a finite subextension $E$ and 
$\mathcal{F}_{E}\in D^b_c(X,\mathcal{O}_E)$ with $\mathcal{F}_{E}\otimes^L_{\mathcal{O}_E}\Zl\cong
\mathcal{F}$. We define $SS(\mathcal{F},X/S):=SS(\mathcal{F}_{E},X/S)$ and, when $S$ is the spectrum of a perfect field, 
$CC(\mathcal{F}):=CC(\mathcal{F}_{E})$ in the sense of $1$. 
\end{enumerate}
{\rm These are independent of auxiliary choices by 
Theorems \ref{relSS}.3 and \ref{CC}.2. }
\end{df}

\begin{rmk}
According to \cite[5.3.]{UYZ}, the 
characteristic cycle can be defined for an object in $D^b_c(X,\Ql)$. 
\end{rmk}

\subsection{Reminder on good pencils}
We recall the theory of the universal hyperplane sections 
(\cite[Section 3.2]{Sai17}) and the notion of good pencils (\cite{SY}). Let $X$ be a quasi-projective smooth scheme over a field $k$. Let $\mathcal{L}$ be an ample invertible $\mathcal{O}_X$-module. 
Let $V$ be a $k$-vector space of finite dimension and $V\to\Gamma(X,\mathcal{L})$ be a $k$-linear mapping that induces a surjection 
$V\otimes_k\mathcal{O}_X\to\mathcal{L}$. Suppose that the 
induced map $h\colon X\to \mathbb{P}=
\Proj(V^\vee)$ is an immersion. Here we use a contra-Grothendieck notation for a 
projective space, i.e., $\Proj(V^\vee)$ parametrizes sub line bundles of $V^\vee$. 
Let $\Proj^\vee:=\Proj(V)$ be the dual projective space. The universal hyperplane 
$Q\subset\Proj\times\Proj^\vee$ parametrizes pairs $(x,H)$ consisting of 
points $x\in \Proj$ and hyperplanes $H\in\Proj^\vee$ which contain 
$x$. Since the kernel of the tautological surjection 
$V\otimes_k\mathcal{O}_\Proj(-1)\to\mathcal{O}_\Proj$ is canonically 
isomorphic to the cotangent bundle $\Omega^1_\Proj$, $Q$ is 
identified with the projective space bundle $\Proj(T^\ast\Proj)$. 
The composition $T^\ast_Q(\Proj\times\Proj^\vee)\to
Q\times_{\Proj\times\Proj^\vee}T^\ast(\Proj\times\Proj^\vee)\to
Q\times_\Proj T^\ast\Proj$ is the universal sub line bundle on $Q$.

Consider the following commutative diagram 
\begin{equation}\label{radon}
\xymatrix{X\ar[d]^-h
&X\times_\Proj Q\ar[l]_-{p}
\ar[r]^-{p^\vee}\ar[d]
&\Proj^\vee\\
\Proj&Q\ar[l]^-{{\bm p}}\ar[ru]_-{{{\bm p}}^\vee}
}
\end{equation}
where the square is cartesian. 
We have $X\times_\Proj Q=\Proj(X\times_\Proj T^\ast\Proj)$. 

Let $L\subset \Proj^\vee$ be a line in $\Proj^\vee$. We have a 
commutative diagram 
\begin{equation}\label{Lefpen}
\xymatrix{
&X_L\ar[r]^f\ar[d]\ar[ld]_\pi&L\ar[d]\\X&X\times_\Proj Q
\ar[r]_{\ \ p^\vee}\ar[l]^{p\ }&\Proj^\vee
}
\end{equation}
where the right square is cartesian. Write $A_L$ for the axis of $L$ 
in $\Proj$. This is a subspace of $\Proj$ of codimension $2$. The 
$\Proj$-scheme $\Proj_L=Q\times_{\Proj^\vee}L$ is 
isomorphic to the blow-up of $\Proj$ along $A_L$. Hence, 
if $X$ and $A_L$ meet transversally, then $X_L$ is the blow-up of $X$ along 
the smooth subvariety $X\cap A_L$.

\begin{df}\label{goodpen}
Let $X\subset\Proj$ be a closed smooth subvariety purely of dimension $n$ over $k$. 
Let $C$ be a closed conical subset of $T^\ast X$ whose irreducible 
components are of dimension $n$. We call the pair $(\pi,f)$ as in (\ref{Lefpen}) 
a good pencil with respect to $C$ if the following conditions hold. 
\begin{enumerate}
\item $X$ and $A_L$ meet transversally (then $X_L$ is smooth). 
\item The morphism $\pi$ is properly $C$-transversal. 
\item The morphism $f$ has at most isolated $\pi^\circ C$-characteristic points. 
\item For every closed point $y\in L$, there exists at most one 
$\pi^\circ C$-characteristic point on the fiber $f^{-1}(y)$. 
\item No isolated characteristic points of $f$ are contained in the exceptional 
locus of $\pi$. 
\item For every irreducible component $C_a$ of $C$, there is an isolated 
$\pi^\circ C$-characteristic point $x\in X_L$ at which $df$ meets $C_a$, but does not meet any other component. For every isolated 
$\pi^\circ C$-characteristic point $x$, $df$ meets exactly one irreducible component at $x$. 
\item For every isolated $\pi^\circ C$-characteristic point $x\in X_L$ of $f$, 
the morphism $x\to f(x)$ of the spectra of fields is purely inseparable. 
\end{enumerate}
\end{df}

The existence of good pencils after a finite field extension is proved in \cite[Lemma 4.2.7.]{UYZ} by using 
\cite[Lemma 2.3.]{SY}. 
\begin{lm}\label{exgoodpen}(\cite[Lemma 4.2.7.]{UYZ}, \cite[Lemma 2.3.]{SY}) 
Let $X$ and $C$ be as in Definition \ref{goodpen}. 
Suppose that $k$ is an infinite field. 
Let ${\rm Gr}(1,\Proj^\vee)$ be the Grassmanian variety parametrizing lines in 
$\Proj^\vee$. Then, after composing $X\hookrightarrow\Proj$ and 
the Veronese embedding $\Proj\hookrightarrow\Proj'$ of $\deg\geq3$ if necessary, 
there exists a dense 
open subset $U\subset{\rm Gr}(1,\Proj^\vee)$ such that, 
for every 
$k$-rational 
point $L\in U(k)$, the pair $(\pi,f)$ in (\ref{Lefpen}) is a good pencil. 
\end{lm}

\section{Local Epsilon Factors (cf. \cite{Del}, \cite{Y1}, \cite{Y3})}\label{LEF}

In the preliminary subsection \ref{3.1}, we review 
Yasuda's generalization of the theory of local epsilon factors for henselian traits of equal-characteristic. In the last subsection \ref{redpos}, we give a way to reduce various computations for $\ell$-adic sheaves on varieties of characteristic $0$ to the case of positive characteristic. In $3.4$, we compute the local epsilon factors of the convolutions of vanishing cycles.

\subsection{Generalities on local epsilon factors}\label{3.1}

Let $k$ be a perfect 
field of characteristic $p$ and $T$ be a henselian trait which is 
isomorphic to the henselization of $\mathbb{A}^1_k$ at the origin.  We denote its generic point by $\eta$. We denote the absolute Galois group of $k$ (resp. $\eta$) 
by $G_k$ (resp. $G_\eta$). 

Let $(\rho,V)$ be a pair of 
a finite free $\Lambda$-module $V$ together with a continuous group homomorphism $\rho\colon G_\eta\to{\rm GL}(V)$, where ${\rm GL}(V)$ is equipped with the discrete topology. To such a pair, Yasuda \cite{Y3} attaches a continuous character $\tilde{\varepsilon}_{0,\Lambda}(V,\tilde{\psi})\colon G_s\to\Lambda^\times$, called the local epsilon factor,  as a generalization of Langlands--Deligne's theory \cite{Del}. This character depends on 
a choice of $\tilde{\psi}$, which is what Yasuda calls a {\it non-trivial additive character sheaf} in loc.~cit. In the sequel, we explain how to attach
an invertible additive character sheaf $\tilde{\psi}_\omega$ to 
$\psi\colon\mathbb{F}_p\to\Lambda^\times$ and $\omega\in\Omega^1_\eta$. Then we put 
$\varepsilon_{0,\Lambda}(T,V,\omega)=\tilde{\varepsilon}_{0,\Lambda}(V,\tilde{\psi}_\omega)$; we use a similar notation as Laumon's one  given in \cite{Lau}. 

Let $F$ be the completion of the fraction field of $T$. For the definition of non-trivial 
 additive character sheaf on $F$, we refer to \cite[4.1]{Y3}. 
Let $\omega$ be a non-zero rational $1$-form on $\eta$. We define a non-trivial additive character sheaf 
$\tilde{\psi}_\omega$ as follows. 
When the residue field $k$ is finite, $\omega$ defines an additive character $\psi_\omega\colon F\to\Lambda^\times$ by 
$a\mapsto\psi({\rm Tr}_{k/\mathbb{F}_p}\circ {\rm Res}(a\omega))$. Then $\tilde{\psi}_\omega$ is the 
additive character sheaf corresponding to $\psi_\omega$, which is constructed in \cite[Corollary 4.3]{Y3}. 
In general, take a uniformizer 
$\pi\in F$, which defines an inclusion $\mathbb{F}_p((\pi))\to F$. When  $\omega=d\pi$, define the 
non-trivial additive character sheaf 
$\tilde{\psi}_\omega$ on $F$ to be the pull-back of $\tilde{\psi}_{d\pi}$ on 
$\mathbb{F}_p((\pi))$. When $\omega=ad\pi$ for a non-zero element $a\in F$, 
define $\tilde{\psi}_\omega$ to be the pull-back of $\tilde{\psi}_{d\pi}$ on $F$
 by the multiplication-by-$a$ map 
$F\to F$. Using \cite[Corollary 4.4.2]{Y3}, we can show that this definition of $\tilde{\psi}_\omega$ is 
independent of an auxiliary choice of $\pi$.

The assignment $(T,(\rho,V),\omega)\mapsto \varepsilon_{0,\Lambda}(T,V,\omega)$ satisfies the following properties. 

\begin{thm}(\cite{Y1}, \cite[4.12]{Y3})\label{torloc}
Let the notation be as above. For a triple $(T,(\rho,V),\omega)$ where $V$ is a finite free $\Lambda$-module with a continuous 
group homomorphism $\rho\colon G_\eta\to{\rm GL}(V)$ and $\omega\in\Omega^1_\eta$ is a non-zero rational $1$-form, the local epsilon factor $\varepsilon_{0,\Lambda}(T,V,\omega)\colon G_k^{ab}\to\Lambda^\times$ satisfies the following properties. 
\begin{enumerate}
\item The character only depends on the isomorphism class of $(T,(\rho,V),\omega)$. 
\item For a short exact sequence $0\to V'\to V\to V''\to0$ of representations of $G_\eta$, we have 
\begin{equation*}
\varepsilon_{0,\Lambda}(T,V,\omega)=\varepsilon_{0,\Lambda}(T,V',\omega)\cdot\varepsilon_{0,\Lambda}(T,V'',\omega). 
\end{equation*}
\item For a local ring homomorphism $f\colon \Lambda\to\Lambda'$, we have 
\begin{equation*}
f\circ\varepsilon_{0,\Lambda}(T,V,\omega)=\varepsilon_{0,\Lambda'}(T,V\otimes_\Lambda\Lambda',\omega)
\end{equation*}
as characters $G_k\to\Lambda'^\times$. Here the character $\mathbb{F}_p\to\Lambda'^\times$ used in defining $\varepsilon_{0,\Lambda'}$ is taken to be $f\circ\psi$. 
\item We have 
\begin{equation*}
\varepsilon_{0,\Lambda}(T,V,\omega)\cdot\varepsilon_{0,\Lambda}(T,V,\omega')^{-1}
=\det(V)_{[
\frac{\omega}{\omega'}]}\chi_\cyc^{({\rm ord}(\omega')-{\rm ord}(\omega)){\rm rk}V}. 
\end{equation*}
Here $\chi_\cyc$ denotes the composition of the $\ell$-adic cyclotomic character $G_k^{ab}\to\mathbb{Z}_\ell^\times$ and the canonical map 
$\mathbb{Z}_\ell^\times\to\Lambda^\times$ and 
$k(\eta)^\times\times{\rm H}^1(\eta,\Lambda^\times)\to{\rm H}^1(k,\Lambda^\times)$, $(a,\chi)\mapsto\chi_{[a]}$ is  
the pairing defined in \cite[4.2]{Y3}, \cite[Definition 3.12]{CLE}. 
\item Let $W$ be an unramified representation of $G_\eta$ on a finite free $\Lambda$-module. 
We have 
\begin{equation*}
\varepsilon_{0,\Lambda}(T,V\otimes W,\omega)=
\det(W)^{\otimes a(T,V,\omega)}\cdot
\varepsilon_{0,\Lambda}(T,V,\omega)^{{\rm rk}W}. 
\end{equation*}
Here $a(T,\mathcal{F},\omega):={\rm Sw}V+\rk V({\rm ord}\omega+1)$. 
\item Assume that $k$ is finite and that there exists a local ring morphism 
$f\colon \mathcal{O}_E\to\Lambda$ from the ring of integers of a finite extension $E/\mathbb{Q}_\ell$ 
such that $V$ comes from a representation on $\mathcal{O}_E$, i.e. there is a representation $V'$ of $G_\eta$ 
on a finite free $\mathcal{O}_E$-module such that $V'\otimes_{\mathcal{O}_E}\Lambda\cong V$. 
Then we have 
\begin{equation*}
\varepsilon_{0,\Lambda}(T,V,\omega)({\rm Frob_k})=(-1)^{\rk V+{\rm Sw}V}f(\varepsilon_{0}(T,V'\otimes_{\mathcal{O}_E}E,\omega)). 
\end{equation*}
Here the local epsilon factor in the right hand side is the classical one  in \cite[Th\'eor\`eme (3.1.5.4)]{Lau}. 
\end{enumerate}
\end{thm}

If no confusions occur, we omit the subscript $\Lambda$ in $\varepsilon_{0,\Lambda}(T,V,\omega)$.
In what follows, we identify $D(\eta,\Lambda)$ with the derived category of discrete $\Lambda$-modules with continuous $G_\eta$-actions as the \'etale topos of $\eta$ is canonically equivalent to the category of discrete sets with continuous $G_\eta$-actions. 
By the multiplicativity in Theorem \ref{torloc}.2, we can  define $\varepsilon_{0,\Lambda}(T,K,\omega)$ for 
a constructible complex $K\in D_{\rm ctf}(\eta,\Lambda)$. In the following lemma, we recall  the precise definition and prove several basic properties. 
Before proceeding, let us recall that $K$ is said to be 
{\it strictly perfect} if $K^i$ are finite free $\Lambda$-modules for all $i$ and are $0$ for all but finitely many $i$. 
\begin{lm}\label{consep}
Let $K\in D_{\rm ctf}(\eta,\Lambda). $
\begin{enumerate}
\item There exists a strictly perfect complex $K'\in D_{\rm ctf}(\eta,\Lambda)$ that is quasi-isomorphic to $K$. Further if $K$ is bounded below (i.e., $K^i=0$ for $i\ll0$), then we can find 
 an actual morphism $K'\to K$ that is a quasi-isomorphism. Here an actual morphism means a morphism of complexes (not in the homotopy nor  derived category). 
\item We put $\varepsilon_{0,\Lambda}(T,K,\omega):=
\prod_i\varepsilon_{0,\Lambda}(T,(K')^i,\omega)^{(-1)^i}$, where $K'$ is a complex chosen in $1$. Then, $\varepsilon_{0,\Lambda}(T,K,\omega)$ is independent of a choice of $K'$. 
\item Let $K_1\to K_2\to K_3\to$ be a distinguished triangle in $D_{\rm ctf}(\eta,\Lambda)$. Then we have 
\begin{equation*}
\varepsilon_{0,\Lambda}(T,K_1,\omega)=\varepsilon_{0,\Lambda}(T,K_1,\omega)\cdot\varepsilon_{0,\Lambda}(T,K_3,\omega). 
\end{equation*}
\item Let $f\colon \Lambda\to\Lambda'$ be a local homomorphism between finite local rings. Then we have 
\begin{equation*}
f\circ\varepsilon_{0,\Lambda}(T,K,\omega)=\varepsilon_{0,\Lambda'}(T,K\otimes_\Lambda^L\Lambda',\omega). 
\end{equation*}
\end{enumerate}
\end{lm}
\proof{
$1$.  Let $M$ be a $\Lambda$-module with a $G_\eta$-action that is continuous with respect to the discrete topology on $M$. For a finite sequence $v_1,\dots,v_n$ of elements in $M$, there exists an 
open subgroup $H<G_\eta$ such that the map $\Lambda[G]^{\oplus n}\to M$ of $\Lambda$-representations sending the basis $e_i$ to $v_i$ factors through the discrete quotient $\Lambda[G/H]^{\oplus n}$. Using this fact and the finiteness of $H^i(K)$, we can find a bounded above  complex $K_1$ of $G_\eta$-representations on finite free $\Lambda$-modules endowed with a map $K_1\to K$ that is a quasi-isomorphism. Suppose 
that $K$ has tor-amplitude in $[a,b]$. 
Then, for an integer $c\leq a$, 
the truncation $\tau_{\geq c}K_1$ has the desired  properties. 
Further suppose that $K$ is bounded below. Take 
an integer $c\leq a$ such that $K^i=0$ for $i<c$. Then, the map $K_1\to K$ factors through $\tau_{\geq c}K_1$. Therefore we find an actual morphism  $K':=\tau_{\geq c}K_1\to K$. 

$2.$ Let $K''$ be another choice. We must show the equality 
\begin{equation*}
\prod_i\varepsilon_{0,\Lambda}(T,(K')^i,\omega)^{(-1)^i}
=\prod_i\varepsilon_{0,\Lambda}(T,(K'')^i,\omega)^{(-1)^i}. 
\end{equation*}
In the derived category, there is a quasi-isomorphism $K'\to K''$. This means that there are actual morphisms of complexes $K'\leftarrow L\to K''$ that are quasi-isomorphisms. Since $K'$ and $K''$ are bounded, we may take $L$ to be bounded. Then, by $1$, there exists an actual morphism $L'\to L$ from  a 
strictly perfect $L'$ that is a quasi-isomorphism. 
Thus, we find actual morphisms of strictly perfect  complexes $K'\leftarrow L'\rightarrow K''$ that are quasi-isomorphisms. 

Therefore, we may assume that there exists an actual morphism $K'\to K''$ that is a quasi-isomorphism. Let 
$C$ denotes its mapping cone. Since $C$ is acyclic, 
$C$ decomposes into finitely many short exact sequences of finite free $\Lambda$-modules with $G_\eta$-actions.  Then Theorem \ref{torloc}.2 implies that $\prod_i\varepsilon_{0,\Lambda}(T,C^i,\omega)^{(-1)^i}=1$. As we have $C^i=(K'')^i\oplus(K')^{i+1}$, the equality follows.

$3$. By the same argument in $2$, we may assume that the map $K_1\to K_2$ in $D_{\rm ctf}(\eta,\Lambda)$ comes from an actual morphism. Then, 
$K_3$ is quasi-isomorphic to the mapping cone of this map. The assertion follows from the well-definedness in $2$. 

$4$. Let $K'$ be a strictly perfect complex that is quasi-isomorphic to $K$. Then, $K\otimes_\Lambda^L\Lambda'$ is quasi-isomorphic to 
$K'\otimes_\Lambda\Lambda'$. The assertion is a consequence of Theorem \ref{torloc}.3. 
\qed
}

Let $\mathcal{O}_E$ be the ring of integers of a finite extension $E/\mathbb{Q}_\ell$. Take and fix a non-trivial 
character $\psi\colon\F_p\to\mathcal{O}_E^\times$. Using Lemma \ref{consep}, 
we define a local epsilon factor $\varepsilon_{0,\Lambda}(T,\mathcal{F},\omega)\colon G_k\to\mathcal{O}_E^\times$ for 
$\mathcal{F}\in D^b_c(\eta,\mathcal{O}_E)$ as follows: Put $\Lambda_n:=\mathcal{O}_E/\ell^{n+1}$. 
Then the reduction $\mathcal{F}_n:=\mathcal{F}\otimes^L_{\mathcal{O}_E}\Lambda_n$ belongs to $D_{\rm ctf}(\eta,\Lambda_n)$ 
and $\psi$ induces a non-trivial character $\F_p\to\Lambda_n^\times$. By Lemma \ref{consep}, we have a character $\varepsilon_{0,\Lambda_n}
(T,\mathcal{F}_n,\omega)\colon G_k\to\Lambda_n^\times$. Since 
$\varepsilon_{0,\Lambda_n}(T,\mathcal{F}_n,\omega)$ and $ \varepsilon_{0,\Lambda_{n+1}}
(T,\mathcal{F}_{n+1},\omega)$ are compatible with the quotient map $\Lambda_{n+1}\to\Lambda_n$ by Lemma \ref{consep}.4, it is allowed to 
define $\varepsilon_{0,\mathcal{O}_E}(T,\mathcal{F},\omega):=\varprojlim_n\varepsilon_{0,\Lambda_n}
(T,\mathcal{F}_n,\omega)$. This is a continuous character of $G_k$ valued in ${\cal O}_E^\times$. 
Finally we give the  definitions of local epsilon factors for $\Zl$-sheaves. 
\begin{df-lm}\label{zlep}
Let the notation be as above. 
Let $\Zl$ be the integral closure of $\mathbb{Z}_\ell$ in an algebraic closure $\Ql$ of $\mathbb{Q}_\ell$. 
\begin{enumerate}
\item 
For a complex 
$\mathcal{F}\in D^b_c(\eta,\Zl)$, we define a character $\varepsilon_{0,\Zl}(T,\mathcal{F},\omega)\colon 
G_k^{ab}\to\Zl^\times$ as follows. 
By definition, there exists a finite subextension $E$ of $\Ql/\mathbb{Q}_\ell$ such that 
$\mathcal{F}$ is defined over $\mathcal{O}_E$, i.e., there exists a complex $\mathcal{G}\in D^b_c(\eta,\mathcal{O}_E)$ 
with $\mathcal{G}\otimes^L_{\mathcal{O}_E}\Zl\cong\mathcal{F}$. We define $\varepsilon_{0,\Zl}
(T,\mathcal{F},\omega)$ to be the composition of 
\begin{equation*}
G_k\xrightarrow{\varepsilon_{0,\mathcal{O}_E}(T,\mathcal{G},\omega)}\mathcal{O}_E^\times\to\Zl^\times. 
\end{equation*}
 This does not depend on choices of $E$ and $\mathcal{G}$. 
 \item For $\mathcal{F}\in D^b_c(T,\Zl)$, we define 
$\varepsilon(T,\mathcal{F},\omega)\colon G_k^{ab}\to\Zl^\times$ to be the product 
\begin{equation*}
\varepsilon_{0}(T,\mathcal{F}_\eta,\omega)\cdot \det(\mathcal{F}_s)^{-1}, 
\end{equation*}
where $s$ denotes the closed point of $T$. 
\end{enumerate}
\end{df-lm}
\proof{
We verify the independence of a choice of $(E,{\cal G})$. Take another choice $(E',{\cal G}')$. Then there exists a finite extension $E''$ of $E$ and $E'$ over which ${\cal G}$ and ${\cal G}'$ are isomorphic, i.e., ${\cal G}\otimes_{{\cal O}_E}{\cal O}_{E''}\cong
{\cal G}'\otimes_{{\cal O}_{E'}}{\cal O}_{E''}$. Then the independence follows from Lemma \ref{consep}.4. 
\qed}

\begin{rmk}\label{eqepep}
Recently, Q. Guignard gives another definition and construction of local epsilon factors \cite{geomep}, using 
Gabber-Katz canonical extension. The local epsilon factors given in Definition \ref{zlep}.2 are the same as his, 
because both of them coincide with the one defined from Laumon's local Fourier transform. 
See \cite[Theorem 11.8]{geomep} and  \cite[Proposition 8.3]{Y2}. 
\end{rmk}

Using these local epsilon factors, we have a generalization of the product formula of Laumon \cite[Th\'eor\`eme (3.2.1.1)]{Lau} to 
a general perfect field. 

\begin{thm}(\cite[Theorem 11.1]{geomep},\cite[Theorem 4.50]{Y3})\label{clprod}
Let $X$ be a connected projective smooth curve over a perfect field $k$ 
of characteristic $p>0$. Let 
$\mathcal{F}\in D^b_c(X,\Zl)$ be a  constructible complex on $X$. Fixing a  
non-zero rational $1$-form $\omega$ on $X$, we have 
\begin{equation*}
\det(R\Gamma(X_{\bar{k}},\mathcal{F}))^{-1}=
\chi_{\cyc}^{-\frac{1}{2}\chi(X)\cdot\rk\mathcal{F}}\prod_{x\in|X|}
\delta_{x/k}^{a(X_{(x)},\mathcal{F})}\varepsilon(X_{(x)},\mathcal{F},
\omega)\circ {\rm tr}_{x/k} 
\end{equation*}
as a character of the absolute Galois group $G_k$ of $k$. 
Here $\chi(X)=\sum_i(-1)^i{\rm dim}{\rm H}^i(X_{\bar
{k}},\Ql)$ is the Euler-Poincar\'e characteristic, $\rk\mathcal{F}$ is the generic rank of $\mathcal{F}$, 
$a(X_{(x)},\mathcal{F})=\rk\mathcal{F}_\eta+{\rm Sw}_x\mathcal{F}-
\rk\mathcal{F}_x$ is  the Artin conductor. ${\rm tr}_{x/k}\colon G_k^{ab}
\to G_{k(x)}^{ab}$ denotes the transfer morphism and $\delta_{x/k}$ denotes the determinant character of the 
induced representation ${\rm Ind}_{G_{k(x)}}^{G_k}
\mathbb{Q}_\ell$ of the trivial representation. 
\end{thm}
\proof{
In \cite[Theorem 11.1]{geomep}, Guignard considers (a complex of) $\Ql$-sheaves. To prove the assertion, we can simply apply his result to ${\cal F}\otimes_{\Zl}^L\Ql$, once we know that the local epsilon factors given by Guignard and Yasuda coincide (see Remark \ref{eqepep}). 
\cite[Theorem 4.50]{Y3}, Yasuda considers a finite coefficient case and proves the formula when $\cal F$ is of the form $j_!{\cal G}$ where $j\colon U\to X$ is a dense open immersion and $\cal G$ is a locally constant sheaf on $U$, from which our result also follows by multiplicativity. 
\qed}

We generalize local epsilon factors  for tamely ramified representations 
to the case of characteristic $0$, which is done as below at the expense of ignoring roots of unity. 

Let $S$ be an affine (not necessarily noetherian) 
normal scheme in which $\ell$ is invertible. 
Consider a pair $(T,\chi)$ such that $T=(T_i)_i$ is a finite family of finite \'etale 
coverings of $S$ 
and $\chi=(\chi_i)_i$ is a family of characters $\chi_i\colon\mathbb{Z}
/d_i\mathbb{Z}(1)\to
\Ql^\times$ of \'etale sheaves on $T_i$ 
where $d_i$ are integers $\geq1$ invertible in $S$ such that $\mathbb{Z}/d_i\mathbb{Z}(1)
\cong \mathbb{Z}/d_i\mathbb{Z}$ as \'etale sheaves on $T_i$. 
Denote by $N_{T_i/S}(\chi_i)$ the character 
$\mathbb{Z}/d_i\mathbb{Z}(1)\to\Ql^\times$ of \'etale sheaves on $S$ 
defined by the composition 
$\mathbb{Z}/d_i\mathbb{Z}(1)\to f_{i\ast}\mathbb{Z}/d_i\mathbb{Z}(1)
\xrightarrow{f_{i\ast}\chi_{i}}f_{i\ast}\Ql^\times\xrightarrow{{\rm {\rm tr}}}\Ql^\times$ 
where $f_{i}\colon T_{i}\to S$ is the structure morphism and the former arrow is given by the adjunction. For an integer $N\geq1$ 
which is a multiple of $d_i$, we regard $N_{T_i/S}(\chi_i)$ as a character 
of $\mathbb{Z}/N\mathbb{Z}(1)$ via the surjection 
$\mathbb{Z}/N\mathbb{Z}(1)\to\mathbb{Z}/d_i\mathbb{Z}(1)$, 
$a\mapsto a^{\frac{N}{d_i}}$. 

Assume that $\prod_iN_{T_i/S}(\chi_i)$ is trivial where the product is 
taken as characters of $\mathbb{Z}/N\mathbb{Z}(1)$ 
for some common multiple 
$N$ of $d_i$. 
 In this case, $(T,\chi)$ is called a 
Jacobi datum in \cite[Section 1]{Jacob}. 
When $S$ is the spectrum of a finite field $k$ with $q$ elements, 
Saito attaches a Jacobi sum $j_\chi\in \Ql^\times$ 
to a Jacobi datum $(T,\chi)=((T_i)_i,(\chi_i)_i)$ in 
\cite[Section 2]{Jacob} by the following formula: 
\begin{equation}\label{jas}
j_\chi:=\prod_{i}(\prod_j\tau_{k_{ij}}(\bar{\chi}_{ij},\psi_0\circ{\rm Tr}_{k_{ij}/k})). 
\end{equation}
Here $T_i=\coprod_j{\rm Spec}(k_{ij})$ for finite fields $k_{ij}$ with $q_{ij}$ 
elements, $\bar{\chi}_{ij}\colon k_{ij}^\times\to\Ql^\times$ is defined 
by $a\mapsto\chi_i(a^{(q_{ij}-1)/d_i})$, and $\psi_0\colon k\to\Ql^\times$ 
is a nontrivial character (note that, since $\mathbb{Z}/d_i\mathbb{Z}(1)$ is assumed to be trivial on $T_i$, $q_{ij}-1$ is divisible by $d_i$). 
 The Gauss sums are defined by $\tau_{k}(\chi,\psi)=-\sum_{a\in k}\chi^{-1}(a)\psi(a)$. Since $\prod_iN_{T_i/S}(\chi_i)$ is trivial, 
the Jacobi sum $j_\chi$ is independent of a choice of $\psi_0$.

Let $(T,\chi)$ be a Jacobi datum on an affine normal scheme $S$. In 
 \cite[Proposition 2.]{Jacob}, Saito constructed a 
 smooth $\Ql$-sheaf $J_\chi$ of rank $1$ on $S$ from the Jacobi datum, 
 which is called a Jacobi sum character. 
This 
is characterized by the following properties. 
\begin{itemize}
\item For every morphism $f\colon S'\to S$ of an affine normal schemes, we have 
$f^\ast J_\chi\cong J_{f^\ast\chi}$. 
\item If $S$ is the spectrum of a finite field $\F_q$, then the action 
of the geometric Frobenius on $J_\chi$ is given by the multiplication by 
$j_\chi$ (\ref{jas}). 
\end{itemize}

Let $k$ be a perfect field of characteristic $p\geq0$  different from $\ell$. 
Take and fix an algebraic closure $\bar{k}$ of $k$ and 
let $I:=\varprojlim_{n}\mu_n(\bar{k})$, 
where $n$ runs through integers $\geq1$ invertible in $k$ and 
$\mu_n(\bar{k})$ is the group of $n$-th roots of unity in $\bar{k}$. 
The group $I$ carries an action of $G_k={\rm Gal}(\bar{k}/k)$. 

Let $V$ be a finite free $\Zl$-module equipped with a 
continuous homomorphism $\rho\colon I\to{\rm GL}(V)$, where ${\rm GL}(V)$ is endowed with the topology coming from the $\ell$-adic topology on $\Zl$. 
Put $V_{\Ql}:=V\otimes_{\Zl}\Ql$. 
For an element $\sigma\in G_k$, 
let $\sigma^\ast V_{\Ql}$ denote the representation of $I$ defined by 
$\rho\circ\sigma$. When $\rho$ factors through the quotient 
$I\to\mu_n(\bar{k})$ for some $n$, 
the twist $\sigma^\ast V_{\Ql}$ depends only on the image of $\sigma$ in 
${\rm Gal}(k(\mu_n(\bar{k}))/k)$. In this case, for $\tau\in{\rm Gal}(k(\mu_n(\bar{k}))/k)$, 
$\tau^\ast V_{\Ql}$ means the twist $\sigma^\ast V_{\Ql}$ for any lift $\sigma\in G_k$ of $\tau$. 

Assume that, for each $\sigma\in G_k$, we have 
$\sigma^\ast V_{\Ql}\cong V_{\Ql}$ and that $V_{\Ql}$ is potentially unipotent, 
i.e., there exists an open subgroup $I'\subset I$ which acts on $V_{\Ql}$ unipotently. 
Then the semi-simplification $V^{ss}_{\Ql}$ 
decomposes into a direct sum 
\begin{equation}\label{sstame}
V^{ss}_{\Ql}\cong\bigoplus_i(\bigoplus_{\tau\in{\rm Gal}(k_i/k)}\tau^\ast\chi_i)
\end{equation}
where $\chi_i\colon\mu_{d_i}(\bar{k})\hookrightarrow\Ql^\times$ is an 
injective character 
and we put $k_i:=k(\mu_{d_i}(\bar{k}))$. 
Such a decomposition is unique up to permutation. 
Note that the determinant $\det(V_{\Ql})=\det(V)$ equals to 
$\prod_iN_{k_i/k}(\chi_i)$. 
\begin{df}\label{js}
Let the notation be as above. Assume that $\sigma^\ast V_{\Ql}$ is isomorphic 
to $V_{\Ql}$ for each $\sigma\in G_k$ and that 
$V_{\Ql}$ is potentially unipotent. Let $\mu$ denotes the group of roots of unity in $\Zl$. 
\begin{enumerate}
\item When the determinant $\det(V)$ is the trivial character of $I$, 
we denote by $J(V)$ the Jacobi sum character of the Jacobi datum 
$(({\rm Spec}(k_i))_i,(\chi_i)_i)$. 
This is a character $G_k\to\Zl^\times$. 
\item In general, we define $J(V)$ to be 
$\iota_n\circ J(V^{\oplus n})$ where $n$  is an integer $\geq1$ such that 
$\det(V^{\oplus n})$ is trivial and $\iota_n\colon \Zl^\times\to\Zl^\times/\mu$ is the map defined by $x\mapsto\sqrt[n]{x}$. 
This is a group homomorphism 
$G_k\to \Zl^\times/\mu$ and independent of a choice of $n$. 
\end{enumerate}

$J(V)\colon G_k\to\Zl^\times/\mu$ is admissible in the terminology of Subsection 4.1. 
\end{df}

Let $T$ be the henselization of $\mathbb{A}^1_k$ at the origin. Let $\eta$ be the generic point of $T$ 
and fix a separable closure $\overline{k(\eta)}$ of $k(\eta)$. 
We take $\bar{k}$ to be the algebraic closure of $k$ in $\overline{k(\eta)}$. 
Let $I$ be the tame inertia group of $G_\eta={\rm Gal}(\overline{k(\eta)}/k(\eta))$. 
This is canonically isomorphic to $\varprojlim_{p\nmid n}\mu_n(\bar{k})$. 
Let $V$ be a smooth sheaf of finite free $\Zl$-modules  on $\eta$ which is tamely ramified. 
Then, as a representation of $I$, $V_{\overline{\eta}}$ is isomorphic 
to $\sigma^\ast V_{\overline{\eta}}$ for $\sigma\in G_k$. 

We give a definition of local epsilon factors modulo roots of unity especially in the case characteristic $0$. 

\begin{df}\label{epchar00}
Let the notation be as above. Let $V$ be 
a smooth sheaf of finite free $\Zl$-modules on $\eta$ which is tamely ramified 
and potentially unipotent. Then there exists an integer $n\geq1$ such that $\det(V)^n$ is unramified. 
We use the same symbol $\det(V)^n$ for the corresponding character 
$G_k\to\Zl^\times.$ 
We 
define the local epsilon factor $\overline{\varepsilon}_0(T,V)\colon
G_k\to\Zl^\times/\mu$ by setting 
\begin{equation*}
\overline{\varepsilon}_0(T,V):=\iota_n\circ(\det(V)^n) \cdot J(V_{\overline{\eta}}). 
\end{equation*}
Here $\iota_n$ is the map in Definition \ref{js}.2. 
This is independent of a choice of $n$. 
\end{df}
\begin{lm}
Let $V$ and $W$ be smooth sheaves of finite free $\Zl$-modules on $\eta$. Assume that 
$V$ is unramified and that $W$ is tamely ramified and potentially unipotent. 
We have 
\begin{equation*}
\overline{\varepsilon}_0(T,V\otimes W)=
(\det V)^{{\rm dim}W}\cdot
\overline{\varepsilon}_0(T,W)^{{\rm dim}V}. 
\end{equation*}
Here, in the right-hand side, $\det V$ denotes the character of $G_k$ defined 
from the unramified character $\det V$ by abuse of notation. 
\end{lm}
\proof{
Since $V$ is unramified, we have $J((V\otimes W)_{\overline{\eta}})=
J(W_{\overline{\eta}})^{{\rm dim}V}$. On the other hand, we have 
$\det(V\otimes W)=\det V^{{\rm dim}W}\cdot
\det W^{{\rm dim}V}$. 
The assertion follows. 
\qed
}

\begin{lm}\label{sameep}
Assume that $k$ is a perfect field of characteristic $p>0$. 
Let $V$ be a tamely ramified smooth sheaf of finite free $\Zl$-modules on $\eta$ which 
is potentially unipotent. 
Let $\pi$ be a uniformizer of  $T$. 
Then the character $\overline{\varepsilon}_0(T,V)$ in Definition 
\ref{epchar00} coincides with 
the character $\varepsilon_0(T,V,d\pi)$ in Definition \ref{zlep}.1 
followed by $\Zl^\times\to\Zl^\times/\mu$. 
\end{lm}
\proof{
Put $V_{\Ql}:=V\otimes_{\Zl}\Ql$. 
We may assume that $V_{\Ql}$ is irreducible. 
Let $\chi$ be 
a character of $I$ which appears in the decomposition of $V_{\Ql}$ as a representation of $I$. 
Let $n$ be the cardinality of the image of $\chi$. This is an integer $\geq1$ prime to $p$ and $\chi$ decomposes as $I\to \mu_n(\bar{k})\hookrightarrow\Ql^\times$, where $\mu_n(\bar{k})$ denotes the group of $n$-th roots of unity. Let 
$\eta_n$ be the unramified extension of $\eta$ with residue field 
$k(\mu_n(\bar{k}))$. The $\chi$-isotypic part $V_\chi$ of 
$V_{\Ql}$ is stable under the action of $G_{\eta_n}$ and 
$V_{\Ql}\cong{\rm Ind}_{G_{\eta_n}}^{G_\eta}V_\chi$. Hence we reduce 
it to the case when $\eta_n=\eta$. In this case, 
we also denote by $\chi$ the character of $G_\eta$ induced from the 
identification 
${\rm Gal}(k(\eta)[\pi^{\frac{1}{n}}]/k(\eta))\cong\mu_n(\bar{k})$. 
Then we have $V_{\Ql}\cong\chi\otimes V_0$ where $V_0$ is an unramified  
representation. By \cite[Proposition (2.5.3.1)]{Lau}, 
we have $F^{(0,\infty)}(V_{\Ql})\cong V_{\Ql}\otimes G(\chi,\psi)$ where $F^{(0,\infty)}$ is the local Fourier transform 
\cite[D\'efinition (2.4.2.3)]{Lau} and 
$G(\chi,\psi)$ is the $1$ dimensional representation  defined in loc.~cit. 
Hence we have $\det F^{(0,\infty)}(V_{\Ql})\cong\det(V_{\Ql})\otimes G(\chi,\psi)^{{\rm dim}V}$. 
On the other hand, for an integer $m\geq1$ such that 
$(\det V)^m$ is unramified, we have 
$J(V^{\oplus m})\cong G(\chi,\psi)^{\otimes m\cdot\rk V}$. 
The assertion follows from Laumon's formula \cite[Th\'eor\`eme (3.5.1.1)]{Lau}. 
\qed
}

\subsection{Reminder on oriented product and local Fourier transform (\cite{CLE})}\label{recOLF}
In this preliminary subsection, we briefly recall the notion of oriented products and the local Fourier transform rephrased in terms of oriented products, which are discussed in detail in \cite{CLE}. 
The results recalled in this subsection are only used in the subsections \ref{redpos} and \ref{locefconv}. 

For morphisms $f\colon X\to S,g\colon Y\to S$ of topoi, the oriented product $X\overset\gets\times_SY$ is a topos equipped with projections $p_X\colon X\overset\gets\times_SY\to X, p_Y\colon 
X\overset\gets\times_SY\to Y$ and a natural transformation $\sigma\colon gp_Y\to fp_X$. 
This topos is characterized by a universality in terms of such a triple. 
We refer to \cite{ori}  
for its precise definition and basic properties. 

Let $\Lambda$ be a finite local ring with characteristic invertible in schemes in what follows. 
Let $f\colon X\to S$ be a morphism of schemes. We use the same symbols $X,S$ for their associated \'etale topoi by abuse of notation. By the universality of 
$X\overset\gets\times_SS$, we get a morphism of topoi $\Psi_f\colon X\to X\overset\gets\times_SS$ which fits into the  $2$-commutative diagram of topoi 
\begin{equation*}
\xymatrix{
X\ar[rd]^-{\Psi_f}\ar[rdd]_-{{\rm id}_X}\ar[rrd]^-f&&\\
&X\overset\gets\times_SS\ar[d]^-{p_X}\ar[r]_-{p_S}&S\\
&X. 
}
\end{equation*}
The derived functor $R\Psi_f\colon D^+(X,\Lambda)\to D^+(X\overset\gets\times_SS,\Lambda)$ is called {\it the nearby cycles functor}.  
The relation $p_X\Psi_f={\rm id}_X$ gives a morphism $p_X^\ast\to R\Psi_f$ by adjunction.  {\it The vanishing cycles functor} 
$R\Phi_f\colon D^+(X,\Lambda)\to D^+(X\overset\gets\times_SS,\Lambda)$ is a triangulated functor which fits into a distinguished triangle 
\begin{equation}\label{distp}
p_X^\ast\to R\Psi_f\to R\Phi_f\to. 
\end{equation}
 A cartesian diagram 
\begin{equation*}
\xymatrix{
X'\ar[r]^-{f'}\ar[d]&S'\ar[d]^-g\\X\ar[r]^-{f}&S
}
\end{equation*}
of schemes induces a morphism of topoi $\overset\gets{g}\colon X'\overset\gets\times_{S'}S'\to
X\overset\gets\times_{S}S$. 
For an element 
 $K\in D^+(X,\Lambda)$, we have a canonical base change map 
 $\overset\gets{g}^\ast R\Psi_f(K)\to R\Psi_{f'}(g^\ast K)$.

The relation of $R\Psi_f, R\Phi_f$ given above with the classical nearby and vanishing cycles formalisms in \cite{SGA7-2} is explained in \cite[1.2]{Ill}. 
When $S$ is a henselian trait with closed and generic points being denoted by $s$ and $\eta$,  the closed  subtopos $X_s\overset\gets\times_S S\subset X\overset\gets\times_SS$, where 
$X_s$ denotes the special fiber, is 
canonically identified with the topos $X_s\times_sS$ appeared in \cite{SGA7-2} and the restrictions $R\Psi_f|_{X_s\overset\gets\times_S\eta}$ and $
R\Phi_f|_{X_s\overset\gets\times_S\eta}$ to the open subtopos 
$X_s\overset\gets\times_S\eta\subset X_s\overset\gets\times_SS$ recover the classical nearby and vanishing cycles functors.

In Definition \ref{rellocfou} below, we recall the generalization of the local Fourier transforms in \cite{Lau} to a relative setting in terms of oriented product. We refer to \cite[Section 3]{CLE} for a detailed account.  

Let $S$ be a noetherian scheme over $\F_p$. 
For a non-trivial character $\psi\colon\mathbb{F}_p\to\Lambda^\times$, 
let $\mathcal{L}_\psi(x)$ be the Artin--Schreier sheaf on $\mathbb{A}^1_{\mathbb{F}_p}$.  We write 
$\mathcal{L}_\psi(x\cdot x')$ for its pull-back by 
the multiplication $\mathbb{A}^1_{\F_p}\times_{\F_p}\mathbb{A}^1_{\F_p}\to \mathbb{A}^1_{\F_p}$, 
where $x$ and $x'$ stand for the standard coordinates of the first and second 
factors respectively.
We also use the same notation  for 
its pull-back by $\mathbb{A}^1_S\times_S\mathbb{A}^1_S\to\mathbb{A}^1_{\F_p}\times_{\F_p}\mathbb{A}^1_{\F_p}$. 
Denote by $0_S$ (resp. $\infty_S$) the section $S\to\mathbb{A}^1_S$ (resp. 
$S\to\Proj^1_S$) at the origin (resp. at the infinity). Consider the morphisms of topoi  
\begin{equation*}
0_S\overset\gets\times_{\mathbb{A}^1_S}\mathbb{A}^1_S
\xleftarrow{\overset\gets p_1}(0_S\times_S\infty_S)\overset\gets\times_{\mathbb{A}^1_S\times_S\Proj^1_S}(\mathbb{A}^1_S\times_S\mathbb{A}^1_S)\xrightarrow{\overset\gets p_2}
\infty_S\overset\gets\times_{\Proj^1_S}\mathbb{A}^1_S
\end{equation*}
 induced from the projections 
$\mathbb{A}^1_S\xleftarrow{p_1}\mathbb{A}^1_S\times_S\Proj^1_S\xrightarrow{p_2}\Proj^1_S$. 
Let $q\colon (0_S\times_S\infty_S)\overset\gets\times_{\mathbb{A}^1_S\times_S\Proj^1_S}(\mathbb{A}^1_S\times_S\mathbb{A}^1_S)\to\mathbb{A}^1_S\times_S\mathbb{A}^1_S$ be the second projection. 

For the notion of constructibility of sheaves on oriented products, see \cite[9.1]{Org}, \cite[1.6]{Ill}. 
\begin{df}\label{rellocfou}(\cite[Definition 3.3]{CLE})
Let $K\in D_{\rm ctf}(0_S\overset\gets\times_{\mathbb{A}^1_S}\mathbb{G}_{m,S},\Lambda)$ be a constructible complex 
of finite tor-dimension whose cohomology sheaves 
are locally constant. 
Define the local Fourier transform $F^{(0,\infty)}(K)$ by 
\begin{equation*}
F^{(0,\infty)}(K):=R\overset\gets p_{2\ast}(\overset\gets p_1^\ast K_!\otimes^L_\Lambda q^\ast 
\mathcal{L}_\psi(x\cdot x'))[1]\in D^b(\infty_S\overset\gets\times_{\Proj^1_S}\mathbb{A}^1_S,\Lambda). 
\end{equation*}
Here $K_!$ denotes the 
$0$-extension of $K$ by $0_S\overset\gets\times_{\mathbb{A}^1_S}\mathbb{G}_{m,S}\hookrightarrow0_S\overset\gets\times_{\mathbb{A}^1_S}\mathbb{A}^1_S$. 
\end{df}

\begin{pr}(cf. \cite[Proposition (2.4.2.2)]{Lau})\label{F}Let the notation be as above. Let $\Lambda_0$ be the residue field of 
$\Lambda$. 
\begin{enumerate}
\item The local Fourier transform $F^{(0,\infty)}(K)$ is of finite tor-dimension. 
Let $\Lambda\to\Lambda\rq{}$ be a ring homomorphism to a finite local ring. Then the canonical map 
\begin{equation*}
F^{(0,\infty)}(K)\otimes^L_\Lambda\Lambda\rq{}\to F^{(0,\infty)}(K\otimes^L_\Lambda\Lambda\rq{})
\end{equation*}
is a quasi-isomorphism. 
\item The formation of $F^{(0,\infty)}(K)$ commutes with arbitrary 
base change $g\colon S'\to S$.  
Namely,  the canonical map 
\begin{equation*}
{g_1}^\ast F^{(0,\infty)}(K)\to F^{(0,\infty)}({g_2}^\ast K)
\end{equation*}
is a quasi-isomorphism 
where $g_1$ denotes the morphism $\infty_{S\rq{}}\overset\gets\times_{\Proj^1_{S\rq{}}}\mathbb{A}^1_{S\rq{}}
\to \infty_S\overset\gets\times_{\Proj^1_S}\mathbb{A}^1_S$ whereas $g_2$ denotes the morphism 
$0_{S\rq{}}\overset\gets\times_{\mathbb{A}^1_{S\rq{}}}\mathbb{G}_{m,{S\rq{}}}\to 0_S\overset\gets\times_{\mathbb{A}^1_S}\mathbb{G}_{m,S}$. 
\item\label{loccf} The local Fourier transform $F^{(0,\infty)}(K)$ is constructible. 
It is locally constant if and only if, 
for each $i\in\mathbb{Z}$, 
the function on $S$ defined by $s\mapsto{\rm dimtot}\mathcal{H}^i(K\otimes^L_{\Lambda}\Lambda_0)|_{\eta_{\bar{s}}}$, 
where $\bar{s}$ is the spectrum of an algebraic closure of $k(s)$ and 
$\eta_{\bar{s}}\cong0\overset\gets\times_{\mathbb{A}^1_{k(\bar{s})}}
\mathbb{G}_{m,k(\bar{s})}$ is the generic point of the 
henselization $\mathbb{A}^1_{k(\bar{s}),(0)}$ at the origin, 
is locally constant. 
\item When $K$ has tor-amplitude in $[0,0]$, so is $F^{(0,\infty)}(K)$. 
\end{enumerate}
\end{pr}

\proof{
See \cite[Proposition 3.4]{CLE}. 
\qed
}

\vspace{0.1in}

Let $k$ be a perfect field of characteristic $p>0$. 
Let $C$ be a smooth curve over $k$ and $x\in C$ be a closed point. Recall that the oriented product $x\overset\gets\times_C(C\setminus\{x\})$ is canonically equivalent to the \'etale topos of the generic point of the henselization $C_{(x)}$ (\cite[Lemma 2.2]{CLE}). 
Let 
 $T$ and $T^{\prime}$ be henselian traits which are isomorphic to the henselization of $\mathbb{A}^1_k$ at the origin. Fix uniformizers $\pi$ and $\pi^{\prime}$ of $T$ and $T^{\prime}$ respectively. By sending 
 $x\mapsto\pi$ (resp. $
x'\mapsto1/\pi'$), where $x$ (resp. $x'$) is the standard coordinate 
of $\mathbb{A}^1_k$ (resp. $\mathbb{A}^1_k\subset\Proj^1_k$), we identify $T$ (resp. $T'$) with the henselization 
of $\mathbb{A}^1_k$ (resp. $\Proj^1_k$) at $0$ (resp. at $\infty$). 
Let $\eta$ and $\eta^{\prime}$ be 
the generic points of $T$ and $T^{\prime}$ respectively. Under the identifications $\eta\cong 0\overset\gets\times_{\mathbb{A}^1_k}\mathbb{G}_{m,k}$ and $\eta'\cong \infty\overset\gets\times_{\mathbb{P}^1_k}\mathbb{A}_{k}^1$, we may regard $F^{(0,\infty)}$ as a triangulated functor $D_{\rm ctf}(\eta,\Lambda)\to
D_{\rm ctf}(\eta',\Lambda)$. In \cite[3.2]{CLE}, it is explained that this $F^{(0,\infty)}$ satisfies similar properties as the one in \cite{Lau}. See loc.~cit. for details. 

\subsection{Reduction to the case of positive characteristic}\label{redpos}

To compute the local epsilon factors of vanishing cycles in the case of characteristic $0$, we give a method to reduce it to the case of positive characteristic. This method may be compared with spreading-out arguments for \'etale sheaves with finite coefficients. 
This subsection is only necessary for the case of characteristic $0$. 
\begin{rmk}
The following technique is needed since we work with $\ell$-adic sheaves. 
If one could develop a theory of epsilon cycles for local epsilon factors without 
taking modulo roots of unity and one could work with $\Lambda$-sheaves for 
a finite local ring $\Lambda$, then the technique seemed unnecessary. 
\end{rmk}

We start with general lemmas. 

Let us fix a prime number $\ell$. Let $R$ be a discrete valuation ring in which $\ell$ is invertible. 
Denote by $K$ and $F$ its function field and residue field respectively. 
Fix a uniformizer $\pi\in R$ and, for an integer $m\geq0$, denote by $R_m$ 
the ring $R[\pi^{1/\ell^m}]$ and by $K_m$ the quotient field 
of $R_m$. We write $R_\infty$ and $K_\infty$ 
for the unions $\cup_{m\geq0}R_m$ and $\cup_{m\geq0}K_m$ respectively. 
The rings $R_m$ (including the case $m=\infty$) are valuation rings with residue field being isomorphic to  $F$. 

Let $\mathcal{X}$ be a scheme over $R$. 
Let $m$ be an integer $\geq0$ or $\infty$. 
Consider the diagram 
\begin{equation}\label{ijm}
X_m\xrightarrow{j_m}\mathcal{X}_m\xleftarrow{i_m}\mathcal{X}_F
\end{equation}
where the left arrow is the base change of 
\begin{equation}\label{uladiag}
X:=\mathcal{X}\times_RK\xrightarrow{j}\mathcal{X}
\end{equation}
by $R\to R_m$ 
and $i_m$ is the lift of the special fiber $
\mathcal{X}_F:=\mathcal{X}\times_RF\xrightarrow{i}\mathcal{X}$. 
\begin{lm}\label{ulainf}
Let the notation be as above. 
Let $\Lambda$ be a finite local ring of residue characteristic $\ell$. 
\begin{enumerate}
\item 
For a bounded below complex $C\in D^+(\mathcal{X},\Lambda)$ 
such that the structure morphism $\mathcal{X}\to{\rm Spec}(R)$ 
is locally acyclic relatively to $C$, the canonical map 
$i^\ast C\to i_\infty^\ast Rj_{\infty\ast}j_\infty^\ast C$ induced from (\ref{ijm}) is 
a quasi-isomorphism. 
\item Assume that $\mathcal{X}$ is of finite type over $R$. 
Then, the functor $Rj_{\infty\ast}$ has finite cohomological dimension. 
For a constructible complex $C\in D_c^b(X,\Lambda)$, 
the complex $i^\ast_\infty Rj_{\infty\ast}( C|_{X_\infty})$ on $\mathcal{X}_F$ is 
constructible. 
\end{enumerate}
\end{lm}
In the situation of $2$, for $C\in D_c^b(X,\Lambda)$, we write $\langle C,-\pi\rangle:=i_\infty^\ast Rj_{\infty\ast}(C|_{X_\infty})$. 
When $\mathcal{X}$ and the special fiber $\mathcal{X}_F$ are regular, and $C$ is a locally constant sheaf on which 
the inertia groups at the generic points of the divisor $\mathcal{X}_F$ act through finite $\ell$-groups, this notion coincides with the one given in 
\cite[Definition 2.11]{CLE}. 
\proof{
We may assume that $R$ is strictly henselian, in particular $F$ is separably closed. 
Let $\overline{K}$ (resp. $\overline{R}$) be a separable closure of $K$ (resp. the 
normalization of $R$ in $\overline{K}$). The residue field $\overline{F}$ of 
$\overline{R}$ is an algebraic closure of $F$. 
Similarly as (\ref{ijm}), we also consider 
\begin{equation*}
\overline{X}\xrightarrow{\bar{j}}\overline{\mathcal{X}}\xleftarrow{\bar{i}}\mathcal{X}_{\overline{F}}, 
\end{equation*}
where $\overline{X}:=X\times_K\overline{K}$, $\overline{\mathcal{X}}:=
\mathcal{X}\times_R\overline{R}$, and $\mathcal{X}_{\overline{F}}:=
\mathcal{X}\times_R\overline{F}$. 
We take and fix an injection $K_\infty\to\overline{K}$ of extensions of $K$. 
Then, they fit into the commutative diagram 
\begin{equation*}
\xymatrix{
\overline{X}\ar[d]_f\ar[r]^{\bar{j}}&\overline{\mathcal{X}}\ar[d]_{\bar{f}}&
\mathcal{X}_{\overline{F}}\ar[l]_{\bar{i}}\ar[d]_{f_F}\\
X_\infty\ar[r]^{j_\infty}&\mathcal{X}_\infty&
\mathcal{X}_F\ar[l]_{i_\infty}. 
}
\end{equation*}
The two squares are cartesian if we replace $\mathcal{X}_{\overline{F}}$ by 
$\mathcal{X}\times_R(\overline{R}\otimes_{R_\infty}F)$, whose \'etale topos is 
canonically isomorphic to that of $\mathcal{X}_{\overline{F}}$. 
Let $I'={\rm Gal}(\overline{K}/K_\infty)$ be the Galois group of $\overline{K}/K_\infty$. 
Note that the functor $\Gamma(I',-)$ on discrete $\Lambda[I']$-modules is exact, since 
all the finite quotients of $I'$ are of order 
prime to $\ell$. 

1. 
By the local acyclicity, the canonical map  
$f_F^\ast i^\ast C\to \bar{i}^\ast R\bar{j}_\ast\bar{j}^\ast C$ is an isomorphism. 
Taking the fixed part $R\Gamma(I',-)=\Gamma(I',-)$, we have an isomorphism 
$\Gamma(I',f_F^\ast i^\ast C)\to \Gamma(I',\bar{i}^\ast R\bar{j}_\ast \bar{j}^\ast C)$. 
The source is isomorphic to $f_F^\ast i^\ast C$ since the action of $I'$ on $f_F^\ast i^\ast C$ is trivial. 
Since we have $\bar{i}^\ast R\bar{j}_\ast \bar{j}^\ast C\cong
f_F^\ast i_\infty^\ast Rj_{\infty\ast}f_\ast f^\ast j^\ast_\infty C$, 
the target is isomorphic to 
\begin{equation*}
f_F^\ast i_\infty^\ast Rj_{\infty\ast}\Gamma(I',
f_\ast f^\ast j_\infty^\ast C)\cong f_F^\ast i_\infty^\ast Rj_{\infty\ast}j^\ast_\infty C, 
\end{equation*}
hence the assertion 1. 

2. 
For an \'etale sheaf $\mathcal{G}$ of $\Lambda$-modules on $X_\infty$, we show that $R^nj_{\infty\ast}\mathcal{G}$ 
is zero for $n>2{\rm dim}X_\infty$. Let $x\to\mathcal{X}_{F}$ be a geometric point. 
We have an isomorphism 
\begin{equation*}
(Rj_{\infty\ast}\mathcal{G})_x\cong
\Gamma(I',R\Gamma(
(X_\infty\times_{\mathcal{X}_\infty}\mathcal{X}_{\infty(x)})\times_{K_{\infty}}
\overline{K},\mathcal{G})). 
\end{equation*}
Since ${\rm H}^n((X_\infty\times_{\mathcal{X}_\infty}\mathcal{X}_{\infty(x)})\times_{K_{\infty}}
\overline{K},\mathcal{G})$ is zero for $n>2{\rm dim}X_\infty$, the first assertion follows.

Let $C\in D^b_c(X,\Lambda)$. We have 
\begin{align*}
f_F^\ast i_\infty^\ast Rj_{\infty\ast}C|_{X_\infty}&\cong \Gamma(I', f_F^\ast 
i_\infty^\ast Rj_{\infty\ast}f_\ast f^\ast C|_{X_\infty})\\
&\cong \Gamma(I',\bar{i}^\ast R\bar{j}_\ast C|_{\overline{X}}). 
\end{align*}
The second assertion follows from the constructibility of the nearby cycles complex $\bar{i}^\ast R\bar{j}_\ast C|_{\overline{X}}$. Indeed, 
since $\Gamma(I',-)$ is exact, the cohomology sheaves 
 $\mathcal{H}^i(f_F^\ast i_\infty^\ast Rj_{\infty\ast}C|_{X_\infty})$ are 
 subsheaves of $\mathcal{H}^i(\bar{i}^\ast R\bar{j}_\ast C|_{\overline{X}})$. 
\qed}

Let $E/\mathbb{Q}_\ell$ be a finite extension with the ring of integers ${\cal O}_E$. 
\begin{cor}\label{constrcor}
Let the notation be as in Lemma \ref{ulainf}. Assume that $\mathcal{X}$ is of finite type over $R$. 
We use the notion and notation in Section \ref{adical}. 
Let $\mathcal{F}\in D(X^{\mathbb{N}^{\rm op}},{\cal O}_{E,\bullet})$ 
be a normalized constructible complex on $X$. 
Then the complex 
$i_\infty^\ast Rj_{\infty\ast}\mathcal{F}\in D(\mathcal{X}_{F}^{\mathbb{N}^{\rm op}},
{\cal O}_{E,\bullet})$ is a normalized ${\cal O}_E$-complex and 
$i_\infty^\ast Rj_{\infty\ast}\mathcal{F}\otimes^L_{{\cal O}_E}{\cal O}_E/\ell{\cal O}_E$ is constructible (where the functor 
$-\otimes_{{\cal O}_E}{\cal O}_E/\ell{\cal O}_E$ is defined in Definition \ref{alpha}.2). 
Namely,  $i_\infty^\ast Rj_{\infty\ast}\mathcal{F}$ defines a constructible complex of ${\cal O}_E$-sheaves on $\mathcal{X}_{F}$, 
in the sense of Definition 
\ref{consadic}.2. 
\end{cor}
We denote this complex $i_\infty^\ast Rj_{\infty\ast}\mathcal{F}$ by 
$\langle\mathcal{F},-\pi\rangle$. 
\proof{
Let $\mathcal{F}_{n}:=\mathcal{F}\otimes^L_{{\cal O}_E}{\cal O}_E/\ell^{n+1}{\cal O}_E$. By the first assertion of  Lemma \ref{ulainf}.2, 
We have $Rj_{\infty\ast}\mathcal{F}_{n+1}\otimes^L_{{\cal O}_E/\ell^{n+2}{\cal O}_E}{\cal O}_E/\ell^{n+1}{\cal O}_E
\cong Rj_{\infty\ast}(\mathcal{F}_{n+1}\otimes^L_{{\cal O}_E/\ell^{n+2}{\cal O}_E}{\cal O}_E/\ell^{n+1}{\cal O}_E)
\cong
Rj_{\infty\ast}\mathcal{F}_{n}$. Hence $i_\infty^\ast Rj_{\infty\ast}\mathcal{F}$ is 
normalized.  
By the second assertion of Lemma \ref{ulainf}.2, the complex $i^\ast_\infty Rj_{\infty\ast}\mathcal{F}\otimes^L_{{\cal O}_E}{\cal O}_E/\ell{\cal O}_E
\cong i^\ast_\infty Rj_{\infty\ast}\mathcal{F}_{0}$ is constructible. 
\qed}

\begin{lm}\label{mendoi}
Let the notation be as in Lemma \ref{ulainf}. 
Let $r\geq0$ be an integer. 
Let $L,N\in D^+(\mathcal{X}_r,\Lambda)$ be bounded below complexes. Assume 
that $i^\ast_rN$ is bounded and constructible and 
that the structure morphism 
$\mathcal{X}_r\to{\rm Spec}(R_r)$ is locally  acyclic 
relatively to $L$. 
Let 
\begin{equation}\label{zxc}
 L|_{X_r}\to M'\to  N|_{X_r}\to
\end{equation}
 be a distinguished triangle on $X_r$. 
Then, for some integer $n\geq r$, there exists a distinguished 
triangle  $L|_{\mathcal{X}_n}\to M\to  N|_{\mathcal{X}_n}\to$ on $\mathcal{X}_n$ whose pull-back to $X_n$ is isomorphic to that of (\ref{zxc}). 
\end{lm}
\proof{
Let $ \phi\colon N|_{X_r}\to L|_{X_r}[1]$ 
be the morphism corresponding to (\ref{zxc}). 
Let $C_n:=i_{n\ast}i_n^!L|_{\mathcal{X}_n}[2]$ be the complex on $\mathcal{X}_n$. It fits into the distinguished triangle 
$L|_{\mathcal{X}_n}[1]\to Rj_{n\ast}L|_{X_n}[1]\to C_n\to$. 
We need to show that, for some $n\geq r$, the composition 
$N|_{\mathcal{X}_n}\to Rj_{n\ast}N|_{X_n}\xrightarrow{Rj_{n\ast}\phi}
Rj_{n\ast}L|_{X_n}[1]\to C_n$ is zero. Since $C_n$ is supported on 
$\mathcal{X}_F$, it is enough to show that the 
restriction $i_r^\ast N\cong i_n^\ast N|_{\mathcal{X}_n}\to i^\ast_nC_n$ is zero. 
Since $\mathcal{X}_r\to{\rm Spec}(R_r)$ is locally acyclic relatively to $L$, 
the colimit $\varinjlim_{n\geq r}i^\ast_nC_n$ is acyclic by applying Lemma \ref{ulainf}.1 to $C=L$. Since 
$i^\ast_rN$ is constructible, the composition is zero 
for a large $n$. 
\qed
}

Let $S$ be a regular connected scheme of finite type over 
$\mathbb{Z}[1/\ell]$. 
Let $k$ be the perfection of the function field of $S$. 
Let $s\in S$ be a closed point. Let $S'$ be the 
blow-up of $S$ at $s$ and let $\mathfrak{s}$ denote  the generic point of the 
exceptional divisor. Let $R$ be the henselization of 
$\mathcal{O}_{S',\mathfrak{s}}$. Fix a uniformizer $\pi\in R$ of $R$. We follow 
the notation given above. That is, 
for an integer $m\geq0$, we write $R_m$ for 
$ R[\pi^{1/\ell^m}]$. We write $K_m$ for the fraction field of $R_m$. 
We write $R_\infty:=\varinjlim_mR_m$ and $K_\infty:=\varinjlim_mK_m$. 
The rings 
$R_m$ are valuation rings whose residue fields are 
isomorphic to $k(\mathfrak{s})$. 
\begin{lm}\label{chebdense}
The conjugates of the images of $G_{K_\infty}\to G_k$, for 
all the closed points $s\in S$ and uniformizers $\pi\in R$, 
topologically generate $G_k$. 
\end{lm}
\proof{
Let $H$ be a finite quotient of $G_k$. 
After shrinking $S$, the quotient map $G_k\to H$ factors through 
$\pi_1(S)$ and $H$ is generated by the geometric Frobeniuses at 
closed points $s\in S$ by the Chebotarev density theorem. Since the composition 
$G_{K_\infty}\to G_k\to\pi_1(S)$ factors as $G_{K_\infty}\to G_{k(\mathfrak{s})}
\to\pi_1(s)\to\pi_1(S)$ and the  map $G_{K_\infty}\to \pi_1(s)$ is surjective, 
the assertion follows.   
\qed
}

Consider a commutative diagram 
\begin{equation}\label{qwqwqw}
\xymatrix{
Z\ar@{^{(}-_>}[r]&U\ar[rd]_g\ar[rr]^{f}&&Y\ar[ld] \ar[r]^t&\mathbb{A}^1_k
\ar[lld]\\
&&{\rm Spec}(k)
}
\end{equation}
of $k$-schemes of finite type and a constructible complex $\mathcal{F}\in D^b_c(U,{\cal O}_E)$ 
with the following properties: 
\begin{enumerate}
\item The morphism $t\colon Y\to\mathbb{A}^1_k$ is  \'etale.  $U$ is smooth over $k$. 
$Z$ is a closed subscheme of $U$ 
finite \'etale over $k$. 
\item $f|_{U\setminus Z}$ is $SS(\mathcal{F})$-transversal, in the sense of 
Definition \ref{SSCCadic}. 
\end{enumerate}
Assume that the data given above except $\mathcal{F}$ are defined over $S$.  
Namely, we have a commutative diagram 
\begin{equation}\label{wq}
\xymatrix{
\mathcal{Z}\ar@{^{(}-_>}[r]&\mathcal{U}\ar[rd]_{\tilde{g}}
\ar[rr]^{\tilde{f}}&&\mathcal{Y}\ar[ld]\ar[r]^{\tilde{t}}&\mathbb{A}^1_S
\ar[lld] \\
&&S
}
\end{equation}
of $S$-schemes of finite type whose base change by 
${\rm Spec}(k)\to S$ is isomorphic to (\ref{qwqwqw}). 
We also assume that there exists $\Tilde{\mathcal{F}}_{0}\in
D_{\rm ctf}(\mathcal{U},{\cal O}_E/\ell{\cal O}_E)$ whose base change to $U$ 
is isomorphic to $\mathcal{F}\otimes^L_{{\cal O}_E}{\cal O}_E/\ell{\cal O}_E$. We further assume that 
they satisfy the following properties. 
\begin{enumerate}
\item The morphism $\tilde{t}\colon\mathcal{Y}\to\mathbb{A}^1_S$ is \'etale. 
$\mathcal{U}$ is smooth over $S$. $\mathcal{Z}$ is a closed subscheme of $\mathcal{U}$ finite \'etale over $S$. 
\item The relative singular support $SS(\Tilde{\mathcal{F}}_{0},\mathcal{U}/S)$ exists and satisfies the condition $2$ in Theorem \ref{relSS}. 
In particular, $\tilde{g}$ is universally locally acyclic relatively to 
$\Tilde{\mathcal{F}}_{0}$ (cf. Remark \ref{rmkss}). 
\item $\tilde{f}|_{\mathcal{U}\setminus \mathcal{Z}}$ is $SS(\Tilde{\mathcal{F}}_{0},\mathcal{U}/S)$-transversal. 
\item The restriction of the vanishing cycles complex $R\Phi_{\tilde{t}\circ\tilde{f}}(\Tilde{\mathcal{F}}_{0})$ 
to $\mathcal{Z}\overset\gets\times_{\mathbb{A}^1_S}(\mathbb{A}^1_S\setminus\tilde{t}\circ\tilde{f}
(\mathcal{Z}))\subset\mathcal{Z}\overset\gets\times_{\mathbb{A}^1_{S}}\mathbb{A}^1_{S}\cong\mathcal{Z}\overset\gets\times_{\mathbb{A}^1_{\mathcal{Z}}}\mathbb{A}^1_{\mathcal{Z}}$ is locally constant. 
For each $i\in\mathbb{Z}$, the function on $\mathcal{Z}$ defined by 
$z\mapsto{\rm dimtot}R^i\Phi_{\tilde{t}\circ\tilde{f}}(\Tilde{\mathcal{F}}_{0}
\otimes^L_{{\cal O}_E/\ell{\cal O}_E}\mathbb{F})
|_{\eta_{\bar{z}}}$ is locally constant (cf.  Proposition \ref{F}.3). Here $\mathbb{F}$ denotes the residue field of ${\cal O}_E$ and
 $\bar{z}$ is the spectrum of an algebraic 
closure of $k(z)$ and 
$\eta_{\bar{z}}$ is the generic point of the strict henselization of 
$\mathbb{A}^1_{k(\bar{z})}$ at $\bar{z}\xrightarrow{\tilde{t}\circ\tilde{f}}\mathbb{A}^1_{k(\bar{z})}$. 
\end{enumerate}

We denote by 
\begin{equation*}
\xymatrix{
\mathcal{Z}_\mathfrak{s}\ar@{^{(}-_>}[r]&\mathcal{U}_\mathfrak{s}\ar[r]^{\tilde{f}_\mathfrak{s}}&\mathcal{Y}_\mathfrak{s}
}
\end{equation*}
the base change of $\mathcal{Z}\hookrightarrow\mathcal{U}\to\mathcal{Y}$ by $\mathfrak{s}\to S$. 
Let $m$ be an integer $\geq0$ or $\infty$. We define $U_m\to\mathcal{U}_m\leftarrow\mathcal{U}_\mathfrak{s}$ 
to be the base change of $\mathcal{U}$ by ${\rm Spec}(K_m)\to{\rm Spec}(R_m)\leftarrow\mathfrak{s}$ 
over $S$. 
\begin{pr}\label{taihen}
Let the notation be as above. For a closed point $s\in S$, we write $S'$ for the blow-up of $S$ at $s$ and $\mathfrak{s}$ for the generic point of the exceptional divisor. We also write $R$ for the henselization of ${\cal O}_{S',\mathfrak{s}}$ and $K_\infty$ for the fraction field of 
$\cup_mR[\pi^{1/\ell^m}]$ where $\pi$ is a fixed uniformizer of $R$.

We have a commutative diagram 
\begin{equation}\label{commepmod}
\xymatrix{
G_k\ar[r]&
\Zl^\times\\G_{K_\infty}\ar[u]\ar[r]&G_{k(\mathfrak{s})}\ar[u]
}
\end{equation}
if $k$ is of positive characteristic. Here the top horizontal arrow is given by 
\begin{equation*}
\prod_{z\in Z}\varepsilon_0
(Y_{(z)},R\Phi_f(\mathcal{F})_z,dt)\circ {\rm tr}_{z/k}
\end{equation*}
 and the right vertical arrow 
is given by $\prod_{z\in 
\mathcal{Z}_\mathfrak{s}}\varepsilon_0
(\mathcal{Y}_{\mathfrak{s},(z)},R\Phi_{\tilde{f}_\mathfrak{s}}(\langle\mathcal{F},-\pi\rangle)_z,d
\tilde{t})\circ {\rm tr}_{z/\mathfrak{s}}$. 
The definition of 
$\langle\mathcal{F},-\pi\rangle$ is given after the statement of Corollary \ref{constrcor}. 

When $k$ is of characteristic $0$, 
we have a similar commutative diagram as (\ref{commepmod}) after replacing 
$\varepsilon_0$ and $\Zl^\times$ 
by $\overline{\varepsilon}_0$ (Definition \ref{epchar00}) and $\Zl^\times/\mu$. 
\end{pr}
\proof{
In the course of the proof, we use the notion of oriented products and the local Fourier transforms in the 
relative settings, which we briefly recall in Subsection \ref{recOLF}. 

Replacing $S$ and $\mathcal{U}$ by $\mathcal{Z}$ and an open 
neighborhood of the graph $\mathcal{Z}\hookrightarrow\mathcal{Z}\times_S\mathcal{U}$, 
we may assume that $\mathcal{Z}\to S$ is an isomorphism. 
We may also replace 
$(\mathcal{Y},\tilde{t})$ by $(\mathbb{A}^1_S,{\rm id})$ as the local epsilon factors in the statement do not change.  We take $S\cong\mathcal{Z}\xrightarrow{\tilde{f}|_{\mathcal{Z}}}
\mathbb{A}^1_S$ as the origin.  
By induction on $m\geq0$, we find a sequence of integers $0=n_0\leq n_1\leq n_2\dots$ and 
construct $
\Tilde{\mathcal{F}}_{m}\in D_{\rm ctf}(\mathcal{U}_{n_m},{\cal O}_E/\ell^{m+1}{\cal O}_E)$ with 
$\Tilde{\mathcal{F}}_{m+1}\otimes^L_{{\cal O}_E/\ell^{m+2}{\cal O}_E}{\cal O}_E/\ell^{m+1}{\cal O}_E\cong\Tilde{\mathcal{F}}_{m}|_{{\cal U}_{n_{m+1}}}$. 
Suppose that $n_m$ and $\Tilde{\mathcal{F}}_{m}$ be given.  Applying Lemma \ref{mendoi} to 
$L=\Tilde{\mathcal{F}}_{m}$ and $N=\Tilde{\mathcal{F}}_{0}$, 
we find a constructible complex $\Tilde{\mathcal{F}}_{m+1}$ of 
${\cal O}_E/\ell^{m+1}{\cal O}_E$-sheaves on  
$\mathcal{U}_{n}$ for some integer $n\geq n_m$ 
which fits into a distinguished triangle 
\begin{equation*}
\Tilde{\mathcal{F}}_{0}\to\Tilde{\mathcal{F}}_{m+1}\to\Tilde{\mathcal{F}}_{m}\to
\end{equation*}
whose restriction to $U_n$ is isomorphic to the pull back of 
\begin{equation*}
\mathcal{F}_{0}\to\mathcal{F}_{m+1}\to\mathcal{F}_{m}\to. 
\end{equation*} 
We show that $\Tilde{{\cal F}}_{m+1}$ and $n_{m+1}=n$ have the desired properties. 
By induction, ${\cal U}_n\to{\rm Spec}(R_n)$ is universally locally acyclic relatively to $\Tilde{{\cal F}}_{m},\Tilde{{\cal F}}_{0}$, hence to $\Tilde{{\cal F}}_{m+1}$. 
We claim that the canonical morphism $\tilde{\phi}\colon\Tilde{\mathcal{F}}_{m+1}\otimes^L_{{\cal O}_E/\ell^{m+2}{\cal O}_E}
{\cal O}_E/\ell^{m+1}{\cal O}_E\to
\Tilde{\mathcal{F}}_{m}$ is a quasi-isomorphism, which also shows that $\Tilde{{\cal F}}_{m+1}\in D_{\rm ctf}(\mathcal{U}_{n},{\cal O}_E/\ell^{m+1}{\cal O}_E)$. 
Indeed, the restriction $\tilde{\phi}|_{U_n}$ is isomorphic to the canonical one 
$\phi\colon\mathcal{F}_{m+1}\otimes^L_{\Zl/\ell^{m+2}{\cal O}_E}{\cal O}_E/\ell^{m+1}{\cal O}_E\to
\mathcal{F}_{m}$, which is a quasi-isomorphism. On the other hand, by Lemma \ref{ulainf}.1 and the local acyclicity 
of ${\cal U}_n$ over $R_n$ relative to $\Tilde{{\cal F}}_{i}$ for $i=m,m+1$,  
the canonical map 
\begin{equation*}
\Tilde{\mathcal{F}}_{i}|_{\mathcal{U}_{\mathfrak{s}}}\to
i^\ast_\infty Rj_{\infty\ast}\mathcal{F}_{i}=:\langle\mathcal{F}_{i},-\pi\rangle. 
\end{equation*}
is a quasi-isomorphism. 
Therefore, 
the restriction $i_n^\ast\tilde{\phi}$ to 
the special fiber is identified with 
$
i_\infty^\ast Rj_{\infty\ast}\mathcal{F}_{m+1}\otimes_{{\cal O}_E/\ell^{m+2}{\cal O}_E}^L{\cal O}_E/\ell^{m+1}{\cal O}_E
\to i_\infty^\ast Rj_{\infty\ast}\mathcal{F}_{m}$, which further can be identified with 
$i_\infty^\ast Rj_{\infty\ast}\phi$ since the cohomological dimension of $Rj_{\infty\ast}$ is finite. 
Hence we can put $n_{m+1}=n$.  

Take one $m$ and put $n=n_m$. 
Let $\tilde{f}_n\colon \mathcal{U}_n\to\mathbb{A}^1_{R_n}$ be the base change of $\tilde{f}$.  
The restrictions of $R\Phi_{\tilde{f}_n}\Tilde{\mathcal{F}}_{m}$ (which is regarded as a complex on 
$0_{R_n}\overset\gets\times_{\mathbb{A}^1_{R_n}}\mathbb{G}_{m,R_n}$) 
to $0_{\mathfrak{s}}\overset\gets\times_{\mathbb{A}^1_{\mathfrak{s}}}\mathbb{G}_{
m,\mathfrak{s}}$ and $0_{K_n}\overset\gets\times_{\mathbb{A}^1_{K_n}}\mathbb{G}_{
m,K_n}$ are isomorphic to $R\Phi_{\tilde{f}_{\mathfrak{s}}}\langle\mathcal{F}_{m},
-\pi\rangle$ and $R\Phi_{f}\mathcal{F}_{m}$ respectively by the commutativity of the formation of $R\Phi_f$  \cite[Proposition 6.1]{Org}.  
By the assumption $4$,  
$R\Phi_{\tilde{f}_n}\Tilde{\mathcal{F}}_{m}$ 
is locally constant on $0_{R_n}\overset\gets\times_{\mathbb{A}^1_{R_n}}\mathbb{G}_{{ m},R_n}$  and its total dimension is locally constant. 
In particular, if the generic fiber $R\Phi_{f}\mathcal{F}_{m}$ is tamely ramified, 
so is $R\Phi_{\tilde{f}_{\mathfrak{s}}}\langle\mathcal{F}_{m},
-\pi\rangle$. 

We show the assertion for the case when $k$ is of positive characteristic. 
By Proposition \ref{F}, $F^{(0,\infty)}(R\Phi_{\tilde{f}}\Tilde{\mathcal{F}}_{m})$ is locally 
constant and its restrictions to $\infty_{\mathfrak{s}}\overset\gets\times_{\mathbb{P}^1_{\mathfrak{s}}}\mathbb{A}^1_{\mathfrak{s}}$ and $\infty_{K_n}\overset\gets\times_{\mathbb{P}^1_{K_n}}\mathbb{A}^1_{K_n}$ are isomorphic to $F^{(0,\infty)}(R\Phi_{\tilde{f}_\mathfrak{s}}\langle\mathcal{F}_{m},
-\pi\rangle)$ and $F^{(0,\infty)}(R\Phi_{f}\mathcal{F}_{m})$ respectively. 
By \cite[Corollary 3.9]{CLE}, the determinant $\det F^{(0,\infty)}(R\Phi_{\tilde{f}}\Tilde{\mathcal{F}}_{m})$ is tamely ramified. 
Applying the construction in  \cite[Definition 2.16]{CLE}, we get a character $\langle\det F^{(0,\infty)}(R\Phi_{\tilde{f}}\Tilde{\mathcal{F}}_{m}),1/x\rangle$ 
on $R_n$ (with the notation being in loc.~cit.), where $x$ is the standard coordinate on $\mathbb{A}^1\subset\Proj^1$. 
The assertion follows from \cite[Lemma 2.19]{CLE} and Laumon's cohomological interpretation \cite[Theorem 3.10]{CLE} in this case. 

The assertion for the case of characteristic $0$ is proved as follows. 
Let $x$ be the standard coordinate of $\mathbb{A}^1$.  
It suffices to show the commutativity of the diagram 
\begin{equation*}
\xymatrix{
G_k\ar[r]&
\Zl^\times/\mu\\G_{K_\infty}\ar[u]\ar[r]&G_{k(\mathfrak{s})},\ar[u]
}
\end{equation*}
where the top horizontal arrow is given by  $J(R\Phi_f\mathcal{F})$ and 
the right vertical arrow is 
given by $J(R\Phi_{\tilde{f}_\mathfrak{s}}\langle\mathcal{F},-\pi\rangle)$. 

Fix a separable closure $\overline{K}_\infty$ of $K_\infty$. 
Define $I:=\varprojlim_N\mu_N(\overline{K}_\infty)$ and 
$I^p:=\varprojlim_{p\nmid N}\mu_N(\overline{K}_\infty)$ where 
$p$ is the residue characteristic of $R$. 
Let $\eta_{k(\mathfrak{s})}$ (resp. $\eta_{K_\infty}$
) be the function field of 
the henselization $\mathbb{A}^1_{k(\mathfrak{s}),(0)}$ (resp. 
$\mathbb{A}^1_{K_\infty,(0)}$). Fix geometric points $\overline{\eta}_\mathfrak{s}$ and 
$\overline{\eta}_\infty$ over $\eta_{k(\mathfrak{s})}$ and $\eta_{K_\infty}$, 
and also fix a specialization $\overline{\eta}_\infty\to \overline{\eta}_\mathfrak{s}$ as 
geometric points of the henselization $\mathbb{A}^1_{R_n,(0_\mathfrak{s})}$. 
The group $I^p$ (resp. $I$) can be naturally regarded as 
the tame inertia group $I^t_{\eta_{k(\mathfrak{s})}}$ (resp. $I^t_{\eta_{K_\infty}}$). 
Let $\pi_1^t$ denote the fundamental group classifying finite \'etale coverings of 
$\mathbb{A}^1_{R_n,(0_\mathfrak{s})}\setminus 0_{R_n,(0_\mathfrak{s})}$ tamely ramified 
along $0_{R_n,(0_\mathfrak{s})}$. 
 We also have a natural embedding $I^p\hookrightarrow \pi^t_1$ and a commutative diagram 
 \begin{equation*}
 \xymatrix{
 I\ar[d]^\cong\ar[r]&I^p\ar[d]\ar[dr]^\cong&\\
 I^t_{\eta_{K_\infty}}\ar[r]&\pi^t_1&I^t_{\eta_{k(\mathfrak{s})}},\ar[l]
 }
 \end{equation*}
where the top horizontal arrow is the projection $I\to I^p$ and 
the bottom horizontal arrows are the canonical ones. 
 
Since the generic characteristic of $R_n$ is zero, 
the restriction of $R\Phi_{\tilde{f}_n}\Tilde{\mathcal{F}}_{m}$ to 
$\mathbb{A}^1_{R_n,(0_\mathfrak{s})}\setminus 0_{R_n,(0_\mathfrak{s})}
\cong0_{\mathfrak{s}}\overset\gets\times_{\mathbb{A}^1_{R_n}}
\mathbb{G}_{{ m},R_n}\subset0_{R_n}\overset\gets\times_{\mathbb{A}^1_{R_n}}
\mathbb{G}_{{ m},R_n}$ is locally constant with tamely ramified cohomology sheaves. 
Since it is locally constant, the specialization 
$(R\Phi_{\tilde{f}_\mathfrak{s}}\langle\mathcal{F}_{m},-\pi\rangle)_{\overline{\eta}_\mathfrak{s}}
\cong (R\Phi_{\tilde{f}_n}\Tilde{\mathcal{F}}_{m})_{\overline{\eta}_\mathfrak{s}}\to 
(R\Phi_{\tilde{f}_n}\Tilde{\mathcal{F}}_{m})_{\overline{\eta}_\infty}\cong 
(R\Phi_f\mathcal{F}_{m})_{\overline{\eta}_\infty}$ is an isomorphism, which we regard as 
an isomorphism of complexes of $I$-representations. These isomorphisms commute with the transition maps 
$R\Phi_f\mathcal{F}_{m+1}\to R\Phi_f\mathcal{F}_{m}$ and 
$R\Phi_{\tilde{f}_\mathfrak{s}}\langle\mathcal{F}_{m+1},-\pi\rangle\to R\Phi_{\tilde{f}_\mathfrak{s}}\langle\mathcal{F}_{m},-\pi\rangle$. 
Hence $C_\infty:=R\Phi_f\mathcal{F}\otimes_{{\cal O}_E} E$ and $C_\mathfrak{s}:=R\Phi_{\tilde{f}_\mathfrak{s}}\langle\mathcal{F},-\pi\rangle\otimes_{{\cal O}_E}E$ 
are isomorphic as complexes of $I$-representations on finite dimensional $E$-vector spaces. 
Let $I\to\mu_N(\overline{K}_\infty)$ be a finite quotient through which $I$ acts on the semi-simplifications of cohomologies of 
$C_\infty\cong C_\mathfrak{s}$. 
We may assume that $p\nmid N$. For each $i$, 
the semi-simplification of ${ H}^i(C_\infty)^{\oplus N}$ 
gives a Jacobi datum on $K_\infty$ as in Definition \ref{js}, which extends to a Jacobi datum on 
$R_\infty$ since $N$ is invertible in $R_\infty$. 
Since this Jacobi datum on $R_\infty$ is restricted to the one on $k(\mathfrak{s})$ 
constructed from ${ H}^i(C_\mathfrak{s})^{\oplus N}$ and 
the one on $K_\infty$ constructed from 
$H^i(C_\infty)^{\oplus N}$
, the assertion follows. 
\qed
}

\begin{rmk}
Using a similar method as above repeatedly, one can reduce several 
problems on $\ell$-adic sheaves on schemes of finite type over (the perfections of) finitely 
generated fields to cases over finite fields. For example, 
Theorem $1$ in \cite{Jacob} can be proven unconditionally, i.e., 
without the assumption that the sheaf $\mathcal{F}$ in loc.~cit. is defined 
over a scheme of finite type over $\mathbb{Z}$, if the function field of 
the base scheme $S$ is a purely inseparable extension of a finitely 
generated field, although it would be desireble to prove it by 
developing a theory of Jacobi sum characters for representations 
with torsion coefficient.  
\end{rmk}

\subsection{Local epsilon factors of convolutions}\label{locefconv} 
At the end of this section, we compute the local epsilon factors of the 
convolutions of vanishing cycles. To do so, we need to recall the 
Thom--Sebastiani theorem for \'etale sheaves, which is  proved in \cite{Ill}. For basics of oriented products which are used in this subsection, we refer to \cite{ori} and \cite{Ill}. 

Let $k$ be a perfect field of characteristic $p\geq0$. We fix a prime number $\ell$ different from $p$. Let $\Lambda$ be a finite local ring with residue characteristic $\ell$. 
Following the notation in \cite{Ill}, let $A_h,A^2_h$ denote the henselizations of $\mathbb{A}^1_{k}, 
\mathbb{A}_{k}^2$ at $0,(0,0)$ respectively. Let $a\colon A^2_h\to A_h$ be the map induced from the summation $\mathbb{A}^2_{k}
\to\mathbb{A}_{k}$. 
Let $f_1\colon X_1\to A_h$ and $f_2\colon X_2\to A_h$ be two morphisms 
of schemes of finite type. 
Let $X:=(X_1\times X_2)\times_{A_h\times A_h}A^2_h$ and  $f\colon X\to 
A^2_h$ be 
the projection. We regard $X$ as an $A_h$-scheme by  
the composition $X\xrightarrow{f} A^2_h\xrightarrow{a}A_h$. 

\begin{df-lm}(\cite[Definition 4.1]{Ill}, \cite[Proposition 4.3]{Ill})
For each $i=1,2$, let $K_i$ be an object of $D_{\rm ctf}(
X_i\overset\gets\times_{A_h}A_h,\Lambda)$. We define the local convolution ${K}_1\ast^L_\Lambda
K_2\in D(X\overset\gets\times_{A_h}A_h,\Lambda)
$ of $K_1$ and $K_2$ by setting 
\begin{equation*}
K_1\ast^L_\Lambda K_2:=R\overset\gets a_\ast
(\overset\gets {\rm pr}_{1\ast}K_1\otimes^L_\Lambda\overset\gets {\rm pr}_{2\ast}
K_2)[1],
\end{equation*}
where $\overset\gets {\rm pr}_{i}\colon X\overset\gets\times_{A^2_h}A^2_h\to
X_i\overset\gets\times_{A_h}A_h$
 and $\overset\gets a\colon X\overset\gets\times_{A^2_h}A^2_h\to
X\overset\gets\times_{A_h}A_h$ are induced from the $i$-th projections and $a$. If no confusions occur, we omit $\Lambda$ and write $\ast^L$ for $\ast^L_\Lambda$. 
The complex $K_1\ast^LK_2$ belongs to $D_{\rm ctf}(X\overset\gets\times
_{A_h}A_h,\Lambda)$ (\cite[Proposition 4.3]{Ill}). 
\end{df-lm}

We remark that this definition is slightly different from that in \cite{Ill}, 
as the complex is shifted by $1$. 

We show that the formation of $K_1\ast^LK_2$ commutes with base changes on $X_i$, which is used in the proof of Lemma \ref{epconv}. 
\begin{lm}\label{astbc}
Let $g_i\colon Y_i\to X_i$ be morphisms of $A_h$-schemes of finite type for $i=1,2$. Put $Y:=Y_1\times Y_2$ and let 
$g\colon Y\to X$ be the product of $g_i$. Let 
$\overset\gets g_i\colon Y_i\overset\gets\times
_{A_h}A_h\to X_i\overset\gets\times
_{A_h}A_h$ and $\overset\gets g\colon Y\overset\gets\times
_{A_h}A_h\to X\overset\gets\times
_{A_h}A_h$ be the induced morphisms. Then, 
for $K_i\in D_{\rm ctf}(
X_i\overset\gets\times_{A_h}A_h,\Lambda)$, 
the canonical map 
\begin{equation*}
\overset\gets g^\ast(K_1\ast^L K_2)\to (\overset\gets g_1^\ast K_1)\ast^L(\overset\gets g_2^\ast K_2)
\end{equation*}
in $D_{\rm ctf}(
Y\overset\gets\times_{A_h}A_h,\Lambda)$ is an isomorphism. 
\end{lm}
\proof{
The assertion follows if one verifies that the base change map $\overset\gets g^\ast R\overset\gets a_\ast\to R\overset\gets a_\ast \overset\gets g^\ast$ arising from the commutative diagram 
\begin{equation*}
\xymatrix{Y\overset\gets\times_{A_h^2}A_h^2\ar[r]^-{\overset\gets a}\ar[d]^-{\overset\gets g}&Y\overset\gets\times_{A_h}A_h\ar[d]^-{\overset\gets g}\\
X\overset\gets\times_{A_h^2}A_h^2\ar[r]^-{\overset\gets a}&X\overset\gets\times_{A_h}A_h
}
\end{equation*}
is an isomorphism. This is proved in \cite[Corollary 2.5]{CLE}. 
\qed}

The following is the Thom--Sebastiani theorem for \'etale sheaves obtained in \cite{Ill}. 
\begin{thm}(\cite[Theorem 4.5]{Ill})\label{conv}
With the notation being above, let $K_1,K_2$ be objects 
of $D_{\rm ctf}(X_1,\Lambda),D_{\rm ctf}(X_2,\Lambda)$ respectively. 
Let $K:=(K_1\boxtimes^LK_2)|_{X}$. Then, 
there is a functorial isomorphism 
\begin{equation*}
(R\Phi_{f_1}(K_1))\ast^L(R\Phi_{f_2}(K_2))|_{
X_0\overset\gets\times_{A_h}A_h}\cong 
R\Phi_{af}(K)[1]|_{X_0\overset
\gets\times_{A_h}A_h}
\end{equation*}
in $D_{\rm ctf}(X_0\overset\gets\times_{A_h}A_h,\Lambda)$, where $X_0$ is 
the closed fiber of $X\to A_h$. 
\end{thm}

If one takes $X_i$ to be the closed point of $A_h$, then $X_i\overset\gets\times_{A_h}A_h$ are  canonically equivalent to the \'etale topos of $A_h$. Then $\ast^L$ can be regarded as a functor 
\begin{equation*}
D_{\rm ctf}(A_h,\Lambda)\times
D_{\rm ctf}(A_h,\Lambda)\to D_{\rm ctf}(A_h,\Lambda). 
\end{equation*}

We consider the $\ell$-adic variant  
of this convolution functor on the derived category 
$D^b_c(A_h,E)$ for a finite extension $E$ of $\mathbb{Q}_\ell$, which is defined as follows. For notations appearing in the construction of $D^b_c(-,E)$, we refer to Section $6$. 

Let $K_1,K_2\in 
D_{\rm ctf}(A_h,\Lambda)$. Since the cohomological dimension of $R\overset\gets a_\ast$ is finite (this can be proved similarly as \cite[Proposition 3.1]{Org}, by using \cite[Proposition 1.13]{Ill}), the canonical map 
$(K_1\ast^L_\Lambda K_2)\otimes_\Lambda^L\Lambda'\to 
(K_1\otimes_\Lambda^L\Lambda')\ast^L_{\Lambda'}(K_2\otimes_\Lambda^L\Lambda')$ 
for a local ring homomorphism $\Lambda\to\Lambda'$
 between finite local rings is an isomorphism. 
 Therefore, the functor $D^b(A_h^{\mathbb{N}^{\rm op}},{\cal O}_{E\bullet})\times D^b(A_h^{\mathbb{N}^{\rm op}},{\cal O}_{E\bullet})\to D^b(A_h^{\mathbb{N}^{\rm op}},{\cal O}_{E\bullet})$ given by $(K_1,K_2)\mapsto 
 R\overset\gets a_\ast
(\overset\gets {\rm pr}_{1\ast}K_1\otimes^L_\Lambda\overset\gets {\rm pr}_{2\ast}
K_2)[1]$ preserves $D_{c,{\rm norm}}$ (Definition \ref{consadic}.1) and 
induces a functor 
\begin{equation*}
\ast^L\colon D^b_c(A_h,{\cal O}_E)\times D^b_c(A_h,{\cal O}_E)\to D^b_c(A_h,{\cal O}_E). 
\end{equation*}
  Taking $\varinjlim_E$ where $E$ runs through the finite extensions in $\Ql/\mathbb{Q}_\ell$ and  inverting $\ell$, we also get convolution functors on 
$D^b_c(A_h, \Zl)$ and  $D^b_c(A_h,\Ql)$. 
Let $\eta$ be the generic point of $A_h$. 
For ${\cal F}\in D^b_c(\eta,\Ql)$, let ${\cal F}_!$ denote  the $0$-extension of $\cal F$ to $A_h$. As explained in \cite[3.B]{Ill}, the functor $({\cal F}_1,{\cal F}_2)\mapsto (
{\cal F}_{1!}\ast^L{\cal F}_{2!})|_\eta$ is isomorphic to the convolution functor considered in \cite[(2.7)]{Lau}.

To state Lemma \ref{epconv}, we slightly change the notation. Let $k$ be a perfect field as above. Let $f_1\colon X_1\to\mathbb{A}^1_k$ 
and $f_2\colon X_2\to\mathbb{A}^1_k$ be $k$-morphisms of finite type. 
Let $X:=X_1\times_kX_2$ and let $af\colon X\to\mathbb{A}^1_k$ denote 
 the composition of the product 
$f:=f_1\times f_2\colon X_1\times_kX_2\to\mathbb{A}^2_k$ and 
the summation $a\colon\mathbb{A}^2_k\to\mathbb{A}^1_k$. 
For  $i=1,2$, let $\mathcal{F}_{i}\in D^b_c(X_i,\Zl)$ and let 
$x_i\in X_i$ be  at most isolated  $SS(\mathcal{F}_{i})$-characteristic 
$k$-rational points of $f_i$ with $f_i(x_i)=0$. 
\begin{lm}\label{epconv}
Let the notation be as above. Let $x:=(x_1,x_2)\in X$ be the $k$-rational 
point over $x_1$ and $x_2$. Let $\mathcal{F}:=\mathcal{F}_{1}\boxtimes^L\mathcal{F}_{2}$. 
Denote by $t$ the standard coordinate of $\mathbb{A}^1_k$. Let $A_h$ be the henselization of $\mathbb{A}^1_k$ at the origin. 
\begin{enumerate}
\item 
Assume that $p>0$. We have the equality 
\begin{align*}
\varepsilon_0(A_h,R\Phi_{af}&
(\mathcal{F})_x,dt)^{-1}=\\
&\varepsilon_0(A_h,R\Phi_{f_1}(\mathcal{F}_1)_{x_1},dt)^{\dt R\Phi_{f_2}
(\mathcal{F}_2)_{x_2}}\cdot\varepsilon_0(A_h,
R\Phi_{f_2}(\mathcal{F}_2)_{x_2},dt)^{\dt R\Phi_{f_1}(\mathcal{F}_1)_{x_1}}. 
\end{align*}
\item Assume that $k$ is finitely generated over $\mathbb{Q}$. 
We have 
 \begin{align*}
\overline{\varepsilon}_0(A_h,R\Phi_{af}&
(\mathcal{F})_x)^{-1}=\\
&\overline{\varepsilon}_0(A_h,R\Phi_{f_1}(\mathcal{F}_1)_{x_1})^{\dt R\Phi_{f_2}
(\mathcal{F}_2)_{x_2}}\cdot\overline{\varepsilon}_0(A_h,
R\Phi_{f_2}(\mathcal{F}_2)_{x_2})^{\dt R\Phi_{f_1}(\mathcal{F}_1)_{x_1}}. 
\end{align*}
Here $\overline{\varepsilon_0}$ is defined in Definition \ref{epchar00}. 
\end{enumerate}
\end{lm}
\proof{
1. By the $\Zl$-variant of Theorem \ref{conv}, we have an isomorphism 
\begin{equation*}
(R\Phi_{f_1}(\mathcal{F}_{1}))\ast^L(R\Phi_{f_2}(\mathcal{F}_{2}))|_{X_0\overset\gets\times_{A_h}A_h}
\cong R\Phi_{af}(\mathcal{F})[1]|_{X_0\overset\gets\times_{A_h}A_h}
\end{equation*}
of objects in $D^b_c({X_0\overset\gets\times_{A_h}A_h},\Zl).$ Shrinking $X_i$ around $x_i$, we suppose that $X_i\setminus\{x_i\}\to\mathbb{A}^1_k$ are universally locally acyclic relatively to ${\cal F}_{i}$, hence that $R\Phi_{f_1}(\mathcal{F}_{i})$ are supported on $x_i\overset\gets\times_{\mathbb{A}^1_k}\mathbb{A}^1_k\cong A_h$. Since $R\Phi_{f_1}(\mathcal{F}_{i})$ are isomorphic to $j_!R\Phi_{f_1}(\mathcal{F}_{i})_{x_i}$ where $j$ denotes the open immersion $\eta\to A_h$, Lemma \ref{astbc} gives an isomorphism 
\begin{equation*}
(R\Phi_{f_1}(\mathcal{F}_{1})_{x_1})\ast(R\Phi_{f_2}(\mathcal{F}_{2})_{x_2})\cong R\Phi_{af}(\mathcal{F})_x[1], 
\end{equation*}
where $\ast$ denotes the convolution functor in \cite[(2.7)]{Lau} shifted by $1$. 
By \cite[Proposition (2.7.2.2)]{Lau}, we have 
\begin{equation*}
F^{(0,\infty)}(R\Phi_{f_1}(\mathcal{F}_1)_{x_1})\otimes
F^{(0,\infty)}(R\Phi_{f_2}(\mathcal{F}_2)_{x_2})\cong
F^{(0,\infty)}(R\Phi_{af}(\mathcal{F})_{x})[1].
\end{equation*}
Using this isomorphism and \cite[Th\'eor\`eme (3.5.1.1)]{Lau}, the assertion follows. 

2. Apply Proposition \ref{taihen} 
to the commutative diagrams 
\begin{equation*}
\xymatrix{
{\rm Spec}(k)\ar@{^{(}-_>}[r]^{\ \ \ x_i}&X_i\ar[rd]\ar[rr]^{f_i}&&
\mathbb{A}^1_{k}\ar[r]^{{\rm id}}\ar[ld]&\mathbb{A}^1_k\ar[lld] \\
&&{\rm Spec}(k)
}
\end{equation*}
and $\mathcal{F}_{i}$ for $i=1,2$ and the similar diagram for $X$ and 
$\mathcal{F}$. 
Then, the assertion follows from 
1 and Lemma \ref{chebdense}. 
\qed
}

\section{Epsilon Cycles of $\ell$-adic Sheaves}\label{cec}

In this section, we construct epsilon cycles which compute local epsilon factors 
modulo roots of unity. Let $\mu$ denote the subgroup 
of roots of unity in $\Zl^\times$. 
\subsection{Group of characters modulo torsion}

For a field $E$, let $\mu_E$ be the group of roots of unity in $E$. 
\begin{df}\label{dfth}
Let $G$ be a compact Hausdorff abelian group. 
\begin{enumerate}
\item
For a finite extension $E$ of $\mathbb{Q}_\ell$, 
define $\Theta_{G,E}$ to be the group  
${\rm Hom}_{\rm conti}(G,\mathcal{O}_E^\times/\mu_E)$ 
of continuous homomorphisms. Here $\mathcal{O}_E^\times/\mu_E$ is endowed with the topology induced from the $\ell$-adic topology on ${\cal O}_E$. 
\item 
Define the group $\Theta_G$ by setting $\Theta_G:=\varinjlim_E\Theta_{G,E}$,
 where $E$ runs through the finite subextensions in $\Ql/\mathbb{Q}_\ell$. 
 \item
When $G$ is the abelianization of the absolute Galois group of a field $k$, 
$\Theta_{G,E}$ and $\Theta_G$ are also denoted by $\Theta_{k,E}$ and 
$\Theta_k$. 
\end{enumerate}
\end{df}
We usually  identify $\Theta_G$ with a subgroup  of the group 
${\rm Hom}(G,\Zl^\times/\mu)$ consisting of group 
homomorphisms. A group homomorphism $
G\to\Zl^\times/\mu$ is said to be {\it admissible} if 
it belongs to $\Theta_G$. By Lemma \ref{Baire} below, a compact subgroup of $\Zl^\times$ 
 is contained in 
$\mathcal{O}_E^\times$ for some finite subextension $E$ of $\Ql/\mathbb{Q}_\ell$, where $\Zl^\times$  is 
endowed with the topology induced from the valuation of $\Ql$. Therefore 
continuous homomorphisms 
$G\to\Zl^\times$ give admissible  homomorphisms $G\to\Zl^\times/\mu$. 
\begin{lm}\label{Baire}
Let $K\subset{\rm GL}_n(\Ql)$ be a compact subgroup. 
Then, there exists a finite subextension $E$ of $\Ql/\mathbb{Q}_\ell$ 
such that $K\subset{\rm GL}_n(E)$. 
\end{lm}
\proof{
We give a proof for completeness. 
Fix a bijection from the set of integers $\geq0$ to the set of 
finite subextensions of $\Ql/\mathbb{Q}_\ell$, which 
is denoted by $m\mapsto E_m$. For integers $m\geq0$, 
put $K_m:=K\cap{\rm GL}_n(E_m)$. They are closed subgroups 
of $K$ and cover the whole of $K$, i.e., $\cup_mK_m=K$. 
Since $K$ is compact Hausdorff, Baire category theorem 
can be applied. Hence, there exists $m\geq0$ such that 
$K_m$ contains a non-empty open subset of $K$, which 
implies that $K_m$ is an open subgroup. Then, the index $[K:K_m]$ is finite, which implies the assertion. 
\qed
}

\begin{lm}\label{theta}
Let $G$ be a compact Hausdorff abelian group. 
\begin{enumerate}
\item The group $\Theta_G$ is uniquely divisible. 
\item Let ${\rm Hom}_{\rm conti}(G,\Zl^\times)$ be the 
group of continuous group homomorphisms to $\Zl^\times$. Then, the kernel and the 
cokernel of the natural map ${\rm Hom}_{\rm conti}(G,
\Zl^\times)\to\Theta_G$ are torsion. 
\end{enumerate}
\end{lm}
\proof{
1. 
Since the group $\mathcal{O}_E^\times/\mu_E$ is torsion-free, so is 
$\Theta_{G,E}$. Hence $\Theta_G$ is torsion-free. 
Let $\chi\in\Theta_{G,E}$ be a continuous homomorphism. 
For an integer $n\geq1$, we need to find a finite extension $E'$ of $E$ 
and a continuous homomorphism $\xi\colon G\to
\mathcal{O}_{E'}^\times/\mu_{E'}$ so that $\xi^n=\chi$. 
Let $E'$ be a finite  extension  of $E$ which 
contains the $n$-th roots of all elements in 
$\mathcal{O}_E^\times$; such an  $E'$  exists since $\mathcal{O}_E^\times/
(\mathcal{O}_E^\times)^n$ is finite. 
Then, the composite of $\chi$ and the natural 
inclusion $\mathcal{O}_E^\times/\mu_E\to\mathcal{O}_{E'}^\times/\mu_{E'}$ 
factors through the injection  $\mathcal{O}_{E'}^\times/\mu_{E'}\to
\mathcal{O}_{E'}^\times/\mu_{E'}$ given by $a\mapsto a^n$. 
Since this injection is a homeomorphism onto the image, 
we find a desired  $\xi$. 

2. The kernel is torsion since compact subgroups of $\mu\subset\Zl^\times$ 
are finite subgroups by Lemma \ref{Baire}. 

Let $E$ be a finite extension of $\mathbb{Q}_\ell$ and $\chi\colon 
G\to \mathcal{O}_E^\times/\mu_E$ be a 
continuous homomorphism. We find a continuous homomorphism $\xi\colon 
G\to \mathcal{O}_E^\times$ and an integer $n\geq1$ 
such that the composite of $\xi$ and the quotient map $\mathcal{O}_E^\times\to 
\mathcal{O}_E^\times/\mu_E$ equals to $\chi^n$. 
Take an open subgroup $U\subset\mathcal{O}_E^\times$ such that 
$U\cap\mu_E$ is trivial. Then, the composition $U\to\mathcal{O}_E^\times\to
\mathcal{O}_E^\times/\mu_E$ is an isomorphism onto an open subgroup 
of $\mathcal{O}_E^\times/\mu_E$, which we also denote by $U$. 
Let $H$ be the inverse image of $U\subset\mathcal{O}_E^\times/\mu_E$ by 
$\chi$. This is an open subgroup of 
$G$. Let $n:=[G:H]$ be the 
index. Define $\xi$ by the composition $G
\xrightarrow{n}H\xrightarrow{\chi}
U\to\mathcal{O}_E^\times$. Then, $\xi$ and $n$ satisfy 
the condition. 
\qed
}

\subsection{Constructions of epsilon cycles}
In this subsection, we construct epsilon cycles for $\Zl$-sheaves (Theorem \ref{epcygenmil}) by applying Proposition \ref{flcy}. To apply this proposition, we 
need to consider a variation of local epsilon factors in families of isolated characteristic points. In the case of positive characteristic, this is done in \cite{CLE}. In 

 deduce the existence of epsilon cycles in the case of positive characteristic (Lemma \ref{flatepchar}) 
from Proposition \ref{flcy}, we need to consider the variation of local epsilon factors in families of 
isolated characteristic points, which is done in \cite{CLE}.

Let $k$ be a field and let $\ell$ be a prime number invertible in $k$. In Theorem \ref{epcygenmil}, we define epsilon cycles for $\ell$-adic sheaves on smooth varieties over $k$, under the following assumption: 
\begin{center}
{\em $k$ is the perfection of a finitely generated field over its prime field. }
\end{center}
Before proceeding, we recall two results which are key ingredients for our purpose. The first one is the following result due to Katz--Lang; since this is crucially used in several steps, we need to impose the above assumption on $k$. 
\begin{thm}(\cite[Theorem 1]{KL})\label{KLfin}
Let $k$ be the perfection of a finitely generated field of characteristic $p\geq0$. Let $X$ be a geometrically  connected smooth scheme over $k$. Then the map 
\begin{equation*}
\pi_1(X)^{ab}\to \pi_1(k)^{ab}=G_k^{ab}
\end{equation*}
is surjective. The kernel 
is the product of a finite group with a pro-$p$ group when $p>0$. If $p=0$, it is a finite group. 
\end{thm}

The second one is the ``continuity'' of local epsilon factors proved in \cite{CLE}. 
\begin{thm}(\cite[Theorem 4.9.2]{CLE})\label{contiep}
Let $S$ be a connected scheme of finite type over a perfect field $k$ of characteristic $p>0$. Consider a commutative diagram
\begin{equation*}
\xymatrix{
Z\ar@{^{(}-_>}[r]&U\ar[rr]^-{f}\ar[rd]_-{g}&&Y\ar[ld]\\
&&S&
}
\end{equation*}
of $S$-schemes of finite type. Let ${\cal F}\in D^b_c(U,\Zl)$. We assume that 
\begin{itemize}
\item $Y$ is a smooth $S$-curve. 
\item $Z$ is a closed subscheme of $U$ finite over $S$. 
\item $g$ and $f|_{U\setminus Z}$ are universally locally acyclic relatively to ${\cal F}$. 
\end{itemize}
Let $t\colon Y\to \mathbb{A}^1_S$ be an \'etale $S$-morphism. Then there exists a continuous character $\rho_{{\cal F},t}\colon\pi(S)^{ab}\to \Zl^\times$ with the following property: 
for any perfect field $k'$ with a morphism ${\rm Spec}(k')\to S$, the composition of $G_{k'}^{ab}\to\pi(S)^{ab}\xrightarrow{\rho_{{\cal F},t}} \Zl^\times$ is equal to 
\begin{equation*}
\prod_{z\in Z_{k'}}\delta_{k(z)/k'}^{{\rm dimtot}(R\Phi_{f_{k'}}({\cal F}|_{U_{k'}})_z)}\cdot\varepsilon_{0}(Y_{k'(z)},R\Phi_{f_{k'}}({\cal F}|_{U_{k'}})_z,
dt)\circ {\rm tr}_{k(z)/k'}. 
\end{equation*}
Here the subscripts $(-)_{k'}$ mean the base changes to $k'$. 
\end{thm}
Combining these results, first we verify that local epsilon factors modulo roots of unity do not depend on a choice of a 
uniformizer. 
\begin{lm}\label{changecoor}
Let $k$ be the perfection of a field finitely generated over $\mathbb{F}_p$. Let $X$ be a smooth scheme of finite type over $k$. Let $\mathcal{F}$ 
be an object of $D^b_c(X,\overline{\mathbb{Z}}_\ell)$. Let 
\begin{equation}\label{charB}
\xymatrix{
U\ar[r]^f\ar[d]_j&Y\\X
}
\end{equation}
be a diagram as (\ref{char}). Let $u\in U$ be an at most isolated 
$SS(j^\ast\mathcal{F})$-characteristic point of $f$. For two local parameters $t$ and $t'$ of 
$Y$ around $f(u)$, the ratio 
$\varepsilon_0(Y_{(u)},R\Phi_f(\mathcal{F})_u,dt)\cdot
\varepsilon_0(Y_{(u)},R\Phi_f(\mathcal{F})_u,dt')^{-1}
=(\det R\Phi_f(\mathcal{F})_u)_{[\frac{dt}{dt'}]}$ of the characters of 
$G_{k(u)}^{ab}$ in Theorem \ref{torloc}.4 is of finite order. 
\end{lm}
\proof{
We may assume that $u\to{\rm Spec}(k)$ is an isomorphism and that $U\setminus\{u\}\to Y$ is universally locally acyclic relatively to ${\cal F}$. 
We show that $(\det R\Phi_f(\mathcal{F})_u)_{[\frac{dt}{dt'}]}$ is of finite order. Let $n$ be the Swan conductor of 
$\det R\Phi_{f}(\mathcal{F})_u$. 
If the difference $dt-dt'$ vanishes at $f(u)$, then 
the character is killed by the $n$-th power of $p$ (\cite[Lemma 4.8]{Y3}). 
Take $a\in k^\times$ so that $d(at)-dt'$ vanishes at $f(u)$. 
Let 
\begin{equation}\label{vardiff}
\xymatrix{
\mathbb{G}_{{\rm m},k}\ar[rrd]_{{\rm id}}\ar@{^{(}-_>}[r]&U
\times_{k}\mathbb{G}_{{\rm m},k}\ar[rd]\ar[rr]^{f
\times{\rm id}}&&
Y\times_{k}\mathbb{G}_{{\rm m},k}\ar[ld] \\
&&\mathbb{G}_{{\rm m},k}
}
\end{equation}
be the product of 
\begin{equation*}
\xymatrix{
u\ar[rrd]_\cong\ar@{^{(}-_>}[r]&U\ar[rd]\ar[rr]^{f}&&
Y\ar[ld] \\
&&{\rm Spec}(k)
}
\end{equation*}
 and $\mathbb{G}_{{\rm m},k}$. 
Let $x$ denote the standard coordinate of $\mathbb{G}_{{\rm m},k}$ and 
let $t'':=xt$. By applying Theorem \ref{contiep} to 
(\ref{vardiff}) and ${\rm pr}^\ast{\cal F}$, where ${\rm pr}$ denotes 
the projection $U\times_k\mathbb{G}_{{\rm m},k}\to U$, 
we get a continuous character 
$\rho_{{\rm pr}^\ast{\cal F},t''}\colon\pi_1(\mathbb{G}_{{\rm m},k})^{ab}\to\Zl^\times$ with the property explained in the theorem. 
By Theorem \ref{KLfin}, the composite map $\pi_1(\mathbb{G}_{{\rm m},k})^{ab}\xrightarrow{\rho_{{\rm pr}^\ast{\cal F},t''}}\Zl^\times\to\Zl^\times/\mu$ factors through 
$\pi_1(\mathbb{G}_{{\rm m},k})^{ab}\to G_k^{ab}$. 
Specializing $x\mapsto1,a$, we obtain the assertion. 
\qed
}
\begin{df}\label{overlineep}
Let the notation and assumptions be as in Lemma \ref{changecoor}. 
We denote the composition of $
\varepsilon_0(Y_{(u)},R\Phi_f(\mathcal{F})_u,dt)$ and the quotient map $\Zl^\times\to\Zl^\times/\mu$ by $\bar{\varepsilon}_0(Y_{(u)},R\Phi_f(\mathcal{F})_u)$. 
\end{df}
The character $\bar{\varepsilon}_0(Y_{(u)},R\Phi_f(\mathcal{F})_u)$ 
is independent of a choice of a local parameter $t$ by Lemma \ref{changecoor}. 
It belongs to $\Theta_{k(u)}$. 
\begin{lm}\label{flatepchar}
Let $k$ be the perfection of a finitely generated field of positive characteristic. 
Let $X$ be a smooth scheme of finite type over $k$. Let $\mathcal{F}$ 
be an object of $D^b_c(X,\overline{\mathbb{Z}}_\ell)$. Let the singular support of $\mathcal{F}$ be 
denoted by $C$. 
For a diagram as (\ref{char}) and an at most isolated $C$-characteristic point 
$u\in U$ of $f$, 
put $\varphi(f,u):=\bar{\varepsilon}_0(Y_{(u)},R\Phi_f(\mathcal{F})_u)^{-1}\circ
{\rm tr}_{k(u)/k}$. 
This assignment defines a $\Theta_k$-valued function on isolated 
$C$-characteristic points in the sense of Definition \ref{fcnchar}.1. Moreover, this function is flat in the sense of Definition \ref{fcnchar}.2. 
\end{lm}
\proof{
First we verify that $\varphi(f,u)$ is a $\Theta_k$-valued function on isolated 
$C$-characteristic points. 

When $u$ is not an isolated $C$-characteristic point, $\varphi(f,u)$ vanishes 
since $R\Phi_f(\mathcal{F})_u$ vanishes, hence the property $(a)$ in Definition \ref{fcnchar}.1 is satisfied. 

We show that the property $(b)$ is also satisfied. Consider a diagram of $k$-schemes of the form  (\ref{char2}) and an isolated $C$-characteristic point 
$u'\in U'$ of $f'$. Since the restriction of $\bar{\varepsilon}_0(Y_{(u)},R\Phi_f(\mathcal{F})_u)$ to $G_{k(u')}^{ab}$ equals to 
$\bar{\varepsilon}_0(Y'_{(u')},R\Phi_{f'}(\mathcal{F})_{u'})$, the 
assertion follows from the fact that the composition 
$G_{k(u)}^{ab}\xrightarrow{{\rm tr}_{k(u')/k(u)}}G_{k(u')}^{ab}\to G_{k(u)}^{ab}$ 
is the multiplication by $\deg(u'/u)$. 

Next we show the flatness. Consider a diagram of $k$-schemes of the form (\ref{c}). We need to 
show that the function $\varphi_f\colon|Z|\to\Theta_k$ defined by 
$\varphi_f(z)=\varphi(f_s,z)$, where $s\in S$ denotes  the image of $z$, is flat over $S$. Let $g$ denote the map $Z\to S$ and take a closed point $z\in Z$. 
After replacing $S$ by an \'etale neighborhood of $g(z)$ and shrinking $Z$ around $z$, we may assume that $Z$ is finite over 
$S$. Further replacing $S$ and $Y$ by open coverings, 
we may assume that there
exists $t\in\Gamma(Y,\mathcal{O}_Y)$ which defines an \'etale 
morphism $Y\to\mathbb{A}^1_S$. We further assume that $S$ is connected. Then Theorem \ref{contiep} gives us  a continuous group homomorphism
\begin{equation*}
\rho_{{\rm pr}_1^\ast{\cal F},t}\colon \pi_1(S)^{ab}\to\overline{\mathbb{Z}}_\ell^\times
\end{equation*}
(where ${\rm pr}_1$ denotes the map $U\to X$ in (\ref{c})) with the property that is explained in loc.~cit. Let $k'$ be the normalization of $k$ 
in $S$. By Theorem \ref{KLfin}, $\rho_{{\rm pr}_1^\ast{\cal F},t}$ followed 
by the quotient map $\Zl^\times\to\Zl^\times/\mu$ factors through 
$G_{k'}^{ab}$, 
which we denote by $\xi\colon G_{k'}^{ab}\to\Zl^\times/\mu$. 
Then, for a closed point $s\in S$, we have 
$\prod_{z\in Z_s}\varphi_{f}(z)=\prod_{z\in Z_s}\overline{\varepsilon}_0
(Y_{s,(z)},R\Phi_{f_s}(\mathcal{F}_s)_z)^{-1}\circ {\rm tr}_{k(z)/k}=
\xi|_{G_{k(s)}^{ab}}\circ {\rm tr}_{k(s)/k}=(\xi\circ {\rm tr}_{k'/k})^{\deg(k(s)/k')}$, hence the assertion. 
\qed
}

To prove the existence of epsilon cycles in the case of characteristic $0$, 
we need the following lemma, by which we reduce to the positive characteristic case. 

\begin{lm}\label{relcyc}
Let $S$ be a noetherian regular scheme. Let 
$X$ be a smooth scheme purely of relative dimension $n$ over $S$. 
Let $Z\subset X$ be an integral closed subscheme which is flat purely 
of relative dimension $n-c$ over $S$. 
Let $W$ be a smooth scheme purely of relative dimension $m$ over $S$. 
Let $h\colon W\to X$ be a closed immersion of $S$-schemes. Assume that 
each irreducible component $C_a$ of $Z\times_XW=\cup_aC_a$, 
equipped with the 
reduced subscheme structure, is flat purely of relative dimension $m-c$ over 
$S$. Then, after shrinking $S$ to a dense open subscheme, 
there exists a family of integers $(t_a)_a$ indexed by the set of irreducible components of $Z\times_XW$
 such that, for every morphism 
$s\to S$ from the spectrum of a field, we have 
$h_s^![Z_s]=\sum_at_a[C_{a,s}]$ as cycles supported on 
$(Z\times_XW)_s$, where $(-)_s$ means the base change 
$(-)\times_Ss$.  
\end{lm}
\proof{ 
After shrinking $S$, we may assume that the morphisms $C_a\cap C_b\to S$ are of relative dimension $<m-c$ for distinct indices $a,b$. 
Let $K$ be the complex $\mathcal{O}_Z\otimes^L_{\mathcal{O}_X}\mathcal{O}_W$ of coherent $\mathcal{O}_X$-modules. 
This is supported on $Z\times_XW$. 
Note that $K$ is bounded since $h\colon W\to X$ is a 
regular immersion. 
Let $U\subset W$ be an open neighborhood around the generic points 
of $Z\times_XW$ so that $U\cap C_a$ are disjoint 
and $K|_{U\cap C_a}$ are extensions of finite free 
$\mathcal{O}_{U\cap C_a}$-modules. 
Let $\eta_a$ denote the generic points of $C_a$. 
Let $t_a$ be the lengths of $K_{\eta_a}$ as complexes of 
$\mathcal{O}_{X,\eta_a}$-modules, i.e., the alternating sums of the lengths of 
$H^i(K_{\eta_a})$. 
Let $s\to S$ be a morphism from the spectrum of a field. Then we have 
\begin{equation*}
h^!_s[Z_s]|_{U_s}=[K\otimes^L_{\mathcal{O}_S}k(s)]|_{U_s}
=\sum_at_a[U_s\cap C_{a,s}]. 
\end{equation*}
Thus  $(t_a)_a$ has the desired  property. 
\qed}

\begin{thm}\label{epcygenmil}
Let $k$ be the perfection of a finitely generated field of characteristic $p\geq0$. 
Let $X$ be a smooth scheme of finite type over $k$ 
and let 
$\mathcal{F}$ be an object of $D^b_c(X,\overline{\mathbb{Z}}_\ell)$. Let $SS(\mathcal{F})=\cup_aC_a$ be the irreducible decomposition of the singular support. 
Then there exists a 
unique cycle $\mathcal{E}(\mathcal{F})_k=\sum_a\xi_a\otimes[C_a]$ with 
coefficients in 
$\Theta_k$ (Definition \ref{dfth}) 
satisfying the following property. For a diagram as (\ref{char}) and an at most 
isolated $SS(\mathcal{F})$-characteristic point $u\in U$ of $f$, we have 
\begin{equation*}
\bar{\varepsilon}_0(Y_{(u)},R\Phi_f(\mathcal{F})_u)^{-1}\circ {\rm tr}_{k(u)/k}=
(\mathcal{E}(\mathcal{F})_k,df)_u^{\deg(u/k)}. 
\end{equation*}
\end{thm}
\proof{
When $p>0$, 
it follows from Lemma \ref{flatepchar} and Proposition \ref{flcy}. 

Let $p=0$. 
For each irreducible component $C_a$ of $SS(\mathcal{F})$, 
choose a diagram 
\begin{equation*}
\xymatrix{
U_a\ar[d]_{j_a}\ar[r]^{f_a}&Y_a\\X,
}
\end{equation*}
where $j_a$ is \'etale and $Y_a$ is a smooth $k$-curve, and an isolated $SS(\mathcal{F})$-characteristic 
point $u_a\in U_a$ of $f_a$ at which $df_a$ only meets $C_a$. 
 Let $\xi_a\in\Theta_k$ be an element satisfying the 
equality $\xi_a^{\deg(u_a/k)(C_a,df_a)_{u_a}}=\overline{\varepsilon}_0(Y_{a,(u_a)},
R\Phi_{f_a}(\mathcal{F})_{u_a})^{-1}\circ {\rm tr}_{k(u_a)/k}$ (such $\xi_a$ does exist as the character in the right hand side belongs to $\Theta_k$ and $\Theta_k$ is divisible). We show that the cycle 
$\sum_a\xi_a\otimes[C_a]$ satisfies the desired  condition. 
Let 
\begin{equation*}
\xymatrix{
U\ar[d]_{j}\ar[r]^{f}&Y\\X
}
\end{equation*}
be a diagram with $j$ \'etale and $Y$ being a smooth $k$-curve. Let $u\in U$ be an 
at most isolated $SS(\mathcal{F})$-characteristic point of $f$. 
We need to show the equality 
\begin{equation}\label{equality}
\overline{\varepsilon}_0(Y_{(u)},R\Phi_f(\mathcal{F}))^{-\frac{1}{\deg(u/k)}}\circ {\rm tr}_{u/k}=\prod_a\overline{\varepsilon}_0(Y_{a,(u_a)},R\Phi_{f_a}(\mathcal{F}))^{-
\frac{(C_a,df)_{u}}{\deg(u_a/k)(C_a,df_a)_{u_a}}}\circ {\rm tr}_{u_a/k}.
\end{equation}
Replacing $k$ with a finite extension, we may assume that $u_a$ and $u$ are 
$k$-rational. 

Shrinking $Y_a$ and $Y$, we take \'etale $k$-morphisms 
$Y_a\to\mathbb{A}^1_k$ and $Y\to\mathbb{A}^1_k$. 
Applying Proposition \ref{taihen} to the diagram 
\begin{equation*}
\xymatrix{u_a\ar[rrd]_\cong\ar@{^{(}-_>}[r]&U_a\ar[rd]\ar[rr]^{f_a}&&
Y_a\ar[r]\ar[ld]&\mathbb{A}^1_k\ar[lld] \\
&&{\rm Spec}(k)
}
\end{equation*}
and the counterpart for $u\hookrightarrow U\to Y$, we get commutative diagrams 
of topological groups as in the proposition. 
By Lemma \ref{chebdense} 
and Lemma \ref{relcyc}, 
the equality (\ref{equality}) follows from the case of positive characteristic. 
\qed
}

\begin{df}
Let $k$ be the perfection of a finitely generated field. 
We call the cycle $\mathcal{E}(\mathcal{F})_k$ constructed in the above theorem {\rm the epsilon cycle} of $\mathcal{F}$. 
If no confusions occur,  we omit the subscript $k$ and 
denote it by $\mathcal{E}(\mathcal{F})$. 
\end{df}
\begin{rmk}\label{rmkmult}
In what follows, we will write the group law of $\Theta_k\otimes Z_n(T^\ast X)$, where $Z_n$ denotes the group of $n$-cycles, additively, whereas that of $\Theta_k$ is written multiplicatively. We also use the following abbeviation: for $\chi\in\Theta_k$ and $Z=\sum_bn_b[D_b]\in\mathbb{Q}\otimes Z_n(T^\ast X)$, we will write $\chi^Z$ for $\sum_b\chi^{n_b}\otimes[D_b]\in\Theta_k\otimes Z_n(T^\ast X)$. 
\end{rmk}
\begin{df}\label{rattatetwist}
Let $X$ be a smooth scheme of finite type over  $k$. 
For a constructible complex $\mathcal{F}\in D^b_c(X,\Zl)$ 
and a rational number $r$, we define the $r$-twisted epsilon cycle 
$\mathcal{E}(\mathcal{F})(r)$ to be 
\begin{equation*}
\mathcal{E}(\mathcal{F})(r):=\chi_{\rm cyc}^{rCC(\mathcal{F})}+
\mathcal{E}(\mathcal{F})=\sum_a\chi_{\rm cyc}^{rm_a}\xi_a\otimes[C_a]\in\Theta_k\otimes Z_n(T^\ast X)
\end{equation*}
where $\mathcal{E}(\mathcal{F})=\sum_a\xi_a\otimes[C_a]$ and $m_a$ are coefficients of $CC({\cal F})$, i.e.,  $CC(\mathcal{F})
=\sum_am_a[C_a]$.
\end{df}
\subsection{Properties of epsilon cycles}
Let $f\colon X\to Y$ be a morphism of smooth $k$-schemes. 
Let $C\subset 
T^\ast X$ be a closed conical subset. Assume that 
of $X$ and $C$ are purely of dimension $n$ and that  $Y$ is purely of dimension $m$. 
Further assume that $f$ is proper on the base of $C$ and that  $f_\circ C$ is purely of dimension $m$. Under these conditions, Saito \cite{properp} defines a group homomorphism 
\begin{equation}\label{prps}
f_!\colon Z_n(C)\to Z_m(f_\circ C)
\end{equation}
as follows. Consider the diagram 
\begin{equation*}
\xymatrix{
C\ar[d]&C'\ar[d]\ar[r]\ar[l]&f_\circ C\ar[d]\\
T^\ast X&T^\ast Y\times_YX\ar[l]_-{df}\ar[r]^-{{\rm pr}}&T^\ast Y
}
\end{equation*}
where the left square is cartesian. Then the intersection theory defines the pull-back $(df)^!\colon Z_n(C)={\rm CH}_n(C)\to {\rm CH}_m(C')$ and the pushforward ${\rm pr}_\ast\colon {\rm CH}_m(C')\to {\rm CH}_m(f_\circ C)=Z_m(f_\circ C)$. The map $f_!$ then is defined to be ${\rm pr}_\ast(df)^!$. 
We also use the same symbol $f_!$ for the map $A\otimes Z_n(C)\to A\otimes 
Z_m(f_\circ C)$ tensored with an abelian group $A$. 

\begin{lm}\label{lm}
Let $k$ be the perfection of a finitely generated field. 
Let $X$ be a smooth scheme of finite type over $k$ and $\mathcal{F}\in 
 D^b_c(X,\Zl)$ be a constructible complex of $\Zl$-sheaves on $X$. 
\begin{enumerate}
\item 
Let $\mathcal{G}$ be a smooth $\Zl$-sheaf of finite free $\Zl$-modules on $X$. 
Assume that $X$ is connected. 
Let $k'$ be the normalization of $k$ in the function field of $X$. By Theorem \ref{KLfin}, the composite map $\pi^{ab}_1(X)\xrightarrow{\det(\mathcal{G})}\Zl^\times\to \Zl^\times/\mu$ factors through $G_{k'}^{ab}$; the element of $\Theta_{k'}$ so obtained we denote by the same symbol $\det(\mathcal{G})$. 

Then we have an equality 
\begin{equation*}
\mathcal{E}(\mathcal{G}\otimes^L\mathcal{F})=(\det(\mathcal{G})\circ {\rm tr}_{k'/k})^{\frac{1}{\deg(k'/k)}\cdot
CC(\mathcal{F})}+{\rk\mathcal{G}}\cdot\mathcal{E}(\mathcal{F}). 
\end{equation*}
Here we put 
\begin{equation*}
(\det(\mathcal{G})\circ {\rm tr}_{k'/k})^{\frac{1}{\deg(k'/k)}\cdot
CC(\mathcal{F})}=\sum_a(\det(\mathcal{G})\circ {\rm tr}_{k'/k})^{\frac{m_a}{\deg(k'/k)}}\otimes[C_a]
\end{equation*}
 for $CC({\cal F})=\sum_am_a[C_a]$ (cf. Remark \ref{rmkmult} for our conventions). 

In particular, we have $\mathcal{E}(\mathcal{F}(n))=
\mathcal{E}(\mathcal{F})(n)$.
\item Let $k_1$ be a subfield of $k$ such that $\deg(k/k_1)$ is finite (then $k_1$ is perfect). We also regard $X$ as a smooth scheme over $k_1$. Then we have 
\begin{equation*}
\mathcal{E}(\mathcal{F})_{k}\circ {\rm tr}_{k/k_1}
=
{\deg(k/k_1)}\cdot\mathcal{E}(\mathcal{F})_{k_1}
\end{equation*}
where we put $\mathcal{E}(\mathcal{F})_{k}\circ {\rm tr}_{k/k_1}=\sum_a(\xi_a\circ {\rm tr}_{k/k_1})\otimes[C_a]$ for ${\cal E}({\cal F})=\sum_a\xi_a\otimes[C_a]$. 
\item Let $k'/k$ be an extension of fields which are perfections of finitely generated fields. For an irreducible component $C_a$ of $SS({\cal F})$, let 
$(C_{a,b}')_b$ be the set of irreducible components of $C_a\times_kk'$ and put $[C_a\times_kk']:=\sum_b[C_{a,b}']$. 
Let $\mathcal{E}(\mathcal{F})=\sum_a\xi_a\otimes[C_a]$ be the epsilon cycle. 
Then we have 
\begin{equation*}
\mathcal{E}(\mathcal{F}|_{X_{k'}})=\sum_a\xi_a|_{G_{k'}^{ab}}
\otimes[C_a\times_kk']
\end{equation*}
if one of the following conditions holds: 
\begin{enumerate}
\item $k$ is algebraically closed in $k'$. 
\item $k$ is a finite field. 
\end{enumerate}
\item Let $i\colon X\to X'$ be a closed immersion to a smooth 
$k$-scheme $X'$ of finite type. Then we have 
\begin{equation*}
i_!\mathcal{E}(\mathcal{F})=\mathcal{E}(i_\ast\mathcal{F}). 
\end{equation*}
Note that the conditions used to construct (\ref{prps}) are fulfilled. 
\end{enumerate}
\end{lm}
\proof{
1. 
Take a diagram as (\ref{char}) and an at most isolated 
$SS(\mathcal{F})$-characteristic point $u\in U$ of $f$. 
We have 
\begin{align*}
(\mathcal{E}(\mathcal{G}\otimes^L\mathcal{F}),df)_{u}^{\deg(u/k)}&=
\overline{\varepsilon}_0(Y_{(u)},R\Phi_f(\mathcal{G}\otimes\mathcal{F})_u)
^{-1}\circ {\rm tr}_{u/k}\\&=\overline{\varepsilon}_0(Y_{(u)},\mathcal{G}_{\bar{u}}\otimes R\Phi_f(\mathcal{F})_u)^{-1}\circ {\rm tr}_{u/k}
\\&=(\det\mathcal{G}\circ {\rm tr}_{u/k})^{-{\rm dimtot}R\Phi_f(\mathcal{F})_u}
\overline{\varepsilon}_0(Y_{(u)}, R\Phi_f(\mathcal{F})_u)
^{-{\rm rk}\mathcal{G}}\circ {\rm tr}_{u/k}\\
&=(\det\mathcal{G}\circ {\rm tr}_{k'/k})^{
\deg(u/k')\cdot(CC(\mathcal{F}),df)_{u}}
(\mathcal{E}(\mathcal{F}),df)_u^{\deg(u/k)\cdot{\rm rk}\mathcal{G}}\\
&=((\det\mathcal{G}\circ {\rm tr}_{k'/k})^{\frac{1}{\deg(k'/k)}\cdot 
CC(\mathcal{F})}+
\mathcal{E}(\mathcal{F})^{{\rm rk}(\mathcal{G})},df)^{\deg(u/
k)}_{u}.
\end{align*}

2. 
Consider morphisms of $k_1$-schemes 
$X\xleftarrow{j} U\xrightarrow{f} Y$ where 
$j$ is \'etale and $Y$ is a smooth $k_1$-curve. 
Replacing $Y$ by $Y\times_{k_1}k$ if necessary, 
to calculate local epsilon factors, 
we may assume that the diagram is defined over $k$. 
Then it follows from the characterization of epsilon cycles in Theorem \ref{epcygenmil}. 

3.  First, we consider the case $(a)$. In this case, $C_a\times_kk'$ is already irreducible. Take and fix one  irreducible component $C_a$. By Lemma \ref{exgoodpen} and a limit argument, we can take morphisms of $k$-schemes $X\xleftarrow{j}U\xrightarrow{f}Y$ where $j$ is \'etale and $Y$ is a smooth $k$-curve together with an isolated $SS({\cal F})$-characteristic point $u\in U$ of $f$ at which $df$ only meets $C_a$. Put $k'(u):=k'\otimes_kk(u)$, which is a field extension of $k'$. Then the assertion follows from Theorem \ref{epcygenmil} and the commutativity of the following diagram 
\begin{equation*}
\xymatrix{
G_{k'}^{ab}\ar[r]\ar[d]_-{{\rm tr}_{k'(u)/k'}}&G_k^{ab}\ar[d]^-{{\rm tr}_{k(u)/k}}\\
G_{k'(u)}^{ab}\ar[r]&G_{k(u)}^{ab}. 
}
\end{equation*}

Next, we consider the case $(b)$. 
Let $k''$ be the algebraic closure of $k$ in $k'$; as $k'$ is the perfection of a finitely generated field, this is a finite field. Since the case $k'/k''$ is already treated in $(a)$, we may assume that $k'$ is also a finite field. 

Let ${\cal E}({\cal F}|_{X_{k'}})_{k'}=\sum_{a,b}\xi_{a,b}\otimes[C_{a,b}']$ be the epsilon cycle of ${\cal F}|_{X_{k'}}$ over $k'$. By the assertion $2$., we have $\deg(k'/k)\cdot{\cal E}({\cal F}|_{X_{k'}})_k=\sum_{a,b}(\xi_{a,b}\circ{{\rm tr}_{k'/k}})\otimes[C_{a,b}']$. As the epsilon cycle is \'etale local, we have $\xi_a^{\deg(k'/k)}=\xi_{a,b}\circ{{\rm tr}_{k'/k}}$ for any $b$. As the composition of $G_{k'}^{ab}\to G_k^{ab}\xrightarrow{{\rm tr}}G_{k'}^{ab}$ is given by the multiplication by $\deg{(k'/k)}$ (which is only true when  $k'$ is a finite field), we have $\xi_{a,b}^{\deg(k'/k)}=(\xi_{a,b}\circ{{\rm tr}_{k'/k}})|_{G_{k'}^{ab}}=\xi_a^{\deg(k'/k)}|_{G_{k'}^{ab}}$. The assertion follows.

4. Consider a commutative diagram of $k$-schemes 
\begin{equation*}
\xymatrix{
X\ar[d]_i&U\ar[l]_j\ar[d]_{i'}\ar[rd]^f\\
X'&U'\ar[l]^{j'}\ar[r]^{f'}&Y
}
\end{equation*}
where the left horizontal arrows are \'etale, the square is cartesian, and $Y$ is a smooth 
curve. Since $SS(i_\ast\mathcal{F})=
i_\circ SS(\mathcal{F})$, it suffices to show, for an isolated 
$SS(i_\ast\mathcal{F})$-characteristic point $u'\in U'$ of $f'$, the equality 
\begin{equation*}
\overline{\varepsilon}_0(Y_{(u')},R\Phi_{f'}(i_\ast\mathcal{F})_{u'})=
\overline{\varepsilon}_0(Y_{(u')},R\Phi_{f}(\mathcal{F})_{u'}). 
\end{equation*}
This follows since the canonical map $R\Phi_{f'}(i_\ast\mathcal{F})_{u'}\to
 i'_{\ast} R\Phi_{f}(\mathcal{F})_{u'}$ is an isomorphism. 
\qed
}

\begin{pr}\label{external}
Let $X_1$ and $X_2$ be smooth schemes of finite type over $k$, where $k$ is the perfection of a finitely generated field. 
Take $\mathcal{F}_{i}\in D^b_c(X_i,\Zl)$ for each $i=1,2$. Then 
we have an equality
\begin{equation*}
\mathcal{E}(\mathcal{F}_{1}\boxtimes\mathcal{F}_{2})=
(\mathcal{E}(\mathcal{F}_{1})\boxtimes CC(\mathcal{F}_{2}))
+(CC(\mathcal{F}_{1})\boxtimes
\mathcal{E}(\mathcal{F}_{2})), 
\end{equation*}
where $\mathcal{E}(\mathcal{F}_{1})\boxtimes CC(\mathcal{F}_{2})$ is 
defined as follows. Write $\mathcal{E}(\mathcal{F}_{1})=\sum_a\xi_a\otimes 
[C_a]$ and $CC(\mathcal{F}_{2})=\sum_bn_b\cdot [D_b]$. Then we put 
$\mathcal{E}(\mathcal{F}_{1})\boxtimes CC(\mathcal{F}_{2}):=
\sum_{a,b}\xi_a^{n_b}\otimes[C_a\times D_b]$. The definition of 
$CC(\mathcal{F}_{1})\boxtimes
\mathcal{E}(\mathcal{F}_{2})$ is similar. 
\end{pr}
\proof{

Let $k'/k$ be a finite field extension. 
By Lemma \ref{lm}.2 applied to $k/k_1=k'/k$, it is enough to show the statement after taking the base changes to $k'$. Then we may assume that all irreducible components of $SS({\cal F}_i)$ for $i=1,2$ are geometrically irreducible after replacing $k$ by its 
finite extension if necessary. 
Let $C_{1},C_{2}$ be irreducible components of $SS(\mathcal{F}_{1}),
SS(\mathcal{F}_{2})$ respectively. By Proposition \ref{example}.2, the cycle 
$\mathcal{E}(\mathcal{F}_{1}\boxtimes\mathcal{F}_{2})$ is supported on 
$SS({\cal F}_1\boxtimes{\cal F}_2)=SS(\mathcal{F}_{1})\times SS(\mathcal{F}_{2})$. Hence it suffices to 
compare the coefficients of $[C_1\times C_2]$. 

We may assume that, after a field extension, there exists a 
diagram 
\begin{equation*}
\xymatrix{
U_i\ar[r]^{f_i}\ar[d]&\mathbb{A}^1_{k}\\X_i
}
\end{equation*}
and a $k$-rational isolated 
$SS(\mathcal{F}_{i})$-characteristic point $u_i\in U_i$ 
at which the section $df_i$ 
meets only $C_{i}$ for $i=1,2$. 
Let $f\colon U_1\times U_2\to \mathbb{A}^2_{k}$ be the product of $f_1$ and $f_2$ and let $a\colon\mathbb{A}^2_{k}\to\mathbb{A}
^1_{k}$ be the summation map. 
Let $\xi_i$ be the coefficient of $C_i$ in 
$\mathcal{E}(\mathcal{F}_{i})$ and $\xi$ be that of $C_1\times C_2$ in 
$\mathcal{E}(\mathcal{F}_{1}\boxtimes\mathcal{F}_{2})$. 
Put $u:=(u_1,u_2)\in U_1\times U_2$. Since $u$ is an isolated 
$SS(\mathcal{F}_{1}\boxtimes\mathcal{F}_{2})$-characteristic point 
of $af$, we have 
\begin{equation*}
(\mathcal{E}(\mathcal{F}_{1}\boxtimes\mathcal{F}_{2}),d(af))_{T^\ast (U_1
\times U_2),u}=\overline{\varepsilon}_0(\mathbb{A}^1_{k,(0)},
R\Phi_{af}(\mathcal{F}_{1}\boxtimes\mathcal{F}_{2})_u)^{-1}.
\end{equation*}
Since $d(af)$ only meets $C_1\times C_2$ at $u$, the left hand side equals to 
$\xi^{(C_1,df_1)_{T^\ast U_1,u_1}\cdot(C_2,df_2)_{T^\ast U_2,u_2}}$. 
On the other hand, by Lemma \ref{epconv}, 
the right hand side equals to 
\begin{equation*}
\overline{\varepsilon}_0(\mathbb{A}^1_{k,(0)},R\Phi_{f_1}(\mathcal{F}_1)_{u_1})^{\dt R\Phi_{f_2}
(\mathcal{F}_2)_{u_2}}\cdot\overline{
\varepsilon}_0(\mathbb{A}^1_{k,(0)},
R\Phi_{f_2}(\mathcal{F}_2)_{u_2})^{\dt R\Phi_{f_1}(\mathcal{F}_1)_{u_1}}, 
\end{equation*}
which equals to 
\begin{equation*}
\xi_1^{(C_1,df_1)_{u_1}\cdot(CC(\mathcal{F}_{2}),df_2)_{u_2}}\cdot
\xi_2^{(C_2,df_2)_{u_2}\cdot(CC(\mathcal{F}_{1}),df_1)_{u_1}}, 
\end{equation*}
hence the assertion. 
\qed
}

We recall the definition of the pull-backs of cycles for a properly transversal morphism, which is given in \cite[Definition 7.1.2]{Sai17}. 

Let $X$ and $W$ be smooth schemes over a field $k$ and let $C$ be a closed 
conical subset of $T^\ast X$. 
Assume that every irreducible component of $X$ and $C$ is of 
dimension $n$ and that every irreducible component of $W$ is of 
dimension $m$. 

 Let $h\colon W\to X$ be a properly $C$-transversal $k$-morphism and consider the diagram
\begin{equation*}
\xymatrix{
h^\circ C\ar[d]&W\times_XC\ar[r]\ar[l]\ar[d]&C\ar[d]\\
T^\ast W &W\times_X T^\ast X\ar[r]^-{{\rm pr}}\ar[l]_-{dh}& T^\ast X. 
}
\end{equation*}
By the properly $C$-transversality, $W\times_XC$ is purely of dimension $m$ and the intersection-theoretic pull-back ${\rm pr}^!$ defines a group homomorphism 
\begin{equation*}
{\rm pr}^!\colon Z_n(C)={\rm CH}_n(C)\to{\rm CH}_m(W\times_XC)=Z_m
(W\times_XC). 
\end{equation*}
Note that $dh\colon W\times_XC\to h^\circ C$ is finite by Lemma \ref{ctrfin}. 
\begin{df}\label{!plbk}
Let the notation and assumption be as above. We define the map $h^!\colon Z_n(C)\to Z_m(h^\circ C)$ to be $(-1)^{n-m}$-times the composite map 
\begin{equation*}
Z_n(C)\xrightarrow{{\rm pr}^!}Z_m(W\times_XC)\xrightarrow{dh_\ast} Z_m(h^\circ C). 
\end{equation*}
We also use the same symbol $h^!$ for the map 
$A\otimes Z_n(C)\to A\otimes Z_m(h^\circ C)$ tensored with an abelian group $A$. 
\end{df}

\begin{cor}\label{smpb}
Let $k$ be the perfection of a finitely generated field. 
Let $h\colon W\to X$ be a smooth morphism of smooth schemes of finite type over 
$k$. Assume that each irreducible component of $X$ and $W$ is 
of dimension $n$ and $m$ respectively. Let $\mathcal{F}\in D^b_c(X,\Zl)$. Then we have 
\begin{equation*}
\mathcal{E}(h^\ast\mathcal{F})=h^!(\mathcal{E}(\mathcal{F})(\frac{n-m}{2})), 
\end{equation*}
where the twist $\mathcal{E}(\mathcal{F})(r)$ for a rational number $r$ is defined in Definition \ref{rattatetwist}. 
\end{cor}
\proof{
Since the assertion is \'etale local on $W$, we may assume that $W=X
\times\mathbb{A}^{m-n}_{k}$ and that $h$ is the projection. 
By induction on $m$, we reduce 
the question to the case when 
$W=X\times\mathbb{A}^{1}_{k}$. By Proposition \ref{external}, it is enough to show, for the trivial $\Zl$-sheaf $\mathcal{G}:=\Zl$ on $\mathbb{A}^1_{k}$, 
the equality $\mathcal{E}(\mathcal{G})=\chi_\cyc^{\frac{1}{2}}\otimes[
T^\ast_{\mathbb{A}^1_{k}}\mathbb{A}^1_{k}]$. 

First we prove the assertion when $p\neq2$. 
Consider the Kummer covering $f\colon\mathbb{A}^1_{k}\to 
\mathbb{A}^1_{k}$ defined by 
$t\mapsto t^2$.  The vanishing cycles $R\Phi_f(\mathcal{G})_0$ supported on the origin is concentrated in degree $0$ and its rank is equal to $1$. The corresponding character is a quadratic tamely ramified character. Then, by Definition \ref{epchar00} and Lemma \ref{sameep}, we have 
$\overline{\varepsilon}_0(
\mathbb{A}^1_{k,(0)},R\Phi_f(\mathcal{G})_0)=\chi_{\rm cyc}^{\frac{-1}{2}}$ as characters valued in $\Zl^\times/\mu$. 
 On the other 
hand, the intersection number $(T^\ast_{\mathbb{A}^1_{\mathbb{F}_q}}
\mathbb{A}^1_{\mathbb{F}_q},df)_0$ is $1$, hence the assertion. 

When $p=2$, we argue as follows. Let $S:={\rm Spec}(\mathbb{Z}[\frac{1}{3\ell}])$ and consider the following commutative diagram 
\begin{equation*}
\xymatrix{
S\ar@{^{(}-_>}[r]^{0}&\mathbb{A}^1_S\ar[rd]
\ar[rr]^{\tilde{f}}&&\mathbb{A}^1_S\ar[ld]\ar[r]^{{\rm id}}&\mathbb{A}^1_S
\ar[lld] \\
&&S
}
\end{equation*}
where $\tilde{f}$ is defined by $t\mapsto t^3$. 
This diagram and the trivial $\Zl$-sheaf on $\mathbb{A}^1_S$ satisfy the conditions from $1$ to $5$ given after the diagram (\ref{wq}). 
Then the assertion follows from Lemma \ref{chebdense}, Proposition 
\ref{taihen}, and 
the case when $p\neq2$. 
\qed
}
\begin{cor}\label{sm}
Let $X$ be a connected smooth scheme of finite type over $k$, the perfection of a finitely generated field. 
Put $n:=\dim X$. 
For a smooth $\Zl$-sheaf $\mathcal{F}$ on $X$, we have 
\begin{equation*}
\mathcal{E}(\mathcal{F})=(\det(\mathcal{F})\circ 
{\rm tr}_{k'/k})^{\frac{(-1)^n}{\deg(k'/k)}}\cdot
\chi_\cyc^{\frac{(-1)^{n+1}n\cdot\rk\mathcal{F}}{2}}\otimes[T^\ast_XX]. 
\end{equation*}
Here $k'$ is the normalization of $k$ in $X$. 
\end{cor}
\proof{
This follows from Lemma \ref{lm}.1 and Corollary \ref{smpb}. Recall that we regard the map $\det({\cal F})\colon \pi_1^{ab}(X)\to\Zl^\times/\mu$ as $\pi_1(k')^{ab}\to\Zl^\times/\mu$ due to Katz--Lang (Theorem  \ref{KLfin}). 
\qed
}

\begin{ex}\label{pfexample}
Let $X$ be a smooth connected curve over $k$. Let $\mathcal{F}
\in D^b_c(X,\Zl)$ be 
a constructible complex of $\Zl$-sheaves on $X$. Let $U\subset X$ be an open dense subset 
where $\mathcal{F}$ is smooth. Then we have 
\begin{equation*}
\mathcal{E}(\mathcal{F})=(\det(\mathcal{F}|_U)\circ {\rm tr}_{k'/k})^{\frac{-1}{\deg(k'/k)}}\cdot
\chi_\cyc^{\frac{\rk\mathcal{F}|_U}{2}}\otimes[T^\ast_XX]
+\sum_{x\in X\setminus U}(\overline{\varepsilon}(X_{(x)},\mathcal{F})^{-1}\circ {\rm tr}_{x/k})^{\frac{1}{\deg(x/k)}}
\otimes[T^\ast_xX]. 
\end{equation*}
Here $\overline{\varepsilon}(X_{(x)},\mathcal{F})=
\overline{\varepsilon}_0(X_{(x)},\mathcal{F}_{\eta_x})\cdot
\det(\mathcal{F}_x)^{-1}$. 
\end{ex}
\begin{lm}\label{prodcur}(cf. \cite[Th\'eor\`eme (3.2.1.1)]{Lau}, \cite[Theorem 11.1]{geomep})
Let $X$ be a projective smooth curve over $k$ and $\mathcal{F}\in D^b_c(X,\Zl)$ be a 
constructible complex. Then the product formula 
\begin{equation*}
\det R\Gamma(X_{\bar{k}},\mathcal{F})=(\mathcal{E}(\mathcal{F}),
T^\ast_XX)_{T^\ast X}
\end{equation*}
as an element of $\Theta_k$ holds. 
\end{lm}
\proof{
When $k$ is of positive characteristic, it follows from Theorem \ref{clprod} and Example \ref{pfexample}. 
Let $Z$ be a closed subscheme of $X$ such that $\mathcal{F}$ is 
smooth outside $Z$. 
The case when $k$ is of characteristic $0$ is reduced to the case of positive characteristic by applying Proposition \ref{taihen} to the 
diagram 
\begin{equation*}
\xymatrix{
Z\ar@{^{(}-_>}[r]&X\ar[rd]\ar[rr]^{{\rm id}}&&X\ar[ld] \\
&&{\rm Spec}(k)
}
\end{equation*}
together with Lemma \ref{chebdense}.  
\qed
}

In the following proposition, we show ``$\ell$-independence" of epsilon cycles. 
\begin{pr}\label{indep1}
When $k=\mathbb{F}_q$, a finite field with $q$ elements, 
we identify the group 
$\Theta_{\F_q}$ with $\Zl^\times/\mu\subset\Ql^\times/\mu$ via $\xi\mapsto\xi({\rm Frob}_q)$. 
Let $X$ be a smooth scheme of finite type over $\F_q$. 
Let $F$ be a field of characteristic 
$0$. Let $\mathcal{F}$ and $\mathcal{F}'$ be elements of 
$D^b _c(X,\Zl)$ and $D^b _c(X,\overline{\mathbb{Z}}_{\ell'})$, where $\ell$ and 
$\ell'$ are prime numbers which do not divide $q$.  
Fix embeddings $F\to\Ql,F\to\overline{\mathbb{Q}}_{\ell'}$ of fields. 
Assume that, for all closed points $x$ of $X$, the coefficients of 
the characteristic polynomials $\det(T-{\rm Frob}_x,\mathcal{F}_{\bar{x}})$ 
and $\det(T-{\rm Frob}_x,\mathcal{F}'_{\bar{x}})$ are contained in $F$ and 
give the same elements of $F$. Then the coefficients of the epsilon cycle 
$\mathcal{E}(\mathcal{F})$ (resp. $\mathcal{E}(\mathcal{F}')$) are 
contained in $F^\times\otimes\mathbb{Q}\subset\Ql^\times
\otimes\mathbb{Q}\cong\Ql^\times/\mu$ (resp. $\subset\overline{\mathbb{Q}}_{\ell'}^\times\otimes\mathbb{Q}
\cong \overline{\mathbb{Q}}_{\ell'}^\times/\mu$) and give the same elements of 
$F^\times\otimes\mathbb{Q}$. 
\end{pr}
\proof{

Since the assertion is \'etale local, we may assume that $X$ is affine. Taking an 
immersion $X\to\Proj$ and replacing $\mathcal{F},\mathcal{F}'$ by their 
$0$-extensions, we may assume that $X$ is projective purely of dimension $n$. 
Let $C$ be a closed conical subset of $T^\ast X$ such that 
$\mathcal{F}$ and $\mathcal{F}'$ are micro-supported on $C$ and its 
irreducible components are of dimension $n$. 
By Lemma \ref{exgoodpen}, after replacing $\F_q$ 
by a finite extension, we have a good pencil  
\begin{equation*}
\xymatrix{
X&X_L\ar[l]_\pi\ar[r]^f&\Proj^1. 
}
\end{equation*}
Let $C=\cup_aC_a$ denote the irreducible decomposition. By the properties $5$ and $6$ in Definition \ref{goodpen}, the base $C_a\cap T^\ast_XX$ of every irreducible component $C_a$  is not contained in the exceptional locus of $\pi$. Thus it is enough to show the statement for 
$\pi^\ast\mathcal{F}$ and $\pi^\ast\mathcal{F}'$. Further by the properties $4$, $6$, 
Theorem \ref{epcygenmil}, and \cite[Theorem 2.]{Gabber}, it suffices to show that the epsilon cycles of  the push-forwards 
$Rf_\ast\pi^\ast\mathcal{F}$ and $Rf_\ast\pi^\ast\mathcal{F}'$ have the same coefficients. 

Hence we may assume that $X$ is a projective smooth curve. 
Let $U$ be an open dense subset of $X$ where 
$\mathcal{F}$ and $\mathcal{F}'$ are smooth. Let $x\in X$ be a closed point and $\omega$ be a basis of $\Omega^1_{X_{(x)}}$. We need to show that 
$\varepsilon(X_{(x)},\mathcal{F},\omega)$ and $\varepsilon(X_{(x)},\mathcal{F}',\omega)$ are contained in $F^\times$ and coincide. This follows from 
 \cite[Th\'eor\`eme (3.1.5.4)(iii)]{Lau} and \cite[Th\'eor\`eme 9.8.]{Del}. 
 \qed}

Next we prove a compatibility of the construction of epsilon cycles and 
the pull-back by properly transversal morphism. We follow the method used for characteristic cycles  
in \cite{Sai17}, due to Beilinson. 

We use the theory of the universal hyperplane sections and follow the notation in (\ref{radon}). 
\begin{lm}(\cite[Lemma 3.11]{Sai17})\label{Ltrconical}
We follow the notation in (\ref{radon}). 
Let $\Proj=\Proj^n$ be a projective space and $\Proj^\vee$ be its dual. 
Let $C^\vee\subset T^\ast\Proj^\vee$ be a closed conical subset whose 
irreducible components are of dimension $n$. Let 
$C\subset T^\ast\Proj$ be the closed conical subset given by $C={\bm p}_\circ {\bm p}^{\vee\circ}C^\vee$. 
Then every irreducible component of $C$ is of dimension $n$. 
\end{lm}
\proof{
See loc.~cit.~for the proof. 
\qed
}

\begin{pr}\label{rrrr}
Let $k$ be the perfection of a finitely generated field. 
Let $\mathbb{P}=\mathbb{P}^n_k$ be a projective space, and $\Proj^\vee$ be 
its dual. Let $\mathcal{G}\in D^b_c(\Proj^\vee,\Zl)$ and 
write $\mathcal{F}:=R{\bm p}_\ast{\bm p}^\ast\mathcal{G}$. Let $C^\vee\subset T^\ast
\Proj^\vee$ be a closed conical subset whose irreducible components 
are of dimension $n$. Put $C:={\bm p}_\circ{\bm p}^\circ C^\vee
\subset T^\ast\Proj$. Assume that $\mathcal{G}$ is micro-supported on 
$C^\vee$. 

Let $X$ be a smooth subscheme of $\Proj$ purely of dimension $m$. Suppose that the immersion 
$h\colon X\to\Proj$ is properly $C$-transversal. 
\begin{enumerate}
\item For an element ${\cal E}\in\Theta_k\otimes Z_m(T^\ast X)$, let us write ${\cal E}^0$ for 
${\cal E}-\xi_0\otimes[T^\ast_XX]$ where $\xi_0$ is the coefficient of the $0$-section $T^\ast_XX$ appearing in $\cal E$. The cycles $p_!\mathcal{E}(p^{\vee\ast}\mathcal{G})$ and $p_!p^{\vee !}(\mathcal{E}(\mathcal{G})(\frac{1-\dim X}{2}))$ 
are well-defined and 
we have 
\begin{equation}\label{radtr}
\mathcal{E}(Rp_\ast p^{\vee\ast}\mathcal{G})^0=
(p_!\mathcal{E}(p^{\vee\ast}\mathcal{G}))^0=(p_!p^{\vee !}(
\mathcal{E}(\mathcal{G})(\frac{1-\dim X}{2})))^0.
\end{equation}
In particular, we have 
\begin{equation*}
\mathcal{E}(R{\bm p}_\ast {\bm p}^{\vee\ast}\mathcal{G})^0=
({\bm p}_!\mathcal{E}({\bm p}^{\vee\ast}\mathcal{G}))^0=
({\bm p}_!{\bm p}^{\vee !}(\mathcal{E}(\mathcal{G})(\frac{1-n}{2})))^0.
\end{equation*}
\item We have 
\begin{equation*}
\mathcal{E}(h^\ast\mathcal{F})=h^!(\mathcal{E}(\mathcal{F})(\frac{n-\dim X}{2})).
\end{equation*}
\end{enumerate}
\end{pr}
\proof{
Note that, by \cite[Corollary 3.13.2]{Sai17}, we have $p_\circ p^{\vee\circ}
C^\vee=h^\circ C$. Therefore, by the assumption that $h\colon X\to\Proj$ is properly $C$-transversal, every irreducible component of 
$p_\circ p^{\vee\circ}C^\vee$ has the same dimension as $X$. 
Thus the cycles $p_!\mathcal{E}(p^{\vee\ast}\mathcal{G})$ and $p_!p^{\vee !}(\mathcal{E}(\mathcal{G})(\frac{1-\dim X}{2}))$ 
are well-defined. 

1. First we prove the second equality of (\ref{radtr}). 
By \cite[Corollary 3.13.2]{Sai17}, $p^\vee\colon X\times_\Proj Q\to \Proj^\vee$ is 
$C^\vee$-transversal, which implies that $p^{\vee\ast}\mathcal{G}$ is micro-supported 
on $p^{\vee\circ}C^\vee$. Since $p^\vee\colon X\times_\Proj Q\to\Proj^\vee$ is smooth 
outside $\Delta_X:=\Proj(T^\ast_X\Proj)\subset X\times_\Proj Q$, 
we have $\mathcal{E}(p^{\vee\ast}\mathcal{G})=p^{\vee !}(\mathcal{E}(\mathcal{G})(\frac{1-\dim X}{2}))$ 
outside $\Delta_X$ by Corollary \ref{smpb}. Let $C_1,\dots,C_r$ be the irreducible components of $p^{\vee\circ}C^\vee$ which are contained in $T^\ast(X\times_{\Proj}Q)|_{\Delta_X}$. The cycle 
$\mathcal{E}(p^{\vee\ast}\mathcal{G})-p^{\vee !}(\mathcal{E}(\mathcal{G})(\frac{1-\dim X}{2}))$ is of the form $\sum_{i=1}^r\xi_i\otimes[C_i]$. Then, to show the second equality in (\ref{radtr}), it suffices to check that 
$p_\circ C_i$ are contained in the $0$-section. 
By the assumption that $h\colon X\to\Proj$ is $C$-transversal, 
the pair $(p,p^\vee)$ is $C^{\vee}$-transversal around $\Delta_X\subset
X\times_\Proj Q$ by \cite[Corollary 3.13.1]{Sai17}, hence the assertion. 

We prove the first equality in (\ref{radtr}). Since both of $\mathcal{E}(Rp_\ast p^{\vee\ast}\mathcal{G})$ and $p_!\mathcal{E}(p^{\vee\ast}\mathcal{G})$ are 
supported on $h^\circ C=p_\circ p^{\vee\circ}C^\vee$,  it suffices to show the equality 
\begin{equation}\label{push}
(\mathcal{E}(Rp_\ast p^{\vee\ast}\mathcal{G}),df)^{\deg(u/
k)}_u=
(p_!\mathcal{E}(p^{\vee\ast}\mathcal{G}),df)^{\deg(u/k)}_u 
\end{equation}
for every diagram 
\begin{equation*}
X\xleftarrow{j}U\xrightarrow{f}\mathbb{A}^1_{k}, 
\end{equation*}
where $j$ is \'etale and $f$ is smooth, and every at most 
isolated $h^\circ C$-characteristic point $u\in U$ of $f$. 
By Theorem \ref{epcygenmil}, the left hand side of (\ref{push}) equals to 
$\overline{\varepsilon}_0(\mathbb{A}^1_{k(u)},R\Phi_f(h^\ast\mathcal{F})_u)^{-1}
\circ {\rm tr}_{u/k}$. By \cite[Corollary 3.15]{Sai17}, the composition $fp\colon U\times_\Proj Q\to 
\mathbb{A}_{k}^1$ has finitely many $p^{\vee\circ} C^\vee$-characteristic points. Hence the right hand side of (\ref{push}) 
equals to $\prod_v(\mathcal{E}(p^{\vee\ast}\mathcal{G}),d(fp))_v
^{\deg(v/k)}$ where $v$ runs through $p^{\vee\circ} C^\vee$-characteristic points of $fp$ over $u$. Furthermore, by Theorem \ref{epcygenmil}, this equals to 
$\prod_v\overline{\varepsilon}_0(\mathbb{A}^1_{k(v)},R\Phi_{fp}(p^{\vee\ast}\mathcal{G})_v)^{-1}\circ {\rm tr}_{v/k}$. Thus the equality (\ref{push}) follows from 
the isomorphism 
\begin{equation*}
R\Phi_f(Rp_\ast p^{\vee\ast}\mathcal{G})_u\xrightarrow{\cong}
\bigoplus_v{\rm Ind}_{G_v}^{G_u}R\Phi_{fp}(p^{\vee\ast}\mathcal{G})_v. 
\end{equation*}

2. By the proper base change theorem, we have an isomorphism 
$h^\ast\mathcal{F}\to Rp_\ast p^{\vee\ast}\mathcal{G}$. 

Therefore, by 1, we have 
\begin{align*}
\mathcal{E}(h^\ast\mathcal{F})^0=
(p_!p^{\vee!}(&\mathcal{E}(\mathcal{G})(\frac{1-\dim X}{2})))^0\\&=(h^!{\bm p}_!{\bm p}^{\vee!}(\mathcal{E}(\mathcal{G})(\frac{1-\dim X}{2})))^0=(h^!(
\mathcal{E}(\mathcal{F})(\frac{n-\dim X}{2})))^0. 
\end{align*}
By the assumption that the immersion $h$ is 
properly $C$-transversal, $X$ intersects the smooth locus of $\mathcal{F}$. 
Hence the coefficients of the $0$-section in both of $\mathcal{E}(
h^\ast\mathcal{F})$ and $h^!(\mathcal{E}(\mathcal{F})(\frac{n-\dim X}{2}))$ coincide. Thus the assertion 
follows. 
\qed
}

Before stating Corollary \ref{CCR}, we recall the  definitions of the Radon transform of $\ell$-adic sheaf  
and the Legendre transform of closed conical subset. 

Let $\mathcal{F}$ be an element of $D^b_c(\Proj,\Zl)$. 
We define the Radon transform $R\mathcal{F}$ of $\mathcal{F}$ by setting 
$R\mathcal{F}:=R{\bm p}^\vee_\ast{\bm p}^\ast\mathcal{F}[n-1]\in D^b_c(\Proj^\vee,\Zl)$. 

Let $C$ be a closed conical subset of $T^\ast\Proj$ whose irreducible components $C_a$ are of dimension $n={\rm dim}\Proj$. Let $A:=\sum_a\beta_a\otimes
[C_a]$ be a cycle supported on $C$ with coefficients in $\Theta_k$. We define the Legendre transform 
$LA$ by setting $LA:=(-1)^{n-1}\cdot{\bm p}^\vee_!{\bm p}^!A$. 
Since the definition of ${\bm p}^!A$ involves the sign $(-1)^{n-1}$, 
that of $LA$ does not involve the sign. 
\begin{cor}\label{CCR}
Let $\mathcal{F}$ be an element of $D^b_c(\Proj,\Zl)$. 
We use the abbreviation ${\cal E}^0$ explained in Proposition \ref{rrrr}.1. We have 
\begin{equation*}
\mathcal{E}(R\mathcal{F})^0=(L(\mathcal{E}(\mathcal{F})(\frac{1-n}{2})))^0.
\end{equation*}
\end{cor}
We will show the equality $\mathcal{E}(R\mathcal{F})=
L(\mathcal{E}(\mathcal{F})(\frac{1-n}{2}))$ in Corollary \ref{corLR}. 
\proof{
This is a restatement of Proposition \ref{rrrr}.1, after exchanging $\Proj$ and $\Proj^\vee$. 
\qed
}

\begin{thm}\label{prtr}
Let $k$ be the perfection of a finitely generated field. Let $X$ be a smooth scheme of finite type over $k$. Let $\mathcal{F}\in 
 D^b_c(X,\Zl)$. Let $h\colon W\to X$ be a properly 
$SS(\mathcal{F})$-transversal $k$-morphism from a smooth $k$-scheme $W$ of finite type. Assume that 
every irreducible component of $X$ and $W$ are of dimension $n$ and $m$ 
respectively. Then we have 
\begin{equation*}
\mathcal{E}(h^\ast\mathcal{F})=h^!(\mathcal{E}(\mathcal{F})(\frac{n-m}{2})).
\end{equation*}
\end{thm}
\proof{
Decomposing $h$ into 
$W\to W\times X \to X$, we assume that 
$h$ is either of a smooth morphism or an immersion. 
The smooth case follows from Corollary \ref{smpb}. 

Suppose that $h$ is an immersion. 
First consider the case when $X$ is a projective space $\Proj$. 
If $\mathcal{F}=R{\bm p}_\ast{\bm p}^{\vee\ast}
\mathcal{G}$ for some ${\cal G}\in D^b_c(\Proj^\vee,\Zl)$, then it is proved in 
Proposition \ref{rrrr}.2. Let $\mathcal{F}\in D^b_c(\Proj,\Zl)$ be any object. Let $R^\vee$ denote the dual Radon transform $R{\bm p}_\ast{\bm p}^{\vee\ast}[n-1]$. 
Since $\mathcal{F}$ is isomorphic to $R^\vee R{\cal F}$ up to a smooth sheaf and the assertion for smooth 
sheaves follows from Corollary \ref{sm}, the theorem  follows in the case when $h$ is an immersion 
to $\mathbb{P}$. 

We show the general case. Since the assertion is local on $W$, we may assume that 
$X$ is affine and take an immersion $i\colon X\to\Proj$. Furthermore, after shrinking $W$ if necessary, we may assume that 
there is a smooth subscheme $V\subset \Proj$ such that $X\cap V=W$ and 
the intersection is transversal. Then the immersion $\tilde{h}\colon V\to\Proj$ is 
properly $ SS(i_!\mathcal{F})$-transversal around $W\subset V$. Then the assertion follows from the equality ${\cal E}(\tilde{h}^\ast i_!{\cal F})=\tilde{h}^!({\cal E}(i_!{\cal F})(\frac{n-m}{2}))$, which holds true around $W\subset V$. 
\qed
}

\subsection{Epsilon cycles for tamely ramified sheaves}

Let $k$ be the perfection of a finitely generated field of characteristic $p\neq\ell$. 
In this subsection, we calculate the epsilon cycles of  tamely ramified 
$\Zl$-sheaves.  

Let $X$ be a smooth scheme of finite type 
over $k$ and let $D\subset X$ be a simple normal crossings divisor. 
Let $U$ be the complement of $D$ in $X$. 
Let $(D_a)_{a\in A}$ denote the irreducible components of $D$. 
For a subset $B\subset A$, we denote by $D_B$ the 
intersection $\cap_{a\in B}D_a$. 

For simplicity, we assume that $X$ is 
connected and of dimension $n$. 
Then $D_B$ is a smooth closed subscheme of $X$ 
 purely of dimension $n-\lvert B\rvert$. 

Let $\mathcal{F}$ be a non-zero smooth $\Zl$-sheaf of 
free $\Zl$-modules on $U$  tamely ramified along $D$. 
Let $j\colon U\to X$ be the inclusion. We have 
\begin{align}
&SS(j_!\mathcal{F})=\cup_BT^\ast_{D_B}X\\
&CC(j_!\mathcal{F})=\sum_B(-1)^n{\rm rk}\mathcal{F}[T^\ast_{D_B}X], 
\end{align}
where $B$ runs through the subsets of $A$  (see \cite[4.2, 7.3]{Sai17}). 

For each $a\in A$, let $\xi_a$ be the generic point of $D_a$ and 
denote by $k_a$ the normalization 
of $k$ in the residue field at $\xi_a$. Since $\mathcal{F}$ is tamely ramified,  its restriction 
to the henselization $X_{(\xi_a)}$ gives a representation $V_{a}$ of 
the tame inertia group $I_a$ of the trait $X_{(\xi_a)}$. 
Note that $I_a$ is isomorphic to $\varprojlim_{n\neq p}\mu_n(\bar{k})$, 
where $n$ runs through the integers $\geq 1$ prime to $p$ 
and $\mu_n(\bar{k})$ denotes the group of $n$-th roots of unity in an algebraic 
closure $\bar{k}$ of $k_a$, and that 
we have $\sigma^\ast V_a\cong V_a$ for each 
$\sigma\in{\rm Gal}(\bar{k}/k_a)$.  
Thus we get a character $J(V_a)\colon{\rm Gal}(\bar{k}/k_a)\to\Zl^\times/\mu$ 
as constructed in Definition \ref{js}.2. 
We define 
\begin{equation*}
J_a:=(J(V_a)\circ {\rm tr}_{k_a/k})^{\frac{1}{\deg(k_a/k)}}. 
\end{equation*}
This is an element of $\Theta_k$. 

\begin{pr}
Let the notation assumptions be as above. Assume that $X$ is connected and that $\mathcal{F}$ is 
tamely ramified along $D$. 
For a subset $B$ of $A$, define 
\begin{equation*}
\chi_B:=(\det(\mathcal{F})\circ {\rm tr}_{k'/k})^{\frac{(-1)^n}{\deg(k'/k)}}\cdot\chi_\cyc^{\frac{|B|-n}{2}(-1)^{n}\rk\mathcal{F}}\cdot\prod_{a\in B}
J_a^{(-1)^n}, 
\end{equation*}
where $k'$ is the normalization of $k$ in the function field of $X$. 
Then, we have 
\begin{equation*}
\mathcal{E}(j_!\mathcal{F})=\sum_B\chi_B\otimes[T^\ast_{D_B}X]. 
\end{equation*}
\end{pr}
\proof{
Let $B\subset A$ be a subset and let $m:=|B|$ be the cardinality. 
Let $x\in D_B$ be a closed point which is not contained in 
$D_a$ for any $a\in A\setminus B$. 
For $1\leq i\leq n$, let $E_i\subset\mathbb{A}^n_{k(x)}$ be the 
$i$-th coordinate hyperplane and put $E:=\cup_{1\leq i\leq m}E_i$. 
After replacing $X$ by an \'etale 
neighborhood of $x$, we find an \'etale morphism $f\colon 
X\to\mathbb{A}^n_{k(x)}$ 
such that $x$ maps to the origin and the pull-backs of the divisors 
$(E_i)_{1\leq i\leq m}$ coincide with $(D_a)_{a\in B}$ if we put a suitable  
numbering on $B$. 

Let $\pi_1^{\rm tame}(X_{(x)}\setminus D)$ be the fundamental 
group which classifies finite \'etale coverings of $X_{(x)}\setminus D$ 
tamely ramified along $D$. Let $\pi_1^{\rm tame}
(\mathbb{A}^n_{k(x)}\setminus E)$ be the one which 
classifies finite \'etale coverings of $\mathbb{A}^n_{k(x)}\setminus E$ tamely ramified along both of  $E$ and 
$\Proj^n_{k(x)}\setminus\mathbb{A}^n_{k(x)}$. 
Then, the morphism $\pi_1^{\rm tame}(X_{(x)}\setminus D)\to
\pi_1^{\rm tame}(\mathbb{A}^n_{k(x)}
\setminus E)$ induced from $f$ is an isomorphism. 
Thus we may assume that $X=\mathbb{A}^n_{k'}$ for some finite extension 
$k'$ of $k$ and $D=E$, and that the sheaf 
$\mathcal{F}$ is also tamely ramified along 
$\Proj^n_{k'}\setminus\mathbb{A}^n_{k'}$. 

Fix a geometric point $\overline{\eta}$ over the generic point of $\mathbb{A}^n_{k'}$. 
For $1\leq i\leq m$, let $I_i$ be the tame inertia group of the henselization 
of $\mathbb{A}^n_{k'}$ at the generic point of $E_i$. 
Note that the canonical map $I_i\to\pi_1^{\rm tame}(\mathbb{A}^n_{k'}\setminus E)$ is injective and its image is a normal subgroup.  
Let $V_i$ be the representation of $I_i$ associated with $\mathcal{F}$. After replacing $\mathcal{F}$ by its subquotients, we may assume that ${\cal F}\otimes_{\Zl}\Ql$ is an irreducible smooth $\Ql$-sheaf. Then the representation $V_i\otimes\Ql$ of $I_i$ is semi-simple, i.e., is decomposed into the direct sum of $1$ dimensional representations as $I_i$ can be identified with a normal subgroup of 
$\pi_1^{\rm tame}(\mathbb{A}^n_{k'}\setminus E)$.  
Let $\chi_i\colon I_i\to\Ql^\times$ be a character whose corresponding representation appears in $V_i\otimes\Ql$. Then $V_i\otimes\Ql$ is the direct sum of finitely many copies of $\chi_i$ and its conjugates. 
Since $k'$ contains only finitely many roots of unity, $\chi_i$ factors through a quotient $I_i\to\mu_{d_i}(\bar{k})$ 
for some integer $d_i\geq1$ prime to $p$. 
Further extending $k'$ to a finite extension, we may assume that 
$\mu_{d_i}(\bar{k})$ are contained in $k'$ for all $i$. Then we have a canonical group homomorphism $\pi_1^{\rm tame}(\mathbb{A}^n_{k'}\setminus E)\to \mu_{d_i}(\bar{k})\times{\rm Gal}(\bar{k}/k)$ such that the composition 
$I_i\to\pi_1^{\rm tame}(\mathbb{A}^n_{k'}\setminus E)\to \mu_{d_i}(\bar{k})\times{\rm Gal}(\bar{k'}/k')\to\mu_{d_i}(\bar{k})$, where the last map is the first projection, is the canonical surjection. Therefore, for each $i$, there 
exists a 
smooth $\Zl$-sheaf $\mathcal{G}_{i}$ of rank $1$ of finite order on 
$\mathbb{A}^1_{k'}\setminus0$ tamely ramified at $0,\infty$ whose pull-back to $\mathbb{A}^n_{k'}\setminus E$ via the $i$-th projection ${\rm pr}_i\colon\mathbb{A}^n_{k'}\to\mathbb{A}^1_{k'}$ has the same local monodromy at the generic point of $E_i$ as $\chi_i$. Then ${\cal F}\otimes\bigotimes_i{\rm pr}_i^\ast{\cal G}_i^{-1}$ is unramified along $E$ and it extends to a smooth $\Zl$-sheaf $\mathcal{H}$ 
on $\mathbb{A}^n_{k'}$. 
By Lemma \ref{lm}.1 and Proposition \ref{external}, 
the coefficient of $[T^\ast_{\cap_{1\leq i\leq m}E_i}\mathbb{A}^n_{k'}]$ 
in $\mathcal{E}(j_!\mathcal{F})$ equals to 
\begin{equation*}
(\det(\mathcal{H})\circ{\rm tr}_{k'/k})^{
\frac{1}{\deg(k'/k)}\cdot(-1)^n}\cdot\chi_\cyc^{\frac{m-n}{2}(-1)^n\rk\mathcal{F}}\cdot 
\prod_{1\leq i\leq m}
(\overline{\varepsilon}_0(\mathbb{A}^1_{k',(0)},\mathcal{G}_i)\circ
{\rm tr}_{k'/k})^{\frac{(-1)^n}{\deg(k'/k)}\rk\mathcal{F}}. 
\end{equation*}

Since we have $(\overline{\varepsilon}_0
(\mathbb{A}^1_{k',(0)},\mathcal{G}_i))^{\rk\mathcal{F}}
=J(V_i)$, the assertion follows. 
\qed
}

\section{Radon Transform and Product Formula}\label{rtpf}
Let $k$ be a field. In $5.2$ and $5.3$, we will assume that $k$ is the perfection of a finitely generated field and take a prime number $\ell$ that is invertible in $k$. 
\subsection{Reminder on the Chow groups of projective space bundles}
In this preliminary subsection, we quickly recall necessary results on the Chow groups of projective space bundles, in the same manner as \cite[6.1]{Sai17}. 

Let $X$ be a scheme of finite type over $k$. We write ${\rm CH}_\bullet(X)=\bigoplus_i{\rm CH}_i(X)$ for the Chow group of $X$. Let $\mathbb{Z}[h]$ be the polynomial ring over $\mathbb{Z}$ with one variable $h$ and put 
${\rm CH}_\bullet(X)[h]:={\rm CH}_\bullet(X)\otimes_{\mathbb{Z}}\mathbb{Z}[h]$. We consider ${\rm CH}_\bullet(X)[h]$ as a module over the ring ${\rm End}({\rm CH}_\bullet(X))[h]:={\rm End}({\rm CH}_\bullet(X))\otimes_{\mathbb{Z}}\mathbb{Z}[h]$. 

For a vector bundle $E$ of rank $n+1$ over $X$, put 
$c_h(E):=\sum_{q=0}^{n+1}c_q(E)h^{n+1-q}\in{\rm End}({\rm CH}_\bullet(X))[h]$. 
\begin{lm}\label{Chowprojlm}
Let $X$ be a scheme of finite type over $k$. Let $E$  be a vector bundle of rank $n+1$ over $X$. 
\begin{enumerate}
\item Let $\pi\colon\Proj(E)\to X$ denote the projection.  The map 
\begin{equation*}
\alpha_E\colon{\rm CH}_\bullet(X)[h]\to{\rm CH}_\bullet(\Proj(E))
\end{equation*}
given by $ah^q\mapsto c_1({\cal O}(1))^q\cap\pi^\ast a$ is surjective and the kernel is equal to $c_h(E)\cdot 
{\rm CH}_\bullet(X)[h]$. 
\item(\cite[Lemma 6.2]{Sai17}) Let $i\colon F\to E$ be an injection of vector bundles over $X$. 
\begin{enumerate}
\item The diagram 
\begin{equation*}
\xymatrix{
{\rm CH}_\bullet(\Proj(E))\ar[rr]^-{i^\ast}&&{\rm CH}_\bullet(\Proj(F))\\
&{\rm CH}_\bullet(X)[h]\ar[ru]_-{\alpha_F}\ar[lu]^-{\alpha_E}&
}
\end{equation*}
is commutative. 
\item Let $K$ be the cokernel of $i$. Then the diagram 
\begin{equation*}
\xymatrix{
{\rm CH}_\bullet(\Proj(F))\ar[r]^-{i_\ast}&{\rm CH}_\bullet(\Proj(E))\\
{\rm CH}_\bullet(X)[h]\ar[u]_-{\alpha_F}\ar[r]_-{c_h(K)\cdot}&{\rm CH}_\bullet(X)[h]\ar[u]_-{\alpha_E}
}
\end{equation*}
is commutative. 
\end{enumerate}
\end{enumerate}
\end{lm}
\proof{
1. By \cite[Theorem 3.3.(b)]{Ful}, the assertion follows once we know that $\alpha_E$ kills $c_h(E)\cdot ah^j$ for $a\in
{\rm CH}_\bullet(X)$ and $j\geq0$. By the same theorem together with Proposition 3.1.(a) in loc.~cit., it suffices to show that $\pi_\ast(c_1({\cal O}(1))^{j'}\cap\alpha_E(c_h(E)\cdot ah^j))=\pi_\ast(\sum_{q=0}^{n+1}c_1({\cal O}(1))^{n+1-q+j+j'}\cap\pi^\ast(c_q(E)\cap a))$ is zero for any $j'\geq0$. The latter one is equal to $\sum_{q=0}^{n+1}s_{1-q+j+j'}(E)\cap c_q(E)\cap a$, which is zero as $s(E)c(E)=1$. 
\qed}

\subsection{Epsilon class and product formula}
In this subsection, we introduce the notion of epsilon class, which is an analogue of 
characteristic class \cite[Section 6]{Sai17}. 
Various results for characteristic class can be proved for epsilon class essentially in the same way.  Especially, we describe the epsilon classes of the Radon transforms in Proposition \ref{prRR}, following Beilinson's method in \cite[Section 7]{Sai17}. Using this, we state and prove our product formula in Theorem \ref{thmEP}. 

Let $X$ be a scheme of finite type over $k$. We say that $X$ is {\it embeddable} if there exists a closed immersion $i\colon X\to M$ into a smooth $k$-scheme $M$. 

Let $X$ be a scheme of finite type over $k$ that is embeddable. Let $i\colon X\to M$ be a closed immersion into a smooth $k$-scheme purely of dimension $n$. By Lemma \ref{Chowprojlm}.1, 
we identify ${\rm CH}_\bullet(X)=\oplus_{i=0}^n{\rm CH}_i(X)$ with 
${\rm CH}_n(\Proj((X\times_MT^\ast M)\oplus\mathbb{A}^1_X))$ by the map  
\begin{equation}\label{chowpr}
{\rm CH}_\bullet(X)\to{\rm CH}_n(\Proj((X\times_MT^\ast M)\oplus\mathbb{A}^1_X))
\end{equation}
sending $(a_i)_i$ to $\sum_ic_1(\mathcal{O}(1))^i\cap \pi^\ast a_i$, 
where $\pi\colon\Proj((X\times_MT^\ast M)\oplus\mathbb{A}^1_X)\to X$ denotes the projection. 
Tensoring with $\Theta_k$, we also identify 
$ \Theta_k\otimes{\rm CH}_\bullet(X)$ with 
$\Theta_k\otimes{\rm CH}_n(\Proj((X\times_MT^\ast M)\oplus\mathbb{A}^1_X))$. 

Let $\mathcal{F}$ be an element  of $D^b_c(X,\Zl)$. Let
${\cal E}(i_\ast{\cal F})=\sum_a\xi_a\otimes[C_a]$ be the epsilon cycle of the $0$-extension $i_\ast{\cal F}$ to $M$, where $C_a$ runs through the irreducible components of $SS(i_\ast{\cal F})$. Note that $C_a$ are contained in $X\times_MT^\ast M$ as $i_\ast{\cal F}$ is supported on $X$. 
\begin{df}\label{epclas}
Let the notations be as above. For an element ${\cal F}\in D^b_c(X,\Zl)$, 
we 
define {\rm the 
epsilon class} $\varepsilon_X(\mathcal{F})$ of $\mathcal{F}$ by 
setting 
\begin{equation*}
\varepsilon_X(\mathcal{F})=\sum_a\xi_a\otimes[\Proj(C_a\times_k\mathbb{A}^1_k)]\in
\Theta_k\otimes{\rm CH}_n(\Proj((X\times_MT^\ast M)\oplus\mathbb{A}^1_X))
=\Theta_k\otimes{\rm CH}_\bullet(X). 
\end{equation*}
\end{df}
First we check that this gives a well-defined object. 
\begin{lm}
The element $\varepsilon_X({\cal F})\in\Theta_k\otimes{\rm CH}_\bullet(X)$ constructed in Definition \ref{epclas} is independent of choices of $M$ or $i$.
\end{lm}
\proof{
The proof goes similarly as \cite[Lemma 6.6]{Sai17}. Let $j\colon X\to N$ be another closed immersion into a smooth $k$-scheme $N$ purely of dimension $m$. By considering the product $M\times_kN$ and the projections, we may assume that there exists a smooth $k$-morphism $f\colon M\to N$ that is compatible with $i$ and $j$. Let $df\colon X\times_NT^\ast N\to X\times_MT^\ast M$ be the induced injection of vector bundles. By Lemma \ref{Chowprojlm}.2.(a), the diagram
\begin{equation*}
\xymatrix{
&\Theta_k\otimes{\rm CH}_\bullet(X)\ar[rd]^-{\cong}\ar[ld]_-{\cong}&\\
\Theta_k\otimes{\rm CH}_n(\Proj((X\times_MT^\ast M)\oplus\mathbb{A}^1_X))\ar[rr]^-{{\rm id}_{\Theta_k}\otimes(df)^\ast}&&
\Theta_k\otimes{\rm CH}_m(\Proj((X\times_NT^\ast N)\oplus\mathbb{A}^1_X))}
\end{equation*}
is commutative. Therefore, it suffices to check that each irreducible component of $SS(i_\ast{\cal F})$ meets properly the image of $df$ and that the pull-back of ${\cal E}(i_\ast{\cal F})$ is equal to ${\cal E}(j_\ast{\cal F})$. Since the claim is \'etale local on $X$, we may assume that there exists a section $s\colon N\to M$ of $f$ with $sj=i$. Then the assertion follows from $s_\circ SS(j_\ast{\cal F})=SS(i_\ast{\cal F})$ and $s_! {\cal E}(j_\ast{\cal F})={\cal E}(i_\ast{\cal F})$ (Lemma \ref{lm}.4). 
\qed
}

Let $K(X,\Zl)$ be the Grothendieck group of the triangulated category 
$D^b_c(X,\Zl)$. The epsilon class defines a group homomorphism 
\begin{equation*}
\varepsilon_X\colon K(X,\Zl)\to\Theta_k\otimes{\rm CH}_\bullet(X). 
\end{equation*}

\begin{lm}(cf. \cite[Lemma 6.9]{Sai17})\label{degepcls}
Let $X$ and $\mathcal{F}$ be as in Definition \ref{epclas}. Suppose that $X$ is smooth. \begin{enumerate}
\item The dimension $0$-part $\varepsilon_{X,0}(\mathcal{F})
\in \Theta_k\otimes
{\rm CH}_0(X)$ is equal to the intersection product $(\mathcal{E}(\mathcal{F}),
T^\ast_XX)_{T^\ast X}$ with the $0$-section. 
\item Assume that $X$ is connected and let $n$ be its  dimension. 
Let ${\rm rk}^\circ({\cal F})$ and $\det^\circ(\mathcal{F})$ be the 
rank and the determinant character 
of the restriction of $\mathcal{F}
$ to a dense open subset $U$ where 
$\mathcal{F}$ is smooth. Let $k'$ be the normalization of $k$ in $X$. Theorem \ref{KLfin} implies that the composite map $\pi_1(U)^{ab}\xrightarrow{\det^\circ(\mathcal{F})}\Zl^\times\to\Zl^\times/\mu$ factors through $\pi_1(k')^{ab}$. We let $\det^\circ(\mathcal{F})$ denote the induced map $\pi_1(k')^{ab}\to\Zl^\times/\mu$.  

The dimension $n$-part 
$\varepsilon_{X,n}(\mathcal{F})\in\Theta_k\otimes
{\rm CH}_n(X)=\Theta_k$ is equal to 
$(\det^\circ(\mathcal{F})\circ {\rm tr}_{k'/k})^{\frac{(-1)^n}{\deg(k'/k)}}\cdot \chi_\cyc^{\frac{(-1)^{n+1}n}{2}{\rm rk}^\circ(\mathcal{F})}$. 
\end{enumerate}
\end{lm}
\proof{
1. We may assume that $X$ is purely of dimension $n$. Let $j\colon T^\ast X\to\mathbb{P}(T^\ast X\oplus \mathbb{A}^1_X)$ be the open immersion given by $v\mapsto(v,1)$ and $0_X\colon X\to T^\ast X$ be the $0$-section. Then the projection ${\rm CH}_n(\mathbb{P}(T^\ast X\oplus \mathbb{A}^1_X))\cong {\rm CH}_\bullet (X)\to{\rm CH}_0(X)$ to the dimension $0$-part is equal to the composition $0_X^\ast j^\ast$ since ${\cal O}(1)$ is trivial on $T^\ast X$. The assertion follows. 

2. Let $\pi\colon \mathbb{P}(T^\ast X\oplus \mathbb{A}^1_X)\to X$ denote the projection. 
Since $\pi_\ast(c_1({\cal O}(1))^i\cap\pi^\ast a)$ is zero when $i<n$ and is equal to $a$ when $i=n$, the projection ${\rm CH}_n(\mathbb{P}(T^\ast X\oplus \mathbb{A}^1_X))\cong {\rm CH}_\bullet (X)\to{\rm CH}_n(X)$ is equal to $\pi_\ast$. 

Take a closed point $x\in U$ and form the cartesian square
\begin{equation*}
\xymatrix{
\mathbb{P}^n_x\ar[r]^-\pi\ar[d]_-i&x\ar[d]_-i\\
\mathbb{P}(T^\ast X\oplus \mathbb{A}^1_X)\ar[r]^-\pi&X. 
}
\end{equation*}
Then the assertion follows from $i^\ast\pi_\ast=\pi_\ast i^\ast$ and Corollary \ref{sm}. 
\qed
}

In the following lemma, we compute the epsilon class of 
the pull-back by a properly transversal immersion. 
\begin{lm}\label{trprepcl}
Let $X$ be a smooth scheme of finite type purely of dimension $n$ over $k$ and $\mathcal{F}$ be an element of 
$D^b_c(X,\Zl)$. Let $W$ be a smooth $k$-scheme purely of dimension $m=n-c$ and $h\colon W\to X$ be a properly $SS(\mathcal{F})$-transversal 
closed immersion. Then we have 
\begin{equation*}
\chi_\cyc^{\frac{-c}{2}}\otimes cc_W(h^\ast\mathcal{F})+
\varepsilon_W(h^\ast\mathcal{F})=(-1)^c\cdot c(T_W^\ast X)^{-1}\cap
h^\ast\varepsilon_X(\mathcal{F})
\end{equation*}
in $\Theta_k\otimes {\rm CH}_\bullet(W)$. Here $cc_W$ in the left hand side denotes the characteristic class defined in \cite[Definition 6.7]{Sai17} and $h^\ast\varepsilon_X(\mathcal{F})$ in the right hand side is the image of $\varepsilon_X(\mathcal{F})$ by
$h^\ast\colon\Theta_k\otimes{\rm CH}_\bullet(X)\to
\Theta_k\otimes{\rm CH}_\bullet(W)$. 
\end{lm}
\proof{
As the same proof of \cite[Proposition 7.8]{Sai17} works, we only give a sketch. Let ${\cal E}({\cal F})=\sum_a\xi_a\otimes[C_a]$ denote the epsilon cycle. 
We have a commutative diagram 
\begin{equation*}
\xymatrix{
{\rm CH}_\bullet(X)\ar[r]^-{h^\ast}\ar[d]_-{\cong}&{\rm CH}_\bullet(W)\ar[d]_-\alpha\\
{\rm CH}_n(\mathbb{P}(T^\ast X\oplus\mathbb{A}^1_X))\ar[r]^-{h^\ast}&{\rm CH}_m(\mathbb{P}((W\times_XT^\ast X)\oplus\mathbb{A}^1_W))
}
\end{equation*}
where $\alpha$ is the isomorphism given by $(a_i)\mapsto\sum_ic_1(\mathcal{O}(1))^i\cap \pi^\ast a_i$. 
Hence the cycle class $\alpha(h^\ast\varepsilon_X({\cal F}))$ is represented by the cycle $\sum_a\xi_a\otimes[\mathbb{P}(C_a\oplus\mathbb{A}^1_W)\times_XW]$, where, for a possibly non-reduced scheme $Y$, $[Y]$ denotes the linear combination $\sum_im_i[D_i]$ of the irreducible components $D_i$ of $Y$ with $m_i$ being the multiplicities of $D_i$ in $Y$. 

Let $h^!{\cal E}({\cal F})=\sum_b\eta_b\otimes[C'_b]$ be the pull-back; for the definition of $h^!$,  see Definition \ref{!plbk}. Let 
$\beta$ denote the isomorphism ${\rm CH}_\bullet(W)\to{\rm CH}_m(\mathbb{P}(T^\ast W\oplus\mathbb{A}^1_W))$ considered in (\ref{chowpr}). 
The same computation as done in \cite[Lemma 6.5]{Sai17} shows that the cycle $\sum_b\eta_b\otimes[\mathbb{P}(C'_b\oplus\mathbb{A}^1_W)]$ represents the cycle class $(-1)^c\cdot\beta(c(T^\ast_WX)^{-1}\cap h^\ast\varepsilon_X({\cal F}))$ in ${\rm CH}_m(\mathbb{P}(T^\ast W\oplus\mathbb{A}^1_W))$; the sign $(-1)^c$ appears here since the definition of $h^!$ involves the sign. Then the 
 assertion follows from Theorem \ref{prtr}. 
\qed
}

For the theory of universal family of pencils, we follow the notation in Subsection 2.3. 
Let $\Proj=\Proj^n$ be the projective $n$-space and $\Proj^\vee$ be its dual projective space.  
We identify 
$Q=\Proj(T^*\Proj)$
and let
${\bm p}\colon Q\to \Proj$ and
${\bm p}^\vee\colon Q\to \Proj^\vee$
denote the projections.
The Radon transforms
$R=R{\bm p}^\vee_!
{\bm p}^*[n-1]$ and 
$R^\vee=R{\bm p}_!
{\bm p}^{\vee*}[n-1](n-1)$
induce group homomorphisms
\begin{equation}
R\colon
K(\Proj,\Zl)
\to
 K(\Proj^\vee,\Zl),
\quad
R^\vee\colon
 K(\Proj^\vee,\Zl)
\to
 K(\Proj,\Zl).
\label{dfKR}
\end{equation}
Let $\Tilde{\Theta}_k:={\rm Hom}_{\rm conti}(
G_k^{ab},\Zl^\times)$. For a proper variety $X$ over $k$ and ${\cal F}\in D^b_c(X,\Zl)$, we put $\chi({\cal F}):=\chi(X_{\bar{k}},{\cal F})$ and $\varepsilon^{-1}({\cal F}):=\det(R\Gamma(X_{\bar{k}},{\cal F}))$ to ease notations. 
The pair $(\chi, \varepsilon^{-1})$ induces 
 a group homomorphism $(\chi,\varepsilon^{-1})\colon 
 K(X,\Zl)
\to \mathbb{Z}\times\Tilde{\Theta}_k$. 

\begin{lm}\label{lmPn}
Let $n\geqq 1$ be an integer
and $\Proj=\Proj^n$.
\begin{enumerate}
\item
The diagram
\begin{equation}
\xymatrix{
K(\Proj,\Zl)\ar[r]^{(\chi,\varepsilon^{-1})}\ar[d]_R&
\mathbb{Z}\times\Tilde{\Theta}_k\ar[d]\\
K(\Proj^\vee,\Zl)\ar[r]_{(\chi,\varepsilon^{-1})}&
\mathbb{Z}\times\Tilde{\Theta}_k
}
\end{equation}
is commutative where the right vertical arrow sends $(a,b)$ to 
$((-1)^{n-1}na,\chi_\cyc^{(-1)^{n}\cdot\frac{n(n-1)}{2}a}b^{(-1)^{n-1}n})$. 
\item
The composition of $
K(\Proj,\Zl)\xrightarrow{R^\vee R}
 K(\Proj,\Zl)\xrightarrow{(\chi,\varepsilon^{-1})}\mathbb{Z}\times\Tilde{\Theta}_k$
maps $\mathcal{F}$ to
$(n^2\chi(\mathcal{F}),\varepsilon^{-1}(\mathcal{F})^{n^2})$.
\end{enumerate}
\end{lm}

\proof{
1. 
For  $\mathcal{F}\in D^b_c(\Proj,\Zl)$,
we have
\begin{equation*}
R\Gamma(\Proj^\vee_{\bar{k}},R{\cal F})
=
R\Gamma(Q_{\bar{k}},{\bm p}^\ast{\cal F})[n-1]
=
R\Gamma(\Proj_{\bar{k}},\mathcal{F}\otimes^LR{\bm p}_\ast
\Zl)[n-1]
\end{equation*}
by the projection formula.
Hence the assertion follows from
$R^q{\bm p}_\ast\Zl
=\Zl(-q/2)$ when $q$ is an even integer with $0\leqq q\leqq 2(n-1)$ 
and $=0$ otherwise.

2. Similarly as 1, for $\mathcal{G}\in D^b_c(\Proj^\vee,\Zl)$, we have 
\begin{equation*}
(\chi,\varepsilon^{-1})R^\vee\mathcal{G}=
((-1)^{n-1}n\chi(\mathcal{G}),\chi_\cyc^{(-1)^{n}\frac{n(n-1)}{2}\chi(\mathcal{G})}
\varepsilon^{-1}(\mathcal{G}(n-1))^{(-1)^{n-1}n}). 
\end{equation*}
Combining this, the assertion 1, and the equality $\varepsilon^{-1}(\mathcal{G}(n-1))
=\chi_\cyc^{(n-1)\chi(\mathcal{G})}\varepsilon^{-1}(\mathcal{G})$, 
we get the assertion. 
\qed}
\smallskip

We define the Legendre transform
\begin{equation}
L\colon
{\rm CH}_\bullet (\Proj)
\to
{\rm CH}_\bullet (\Proj^\vee)
\label{eqbarL}
\end{equation}
by setting $L(a)={\bm p}^\vee_\ast(c(T^\ast(Q/\Proj))\cap{\bm p}^\ast a)$. Here 
$c(T^\ast(Q/\Proj))\in{\rm CH}_\bullet(Q)$ is the total Chern class of the relative 
cotangent bundle of $Q/\Proj$ and ${\bm p},{\bm p}^\vee$ are the projections $Q\to\mathbb{P},\Proj^\vee$. Similarly, we define $L^\vee\colon {\rm CH}_\bullet (\Proj^\vee)
\to
{\rm CH}_\bullet (\Proj)$ by setting $L^\vee(a)={\bm p}_\ast(c(T^\ast(Q/\Proj^\vee))\cap{\bm p}^{\vee\ast} a)$. 
We also let $L$ and $L^\vee$ denote the induced maps
$\Theta_k\otimes{\rm CH}_\bullet(\Proj)\to
\Theta_k\otimes{\rm CH}_\bullet(\Proj^\vee)$ and $\Theta_k\otimes{\rm CH}_\bullet(\Proj^\vee)\to
\Theta_k\otimes{\rm CH}_\bullet(\Proj)$.

\begin{pr}\label{prRR}
Let $n\geqq 1$ be an integer
and $\Proj=\Proj^n$. Let $cc_X$ denote the map of characteristic classes defined in \cite[Definition 6.7]{Sai17}. 
\begin{enumerate}
\item 
The diagram
\begin{equation}
\xymatrix{
 K(\Proj,\Zl)
\ar[rr]^{(cc_\Proj,\varepsilon_\Proj)\ \ \ \ \ \ \ \ \ \ }\ar[d]_R&&
(\mathbb{Z}\oplus\Theta_k)\otimes{\rm CH}_\bullet (\Proj)\ar[d]^{\Tilde{L}}
\\
K(\Proj^\vee,\Zl)
\ar[rr]^{(cc_{\Proj^\vee},\varepsilon_{\Proj^\vee})\ \ \ \ \ \ \ \ \ \ }&&
(\mathbb{Z}\oplus\Theta_k)\otimes{\rm CH}_\bullet (\Proj^\vee)
}
\label{eqR}
\end{equation}
is commutative, where $\Tilde{L}$ is given by 
$\Tilde{L}(a,b)=(L(a),L(\chi_\cyc^{\frac{1-n}{2}a}\cdot b))$. 
The diagram with $R$ replaced  by $R^\vee$ and $\Tilde{L}$ by 
$\Tilde{L}^\vee\colon(a,b)\mapsto(L^\vee(a),L^\vee(\chi_\cyc^{\frac{n-1}{2}a}\cdot b))$ is commutative. 
\item
The diagram
\begin{equation}
\xymatrix{
K(\Proj,\Zl)\ar[rrd]_{(\chi,\varepsilon^{-1})}
\ar[rr]^{\!\!\!\!\!\!\!\!\!
(cc_\Proj,\varepsilon_\Proj)}&&
(\mathbb{Z}\oplus\Theta_k)\otimes{\rm CH}_\bullet (\Proj)
\ar[d]^{\deg}\\
&&\mathbb{Z}\oplus\Theta_k
}
\label{eqchP}
\end{equation}
is commutative. 
\end{enumerate}
\end{pr}

\proof{
We prove the assertions by 
induction on $n$.
When $n=1$, the projections
${\bm p}\colon Q\to \Proj$
and
${\bm p}^\vee\colon Q\to \Proj^\vee$
are isomorphisms
and the assertion 1 is obvious.
Since $\deg cc_\Proj{\cal F}
=
(CC{\cal F},T^*_\Proj\Proj)_
{T^*\Proj}$ (\cite[Lemma 6.9.1]{Sai17}) and $\deg\varepsilon_\Proj\mathcal{F}=
(\mathcal{E}(\mathcal{F}),T^*_\Proj\Proj)_
{T^*\Proj}$ (Lemma \ref{degepcls}.1),
the assertion 2 for $n=1$ is nothing but
the Grothendieck--Ogg--Shafarevich
formula and the product formula (Lemma \ref{prodcur}). 

From now on, suppose that $n\geq2$. We show that the assertion 2 for $n-1$
implies the assertion 1 for $n$. As the second part of the assertion 1 
follows from the first one, 
we prove the commutativity of (\ref{eqR}). Put $\theta:=\Tilde{L}\circ(cc_{\Proj},\varepsilon_{\Proj})-
(cc_{\Proj^\vee},\varepsilon_{\Proj^\vee})\circ R$. 
We show $\theta=0$ 
by using the direct sum decomposition 
\begin{equation}
\xymatrix{
{\rm CH}_\bullet(\Proj^\vee)
\cong
{\rm CH}_n(\Proj(T^*\Proj^\vee
\oplus {\mathbb{A}}^1_{\Proj^\vee}))
\ar[r]&
{\rm CH}_{n-1}(\Proj(T^*\Proj^\vee))
\oplus
{\rm CH}_n(\Proj^\vee)=
{\rm CH}_{n-1}(Q)
\oplus {\mathbb{Z}}, 
}
\label{eqCHQZ}
\end{equation}
where the map to the first component is the pull-back induced from $\Proj(T^*\Proj^\vee)\hookrightarrow\Proj(T^*\Proj^\vee
\oplus {\mathbb{A}}^1_{\Proj^\vee})$ and the map to the second one is the projection to the dimension $n$-part. 

First, we check that 
the composition of $\theta$ with the projection ${\rm CH}_\bullet(\Proj^\vee)\to{\rm CH}_{n-1}(Q)$ is zero. Consider the diagram 
\begin{equation*}
\xymatrix{
{\rm CH}_\bullet(\Proj)\ar[r]\ar[d]_-{{\bm p}^\ast}&{\rm CH}_{n-1}(\Proj(T^\ast \Proj))\ar[d]_-{{\bm p}^\ast}&\\
{\rm CH}_\bullet(Q)\ar[d]_-{c(T^\ast(Q/\Proj))\cap}\ar[r]&{\rm CH}_{2n-2}(\Proj(Q\times_{\Proj}T^\ast\Proj))\ar[d]_-{d{\bm p}_\ast}&\\{\rm CH}_\bullet(Q)\ar[d]_-{{\bm p}^\vee_\ast}\ar[r]&{\rm CH}_{2n-2}(\Proj(T^\ast Q))\ar[r]^-{d{\bm p}^{\vee\ast}}&{\rm CH}_{n-1}(\Proj(Q\times_{\Proj^\vee}T^\ast\Proj^\vee))\ar[d]_-{{\bm p}^\vee_\ast}\\
{\rm CH}_\bullet(\Proj^\vee)\ar[rr]&&{\rm CH}_{n-1}(\Proj(T^\ast\Proj^\vee)). 
}
\end{equation*}
Here the horizontal arrows without index are the maps sending $a_i\in{\rm CH}_i$ to $c_1({\cal O}(1))^i\cap \pi^\ast a_i$, where $\pi$ stands for the projections of projective space bundles to base schemes. This diagram is commutative: the commutativity of the middle square comes from Lemma \ref{Chowprojlm}.2.(b). That of the other squares follows from the projection formula. Then 
 Corollary \ref{CCR} and \cite[Corollary 7.5]{Sai17} show that the composition $K(\Proj,\Zl)\xrightarrow{\theta}(\mathbb{Z}\oplus\Theta_k)\otimes{\rm CH}_\bullet(\Proj^\vee)\to
 (\mathbb{Z}\oplus\Theta_k)\otimes{\rm CH}_{n-1}(\Proj(T^\ast\Proj^\vee))$ is zero. 
 
Next, we show that the composition of $\theta$ with the
 projection $
(\mathbb{Z}\oplus\Theta_k)\otimes{\rm CH}_\bullet(\Proj^\vee)
\to {\mathbb{Z}}\oplus\Theta_k$
to the dimension $n$-part is zero. 
Let ${\cal F}\in D^b_c(\Proj,\Zl)$ 
and put $C=SS({\cal F})$. For a hyperplane section  $h\colon H\to \Proj$, $h$ is $C$-transversal if (and only if) the $k$-rational point $i\colon {\rm Spec}(k)
\to \Proj^\vee$ corresponding to $H$ is not contained in the the image of $\Proj(C)\to\Proj^\vee$. 
Therefore, 
after replacing $k$ by 
a finite extension, we can find $H$ so that $h$ is $C$-transversal and the base of any irreducible component of $C$ is not contained in $H$: hence $h$ is properly $C$-transversal. Note that the projection ${\rm CH}_\bullet(\Proj^\vee)
\to{\rm CH}_n(\Proj^\vee)=\mathbb{Z}$ coincides with $i^\ast\colon{\rm CH}_\bullet(\Proj^\vee)
\to{\rm CH}_\bullet({\rm Spec}(k))\cong\mathbb{Z}$. 

By the hypothesis of the induction,
we have
$\deg cc_Hh^*{\cal F}=
\chi (h^*{\cal F})$ and $\deg\varepsilon_Hh^\ast\mathcal{F}
=\varepsilon(h^\ast\mathcal{F})^{-1}$. 
By \cite[Proposition 7.8]{Sai17},
we have
$cc_Hh^*{\cal F}=
-c(\mathcal{O}_H(-1))^{-1}\cap h^\ast cc_\Proj{\cal F}$ and by Lemma \ref{trprepcl}, 
$\chi_\cyc^{\frac{1}{2}}\otimes cc_H(h^\ast\mathcal{F})-\varepsilon_H
(h^\ast\mathcal{F})=c(\mathcal{O}_H(-1))^{-1}\cap h^\ast 
\varepsilon_\Proj({\cal F})$. 
Pulling back the short exact sequence 
\begin{equation*}
0\to T^\ast_Q(\Proj\times\Proj^\vee)\to T^\ast\Proj^\vee\times_{\Proj^\vee}
Q\to T^\ast(Q/\Proj)\to0
\end{equation*}
by $H\hookrightarrow Q$, we have 
\begin{equation*}
c(T^\ast(Q/\Proj)\times_QH)=c(T_Q^\ast(\Proj\times\Proj^\vee)\times_QH)^{-1}
=c(\mathcal{O}_H(-1))^{-1}. 
\end{equation*}
Hence we have a commutative diagram
\begin{equation}
\xymatrix{
{\rm CH}_\bullet(\Proj)\ar[r]^-{h^\ast}\ar[d]_L&
{\rm CH}_\bullet(H)\ar[d]^-{\deg(c(\mathcal{O}_H(-1))^{-1}\cap-)}\\
{\rm CH}_\bullet(\Proj^\vee)\ar[r]^-{i^\ast}&\mathbb{Z}. 
}
\label{eqH}
\end{equation}
Thus we get 
\begin{align*}
i^\ast\Tilde{L}(cc_\Proj,\varepsilon_\Proj)\mathcal{F}
&=i^\ast(L(cc_\Proj(\mathcal{F})),L(\chi_\cyc^{\frac{1-n}{2}}\otimes 
cc_\Proj(\mathcal{F})+
\varepsilon_\Proj(\mathcal{F})))\\&=(\deg(c(\mathcal{O}_H(-1))^{-1}\cap h^\ast cc_\Proj(\mathcal{F})),\deg(c(\mathcal{O}_H(-1))^{-1}\cap h^\ast(\chi_\cyc^{\frac{1-n}{2}}\otimes 
cc_\Proj(\mathcal{F})+
\varepsilon_\Proj(\mathcal{F}))))\\&=
(-\deg cc_H(h^\ast\mathcal{F}),\deg(\chi_\cyc^{\frac{n-1}{2}}\otimes cc_H(h^\ast\mathcal{F})+
 \chi_\cyc^{\frac{1}{2}}\otimes cc_H(h^\ast\mathcal{F})-\varepsilon_H(
h^\ast\mathcal{F})))\\&=
(-\chi(h^\ast\mathcal{F}),\chi_\cyc^{\frac{n}{2}\chi(h^\ast\mathcal{F})}\cdot
\det R\Gamma(H_{\bar{k}},h^\ast\mathcal{F})^{-1}). 
\end{align*}
On the other hand, 
Lemma \ref{degepcls}.2 and \cite[Lemma 6.9.2]{Sai17} show 
\begin{equation*}
i^\ast(cc_{\Proj^\vee},\varepsilon_{\Proj^\vee})R\mathcal{F}
=((-1)^n{\rm rk}^\circ R\mathcal{F},\det(R\mathcal{F})^{\circ
(-1)^n}\cdot \chi_\cyc^{(-1)^{n+1}\frac{n}{2}\cdot{\rm rk}^\circ R\mathcal{F}
}). 
\end{equation*}
Hence the assertion 1 follows.

We show that the assertion 1 for $n\geqq 2$
implies the assertion 2 for $n$.
By the commutative diagrams
(\ref{eqR}),
the endomorphism $R^\vee R$
of $K(\Proj,\Zl)$
preserves the kernel $K(\Proj,\Zl)^\circ$ of $(cc_\Proj,\varepsilon_\Proj)
\colon K(\Proj,\Zl)
\to (\mathbb{Z}\oplus\Theta_k)\otimes
{\rm CH}_\bullet(\Proj)$. 
Take an element $\mathcal{F}$ of $K(\Proj,\Zl)^\circ$. 
There is an element $\mathcal{G}\in K({\rm Spec}(k),\Zl)$ such that 
$a^\ast\mathcal{G}=R^\vee R\mathcal{F}-\mathcal{F}$, where 
$a\colon\Proj\to{\rm Spec}(k)$ is the structure morphism. 
Since $cc_\Proj(a^\ast\mathcal{G})=0$ and $\varepsilon
_\Proj(a^\ast\mathcal{G})=1$, we know that ${\rm rk}\mathcal{G}=0$ and 
$\det(\mathcal{G})=1$. Hence we get 
\begin{equation*}
(\chi,\varepsilon^{-1})R^\vee R\mathcal{F}=(\chi,\varepsilon^{-1})\mathcal{F}, 
\end{equation*}
which is equivalent to 
$(n^2\chi(\mathcal{F}),\varepsilon(\mathcal{F})^{-n^2})=
(\chi(\mathcal{F}),\varepsilon(\mathcal{F})^{-1})$ by Lemma 
\ref{lmPn}.2. This means that $\mathcal{F}$ is contained in the kernel 
of $(\chi,\varepsilon^{-1})$. 
We show that the cokernel of $(cc_\Proj,\varepsilon_\Proj)$ is torsion. Indeed, 
for a linear subspace $\Proj^a\subset\Proj$ and a continuous character $\xi\colon G_k\to\Zl^\times$, let $\xi_{\Proj^a}$ denote the smooth $\Zl$-sheaf of rank $1$ on $\Proj^a$ corresponding to the character $\pi_1(\Proj^a)\to G_k\to\Zl^\times$, which we regard as an element of $D^b_c(\Proj,\Zl)$ by $0$-extension.  We have ${\cal E}(\xi_{\Proj^a}[a])=\xi\cdot\chi_{\rm cyc}^{-\frac{a}{2}}\otimes[T_{\Proj^a}^\ast \Proj]$. Since the classes of $\Proj(T^\ast_{\Proj^a}\Proj\times_k\mathbb{A}^1_k)$ span ${\rm CH}_n(\Proj(T^\ast\Proj\otimes\mathbb{A}^1_{\Proj}))\cong {\rm CH}_\bullet(\Proj)$, the assertion follows from Lemma \ref{theta}.2.  

Thus there exists a unique group homomorphism $\deg'\colon(\mathbb{Z}\oplus\Theta_k)\otimes{\rm CH}_\bullet (\Proj)\to\mathbb{Q}\oplus\Theta_k$ which 
makes the diagram (\ref{eqchP}) with $\deg$ replaced  by $\deg'$ commutative. 
We need to show the equality $\deg=\deg'$. But $(\mathbb{Q}\oplus\Theta_k)\otimes{\rm CH}_\bullet(\Proj)$ is spanned as a $\mathbb{Q}$-vector space by the images of $\xi_{\Proj^a}[a]$ for $\xi\colon G_k\to\Zl^\times$ and linear subspaces $\Proj^a\subset\Proj$, it suffices to check the equality for $\xi_{\Proj^a}[a]$, which is done by using Lemma \ref{degepcls} and \cite[Lemma 6.9.1]{Sai17}. 
\qed}

\begin{cor}\label{corLR}
Let ${\cal F}\in D(\Proj,\Zl)$
be a constructible complex
of $\Zl$-sheaves on $\Proj=\Proj^n$.
Then, for the Radon transform
$R{\cal F}$, we have
\begin{equation}
\mathcal{E}(R{\cal F})=
L(\mathcal{E}(\mathcal{F})(\frac{1-n}{2})).
\label{eqCCL}
\end{equation}
\end{cor}

\proof{
Except for the coefficient
of the $0$-section, it is proved in
Corollary \ref{CCR}.
Since the coefficient
of the $0$-section
is given by 
${\rm id}\otimes\deg\colon \Theta_k\otimes 
{\rm CH}_\bullet(\Proj^\vee)
\to \Theta_k$,
it follows from
Proposition \ref{prRR}.1.
\qed}

Here is the product formula of the global epsilon factors. 
\begin{thm}\label{thmEP}
Let $k$ be the perfection of a finitely generated field. 
Let $X$ be a projective smooth
variety over $k$.
Then, for $\mathcal{F}\in D(X,\Zl)$, we have
\begin{equation}
\det(R\Gamma(X_{\bar{k}},{\cal F}))
=
(\mathcal{E}({\cal F}),
T^{*}_{X}X)_{T^{*}X}
\label{eqEP}
\end{equation}
as elements of $\Theta_k$. 
\end{thm}

\proof{
Since $X$ is assumed projective,
it follows from Lemma \ref{lm}.4 
and Proposition \ref{prRR}.2. 
\qed}

\begin{cor}
Let $X$ be a projective smooth variety over a finite field $\F_q$. Let $K$ be a field endowed with a valuation $\lvert-\rvert\colon K\to\mathbb{R}_{\geq0}$. 
Let $\iota\colon\Ql\to K$ be a field homomorphism. The composition of $\Ql^\times\xrightarrow{\iota}K^\times\xrightarrow{\lvert-\rvert}\mathbb{R}^\times$ induces a group homomorphism $\Ql^\times/\mu\to\mathbb{R}^\times$, for which we write $x\mapsto \lvert\iota(x)\rvert$ by abuse of notation. 

Let $\mathcal{F}$ be 
an element of $D^b_c(X,\Zl)$ and $\mathcal{E}(\mathcal{F})=
\sum_a\beta_a\otimes[C_a]\in\Zl^\times/\mu\otimes Z_n(T^\ast X)$ be the epsilon cycle. 
Here we identify $\Theta_{\F_q}$ and $\Zl^\times/\mu$ via $\xi\mapsto \xi({\rm Frob}_q)$. 
Then we have 
\begin{equation}\label{prodform}
|\iota(\det({\rm Frob}_q,R\Gamma(X_{\mathbb{F}_q},{\cal F})))|=\prod_a
|\iota(\beta_a)|^{\deg(C_a,T^\ast_XX)_{T^\ast X}}. 
\end{equation}
\end{cor}
\qed

\begin{ex}\label{p-adicval}
Let $\F_q$ be a finite field with $q$ elements. Let $X$ be a projective 
smooth scheme over $k$. Let $\mathcal{F}\in D^b_c(X,\Zl)$ be a 
constructible complex on $X$. 
\begin{enumerate}
\item Assume that ${\cal F}_{\Ql}:=\mathcal{F}\otimes_{\Zl}\Ql$ is $\iota$-pure 
of $\iota$-weight $0$ (in the sense of \cite[II.12.7]{KW})  
for a field isomorphism $\iota\colon\Ql\to\mathbb{C}$. 
Then we know that the absolute values $|\iota(\alpha)|$ of 
the eigenvalues $\alpha$ of the geometric Frobenius on ${\rm H}^i
(X_{\overline{\F}_q},\mathcal{F}_{\Ql})$ are equal to $q^{\frac{i}{2}}$. Hence 
the product formula (\ref{prodform}) gives an expression of the weighted Euler--Poincar\'e 
characteristic $\frac{1}{2}\sum_i(-1)^ii\cdot{\rm dim}{\rm H}^i(X_{\overline{\F}_q},\mathcal{F}_{\Ql})$ 
as the intersection number $({\rm log}_q|\iota\mathcal{E}(\mathcal{F})|,
T^\ast_XX)_{T^\ast X}$. 
\item Let $\iota\colon\Ql\to\overline{\mathbb{Q}}_p$ be a field  
isomorphism (where $p$ is the characteristic of $\mathbb{F}_q$). Then the product formula (\ref{prodform}) 
gives an expression of the $p$-adic 
valuation of the global epsilon factor $\varepsilon(X,\mathcal{F})=\det(-{\rm Frob}_q,R\Gamma(X_{\overline{\mathbb{F}}_q},{\cal F}_{\Ql}))^{-1}$ 
in terms of those of local epsilon factors. The $p$-adic valuation of 
the local epsilon factors of tamely ramified representations can be 
computed by Stickelberger's theorem (\cite[Proposition 6.13]{CF}), which is suggested to the 
author by N. Katz. 
\end{enumerate}
\end{ex}

\subsection{An axiomatic description of epsilon cycles}

We give an axiomatic description of epsilon cycles. 
A similar description of characteristic cycles is 
considered in \cite[Proposition 8]{properp}. 

\begin{thm}\label{epcygen}
Let $k$ be the perfection of a finitely generated field of characteristic $p\neq\ell$. 
There exists a unique way to attach, for pairs $(X,\mathcal{F})$ where 
$X$ is a smooth scheme of finite type over $k$ and 
$\mathcal{F}\in D^b_c(X,\Zl)$, 
with a cycle $\mathcal{E}(\mathcal{F})=\sum_a\xi_a\otimes[C_a]$ with 
$\Theta_k$-coefficient and supported on the singular support 
$SS(\mathcal{F})$ which should satisfy the following axioms.
\begin{enumerate}
\item (Normalization) Let $X={\rm Spec}(k')$ be the spectrum of a finite 
extension $k'$ of $k$. 
Then we have 
\begin{equation*}
\mathcal{E}(\mathcal{F})=(\det(\mathcal{F})\circ {\rm tr}_{k'/k})^{\frac{1}{\deg(k'/k)}}
\otimes[T^\ast_XX].
\end{equation*}
\label{Norm}
\item (Tate Twist) We have 
\begin{equation*}
\mathcal{E}(\mathcal{F}(1))=
\chi_\cyc\otimes 
CC(\mathcal{F})+\mathcal{E}(
\mathcal{F}). 
\end{equation*}
\label{Twist}
\item (Additivity) For a distinguished triangle 
\begin{equation*}
\mathcal{F}'\to\mathcal{F}\to\mathcal{F}''\to,
\end{equation*}
we have $\mathcal{E}(\mathcal{F})=\mathcal{E}(\mathcal{F}')+
\mathcal{E}(\mathcal{F}'')$. \label{Mult}
\item (Closed Immersion) For a closed immersion $i\colon X\to P$ of smooth 
$k$-schemes of finite type and $\mathcal{F}\in D^b_c(X,\Zl)$, we have 
$\mathcal{E}(i_\ast\mathcal{F})=i_!\mathcal{E}(\mathcal{F})$. 
\label{Clim}
\item (Pull-Back) For a half integer $r\in\frac{1}{2}\mathbb{Z}$, put $\mathcal{E}(\mathcal{F})(r):=\chi_{{\rm cyc}}^r\otimes CC({\cal F})+
\mathcal{E}(\mathcal{F})$. 
For a properly 
$SS(\mathcal{F})$-transversal morphism $h\colon W\to X$ from a smooth 
$k$-scheme $W$ of finite type, we have 
\begin{equation*}
\mathcal{E}(h^\ast\mathcal{F})=h^!(\mathcal{E}(\mathcal{F})(\frac{\dim X-\dim W}{2})). 
\end{equation*}
Here ${\rm dim}X$ and ${\rm dim}W$ are regarded as  locally constant functions on 
$X$ and $W$. \label{Pullb}
\item (Radon Transform) For a constructible complex $\mathcal{F}\in 
 D^b_c(\Proj,\Zl)$ on a projective space $\Proj=\Proj^n$, we have 
\begin{equation*}
\mathcal{E}(R\mathcal{F})=L(\mathcal{E}(\mathcal{F})(\frac{1-n}{2})). 
\end{equation*}\label{Radtr}
\item (Same Monodromy) \label{Samem}
Let $X$ (resp. $X'$) be a smooth curve over $k$ 
and $x$ (resp. $x'$) be a closed point of $X$ (resp. $X'$). 
Let $\mathcal{F}$ (resp. $\mathcal{F}'$) be an element of 
$D^b_c(X,\Zl)$ (resp. $D^b_c(X',\Zl)$). Assume that there exists an 
isomorphism $f\colon X_{(x)}\xrightarrow{\cong}X'_{(x')}$ of 
$k$-schemes between the henselizations such that the complexes 
$\mathcal{F}|_{X_{(x)}}$ and $f^\ast\mathcal{F}'|_{X'_{(x')}}$ are 
isomorphic. 
Then the coefficient of $[T^\ast_xX]$ 
in $\mathcal{E}(\mathcal{F})$ coincides with 
that of $[T^\ast_{x'}X']$ in $\mathcal{E}(\mathcal{F}')$. 
\end{enumerate}
\end{thm}
To prove the theorem, we give some lemmas. 
\begin{lm}\label{unip}
Let $\mathcal{E}(-)$ be an assignment as in Theorem \ref{epcygen} 
satisfying the axioms there. 
Let $X$ be a smooth curve  of finite type over $k$ and $x\in X$ be a closed point. 
Denote by $U$ the complement of $x$ in $X$. Let $\mathcal{F}$ be a 
smooth $\Zl$-sheaf on $U$. Assume that $\mathcal{F}\otimes_{\Zl} 
\Ql$ 
has a unipotent monodromy at $x$. Then the coefficient of $[T^\ast_xX]$ in 
$\mathcal{E}(j_!\mathcal{F}) $ equals to $(\det(\mathcal{F})_x^{-1}
\circ {\rm tr}_{k(x)/k})^{\frac{1}{\deg(x/k)}}$ where 
$j\colon U\to X$ is the immersion. Note that we can extend $\det(\mathcal{F})$ 
to $X$ smoothly since this is unramified at $x$. 
\end{lm}
\proof{
Let $k'$ be the residue field at $x$. We regard $\mathbb{A}^1_{k'}$ as 
a smooth $k$-scheme. Fix an isomorphism 
$f\colon X_{(x)}\cong \mathbb{A}^1_{k',(0)}$ of $k$-schemes. 
By the Gabber--Katz extension \cite[THEOREM 1.5.6]{GK}, there exists a smooth $\Zl$-sheaf $\mathcal{G}$ on 
$\mathbb{G}_{{\rm m},k'}$ such that $\mathcal{F}|_{\eta_x}$ and $f^\ast\mathcal{G}|_{\eta_0}$ are isomorphic and $\mathcal{G}$ is 
tamely ramified at $\infty$, where $\eta_x$ and $\eta_0$ are the generic points of 
$X_{(x)}$ and $\mathbb{A}^1_{k',(0)}$. Since the monodromy 
of $\mathcal{G}$ at $0$ is unipotent, the semi-simplification of $
\mathcal{G}\otimes_{\Zl}\Ql$ is unramified at $0$ and $\infty$. Hence the assertion 
follows from the axioms (\ref{Norm}), (\ref{Mult}), (\ref{Clim}), and 
(\ref{Samem}). 
\qed
}

\begin{lm}\label{prodgen}
Let $\mathcal{E}(-)$ be an assignment as in Theorem \ref{epcygen} 
satisfying the axioms there. 
Let $X$ be a projective smooth scheme over $k$. Let $\mathcal{F}$ 
be an element of $D^b_c(X,\Zl)$. Then we have an equality 
\begin{equation*}
\det(R\Gamma(X_{\bar{k}},\mathcal{F}))=
(\mathcal{E}(\mathcal{F}),T^\ast_XX)_{T^\ast X}
\end{equation*}
of elements of $\Theta_k$. 
\end{lm}
\proof{
Since $X$ is projective, we may assume that $X=\Proj=\Proj^n\ (n\geq2)$ by the 
axiom (\ref{Clim}). Consider the universal hyperplane section $\Proj
\xleftarrow{{\bm p}}Q\xrightarrow{{\bm p}^\vee}\Proj^\vee$. 
Let $R^\vee=R{\bm p}_\ast{\bm p}^{\vee\ast}[n-1](n-1)$ be the 
inverse Radon transform. By \cite[IV. Lemma 1.4]{KW}, we have a distinguished triangle 
\begin{equation*}
\mathcal{F}\to R^\vee R\mathcal{F}\to\oplus_{i=1}^{n-1}
 R\Gamma(\Proj_{\bar{k}},\mathcal{F})[2i](i)\to, 
\end{equation*}
where $R\Gamma(\Proj_{\bar{k}},\mathcal{F})$ is 
regarded as a complex of geometrically constant sheaves. 
By the axioms (\ref{Norm}), (\ref{Mult}), and (\ref{Pullb}), we have 
\begin{equation}\label{RR}
\mathcal{E}(R^\vee R\mathcal{F})-\mathcal{E}(\mathcal{F})
=\det(R\Gamma(\Proj_{\bar{k}},\mathcal{F}))^{(-1)^n(n-1)}\otimes
[T^\ast_\Proj\Proj].  
\end{equation}
On the other hand, by the axioms $(\ref{Radtr})$ and $(\ref{Twist})$, the left hand side of (\ref{RR}) 
equals to 
\begin{align*}
L^\vee(\mathcal{E}(R\mathcal{F}(n-1))(\frac{1-n}{2}))-\mathcal{E}(\mathcal{F})&=
L^\vee L\mathcal{E}(\mathcal{F})-\mathcal{E}(\mathcal{F})\\&
=
(\mathcal{E}(\mathcal{F}),T^\ast_\Proj\Proj)_{T^\ast\Proj}^{(-1)^n(n-1)}
\otimes[T^\ast_\Proj\Proj]. 
\end{align*}
Since $n\geq2$, we have the assertion. 
\qed
}
\vspace{\fill}\\
(Proof of Theorem \ref{epcygen})

First, we show the uniqueness. Let $\mathcal{E}(-)$ be an 
assignment which satisfies the conditions in the statement of the theorem. 
Let $X$ and $\mathcal{F}$ be as in the theorem. We need to 
determine the coefficients of $\mathcal{E}(\mathcal{F})$ uniquely from the 
axioms. By the axiom (\ref{Pullb}), we may assume that $X$ is affine, and by the 
axioms (\ref{Clim}) and (\ref{Pullb}), we may assume that $X$ is projective. 

Let  $i\colon X\hookrightarrow\Proj=\Proj^n$ be a closed immersion. 
Composing $i$ with the Veronese embedding $\Proj\hookrightarrow\Proj'$ of 
$\deg\geq3$ if necessary, we find a line $
 L\to\Proj^\vee_{k'}$ for a finite extension $k'$ of $k$ such that 
 the pair $(f,\pi)$ in the diagram (\ref{Lefpen}) with  
 $X,X\times_{\Proj}Q,\Proj^\vee$ replaced by the base changes by $k\to k'$ is a good pencil (Definition \ref{goodpen}). 
 In the sequel, we regard $k'$-schemes also as $k$-schemes. 
Let $C_a$ be an irreducible component of 
$SS(\mathcal{F}|_{X_{k'}})$. By the definition of a good pencil, there exists a 
closed point $x_a\in X_{k',L}$ of the blow-up $X_{k',L}$ of $X_{k'}$ such that $x_a$ is the unique isolated $SS(\mathcal{F}|_{X_{k'}})$-characteristic point on the fiber $f^{-1}(f(x_a))$ at which $df$ only meets $C_a$ 
and $x_a$ is not contained in the exceptional locus of the blow-up. 
By the axiom (\ref{Radtr}), we have $\mathcal{E}(R(i_\ast\mathcal{F}))
=L(
\mathcal{E}(i_\ast\mathcal{F})(\frac{1-n}{2}))$. Let $\Proj_{k',L}$ be the blow-up of $\Proj_{k'}$ 
along the axis $A_L$ defined by $L$. Since $\pi\colon X_{k',L}\to X$ is 
$SS(\mathcal{F})$-transversal and $X_{k'} $ meets $A_L$ transversally, $\Proj_{k',L}\to\Proj$ is 
$SS(i_\ast\mathcal{F})$-transversal. 
Let $i'\colon L\to\Proj^\vee$ be the composition of 
$L\hookrightarrow\Proj^\vee_{k'}\to\Proj^\vee$. 
Since we have $SS(R(i_\ast\mathcal{F}))\subset
LSS(i_\ast\mathcal{F})\cup T^\ast_{\Proj^\vee}\Proj^\vee$, 
applying \cite[Lemma 3.9.3]{Sai17} to 
the cartesian diagram 
\begin{equation*}
\xymatrix{
\Proj_{k',L}\ar[r]\ar[d]&L\ar[d]^{i'}\\
Q\ar[r]&\Proj^\vee
}
\end{equation*}
shows that $i'$ is properly $SS(R(i_\ast\mathcal{F}))$-transversal. Hence we have 
\begin{align*}
\mathcal{E}(i'^\ast R(i_\ast\mathcal{F}))&=i'^!(\mathcal{E}(R(i_\ast\mathcal{F}))(\frac{n-1}{2}))\\
&=i'^!L\mathcal{E}(i_\ast\mathcal{F}). 
\end{align*}
Let $\xi_a$ be the coefficient of $[C_a]$ in ${\cal E}({\cal F})$ and let $v_a=f(x_a)$ be the image of $x_a$ in $L$. Then the coefficient of $[T^\ast_{v_a}L]$ in $i'^!L\mathcal{E}(i_\ast\mathcal{F})$ is equal to $\xi_a^{(-1)^{n-1}\cdot(C_a,df)_{x_a}}$. 
In this way, we reduce to the case when $X$ is a projective 
smooth curve. 

Suppose that $X$ is a projective smooth curve. 
By the axioms (\ref{Norm}) and (\ref{Clim}), 
we may assume $\mathcal{F}=j_!\mathcal{G}$ where 
$j\colon U\to X$ is an open immersion from an open dense subset $U$ and 
$\mathcal{G}$ is a smooth $\Zl$-sheaf on $U$. The coefficient of 
$[T^\ast_XX]$ can be determined by the axioms (\ref{Norm}) and (\ref{Pullb}). Let $x\in X$ be a closed point not contained in $U$. 
By weak approximation, we can find a finite morphism $f\colon X'\to X$ from 
a projective smooth curve $X'$ such that $f$ is \'etale around $f^{-1}(x)$ 
and $f^\ast\mathcal{G}$ has unipotent monodromy at points of $X'\setminus
f^{-1}(U\cup\{x\})$. Then we can determine the coefficient of $[T^\ast_xX]$ 
by the axiom (\ref{Pullb}), Lemma \ref{unip}, and Lemma \ref{prodgen}. 

We show that epsilon cycles constructed in Theorem \ref{epcygenmil} 
satisfy the axioms. 
The axioms (\ref{Mult}) and (\ref{Samem}) follow from 
the construction. 
The others are proved in Lemma \ref{lm}, Corollary \ref{sm}, 
Theorem \ref{prtr}, and Corollary \ref{corLR}. 
\qed

\begin{rmk}
According to the proof, we can replace the axiom (\ref{Samem}) 
by Lemma \ref{unip}. 
\end{rmk}

\section{Appendix : Complements on $\ell$-adic formalism}\label{adical}

In this appendix, we review the $\ell$-adic formalism on a noetherian topos, following \cite{adic}. This appendix is included since we need to know the explicit definition of $\ell$-adic sheaves in order to complete proofs in the subsections $3.3$ and $3.4$.

To simplify the 
argument, we restrict ourselves to the construction 
of bounded complexes. 
 For a topos $T$, 
we define the category $T^{\mathbb{N}^{\rm op}}$ as follows. The objects are projective systems $(M_n,\varphi_n)_{n\in\mathbb{N}}$ 
indexed by $\mathbb{N}$ where $M_n$ are objects of $T$ and $\varphi_n\colon 
M_{n+1}\to M_n $ are morphisms in $T$, which are referred to as transition maps. If no confusions occur, we omit 
$\varphi_n$ and simply write $(M_n)_{n}$ for it. 
The morphisms $(M_n)_n\to(M'_n)_n$ are families of morphisms 
$M_n\to M'_n$ compatible with the transition maps. 
The category $T^{\mathbb{N}^{\rm op}}$ is known to be  a topos. Let 
\begin{equation}\label{limitT}
\pi\colon T^{\mathbb{N}^{\rm op}}\to T
\end{equation}
be the morphism of topoi given by $\pi_\ast(M_n)_{n}=\varprojlim_n M_n$. The functor 
$\pi^{-1}M=(M,{\rm id}_M)_{n\in\mathbb{N}}$ gives the left adjoint to $\pi_\ast$. 

\begin{df}\label{essdef}
\begin{enumerate}
\item We say that a commutative group object $(M_n)_{n}
$ in $T$ is {\rm essentially zero} if, for every $n\in \mathbb{N}$, there exists $m\geq n$ such that the transition map $M_m\to M_n$ is zero (in \cite[(1.10)]{Jan}, this is called ML-zero). 
\item We say that a complex $K\in D(T^{\mathbb{N}^{\rm op}},\mathbb{Z})$ of 
sheaves of abelian groups is {\rm essentially zero} if each cohomology sheaf of $K$ is 
essentially zero. 
\item We say that a morphism in $D(T^{\mathbb{N}^{\rm op}},\mathbb{Z})$ is an 
{\rm essential isomorphism} if its mapping cone is essentially zero. 
\item We say that a complex $K\in D(T^{\mathbb{N}^{\rm op}},\mathbb{Z})$ 
is {\rm essentially constant} if there exist complexes $L\in 
D(T^{\mathbb{N}^{\rm op}},\mathbb{Z})$ and $M\in D(T,\mathbb{Z})$ together with morphisms 
\begin{equation*}
K\leftarrow L\rightarrow\pi^{-1}M
\end{equation*}
in $D(T^{\mathbb{N}^{\rm op}},\mathbb{Z})$ which are essential isomorphisms. 
\end{enumerate}
\end{df}
We list some of basic properties which are necessary. 
\begin{lm}\label{eses}
\begin{enumerate}
\item Let $M\in D^b(T,\mathbb{Z})$ be a bounded complex. 
Then the canonical morphism $M\to R\pi_\ast\pi^{-1}M$ is an 
isomorphism. 
\item(\cite[Lemma (1.11)]{Jan}) Let $K\in D^b(T^{\mathbb{N}^{\rm op}},\mathbb{Z})$ be 
an essentially zero complex. Then we have $R\pi_\ast K=0$. 
\item (cf. \cite[Lemma 1.3.iv)]{adic}) Let $K\in D^b(T^{\mathbb{N}^{\rm op}},\mathbb{Z})$ be a bounded complex. 
If $K$ is essentially constant, then $R\pi_\ast K$ is bounded and 
the canonical morphism $\pi^{-1}R\pi_\ast K\to K$ is an essential 
isomorphism. 
\end{enumerate}
\end{lm}
\proof{
For a sheaf $N=(N_n)_n$ of abelian groups on $T^{\mathbb{N}^{\rm op}}$ and an 
object $U\in T$, we have a short exact sequence \cite[Proposition (1.6)]{Jan}
\begin{equation}\label{limex}
0\to R^1\varprojlim_n{\rm H}^{i-1}(U,N_n)\to{\rm H}^i(\pi^{-1}(U),N)\to 
\varprojlim_n{\rm H}^{i}(U,N_n)\to 0. 
\end{equation}

1. We may assume that $M$ is a sheaf. Applying (\ref{limex}) to $N=
\pi^{-1}M$, we know that ${\rm H}^i(\pi^{-1}(U),\pi^{-1}M)$ is isomorphic to 
${\rm H}^i(U,M)$, hence the assertion. 

2. We may assume that $K=(K_n)$ is a sheaf. Since $K$ is essentially zero, so is 
$({\rm H}^i(U,K_n))_n$ for $U\in T$. The assertion follows from the exact sequence 
(\ref{limex}). 

3. Take morphisms $K\leftarrow L\to\pi^{-1}M$ as in Definition \ref{essdef}.4. 
Since $K$ is bounded, we may take $L$ and $M$ to be  also bounded. Then, by $2$, the morphisms 
\begin{equation*}
\pi^{-1}R\pi_\ast K\leftarrow \pi^{-1}R\pi_\ast L\to 
\pi^{-1}R\pi_\ast \pi^{-1}M
\end{equation*}
 are isomorphisms. By $1$, the complex in the right hand side is isomorphic to $\pi^{-1}M$, hence the first assertion. 

To prove the second assertion, consider the following commutative diagram 
\begin{equation*}
\xymatrix{
\pi^{-1}R\pi_\ast K\ar[d]&\pi^{-1}R\pi_\ast L\ar[l]_\cong\ar[r]^{\cong\ \ \ \ }
\ar[d]&
\pi^{-1}R\pi_\ast\pi^{-1}M\ar[d]&\pi^{-1}M\ar[l]_{\ \ \ \ \cong}\ar[ld]^{\rm id}\\
K&L\ar[r]\ar[l]&\pi^{-1}M. 
}
\end{equation*}
Since the top horizontal arrows are isomorphisms, the 
assertion follows from a diagram chasing. 
\qed}

Let $(R,\mathfrak{m})$ be a complete discrete valuation ring. 
Put $R_n:=R/\mathfrak{m}^{n+1}$. Let $R_\bullet:=(R_n)_{n\in\mathbb{N}}$ be the ring object of $T^{\mathbb{N}^{\rm op}}$, where the transition maps are given by the natural 
projections $R_{n+1}\to R_n$. The morphism  (\ref{limitT}) is extended to a morphism of ringed topoi 
\begin{equation}\label{ringtop}
\pi\colon(T^{\mathbb{N}^{\rm op}},R_\bullet)\to(T,R), 
\end{equation}
with $\pi^{-1}R\to R_\bullet$ being the canonical one. 
Let $\pi^\ast M:=R_\bullet\otimes_{\pi^{-1}R}\pi^{-1}M$ for a 
sheaf $M$ of $R$-modules on $T$. Since $R_\bullet$ has finite tor-dimension as a $\pi^{-1}R$-module, $L\pi^\ast$ can be defined on  
$D(T,R)$.

For each $n\in\mathbb{N}$, let $i_n\colon T\to T^{\mathbb{N}^{\rm op}}$ be 
the morphism of topoi defined by $i_{n\ast}M=(\cdots M
\xrightarrow{{\rm id}}\cdots \overset{n\text{-th}}M\to\ast\to\cdots\to\ast)$ and 
$i_n^{-1}(M_n)_{n}=M_n$, where $\ast$ is the 
final object of $T$. 
This induces a morphism of ringed 
topoi $i_n\colon(T,R_n)\to(T^{\mathbb{N}^{\rm op}},R_\bullet)$. 
Note that the morphism 
$i^{-1}_nR_\bullet\to R_n$ is an isomorphism. 

\begin{lm}\label{eses1}
\begin{enumerate}

\item Let $M$ be a sheaf of $R_0$-modules on $T$. 
Then the morphisms $\pi^{-1}M\to L\pi^\ast M$ and $L\pi^\ast M\to\pi^{-1}M$ 
are essential isomorphisms. Here the first one is given by $\pi^{-1}M\cong
\pi^{-1}R\otimes_{\pi^{-1}R}^L\pi^{-1}M\to 
R_{\bullet}\otimes_{\pi^{-1}R}^L\pi^{-1}M=L\pi^\ast M$ and the second one is given by 
$L\pi^\ast M\to H^0(L\pi^\ast M)\cong\pi^{-1}M$. 
\item Let $K,L\in D^-(T^{\mathbb{N}^{\rm op}},R_\bullet)$ be 
bounded above complexes. If either of $K$ or $L$ is essentially zero, then 
$L\otimes^L_{R_\bullet}K$ is also essentially zero. 
\item Let $C\in D^-(T^{\mathbb{N}^{\rm op}},R_\bullet)$ be a bounded above 
complex. 
If $R_0\otimes^L_{R_\bullet}C$ is essentially zero, then so is $C$. If $R_0\otimes^L_{R_\bullet}C$ is acyclic, then so is $C$.  
\end{enumerate}
\end{lm}
\proof{
$1$. Let $L_1=(R)_{n}=\pi^{-1}R$ and $L_2=(\mathfrak{m}^{n+1})_n$ be sheaves of $R$-modules 
on $T^{\mathbb{N}^{\rm op}}$ of which the transition maps are the inclusions. 
We have a short exact sequence $0\to L_2\to L_1\to R_\bullet\to0$, which 
gives an $R$-flat resolution of $R_\bullet$. Hence the mapping cone of 
$\pi^{-1}M\to L\pi^\ast M$ is isomorphic to $L_2\otimes_R\pi^{-1}M[1]$. 
Since the transition maps of $L_2\otimes_R\pi^{-1}M$ are zero, the first 
morphism is an essential isomorphism. The assertion for the second one 
follows since the composition of $\pi^{-1}M\to L\pi^\ast M\to\pi^{-1}M$ 
coincides with the identity. 

$2$. This follows from the spectral sequence 
\begin{equation*}
E_2^{p,q}=\oplus_{i+j=q}{\rm Tor}_{-p}^{R_\bullet}({H}^i(L),{ H}^j(K))
\Rightarrow{ H}^{p+q}(L\otimes^L_{R_\bullet}K). 
\end{equation*}

$3$. 
Put $R'_n:=R_\bullet/\mathfrak{m}^{n+1}R_\bullet$ for $n\geq0$. Let $K_n$ be the kernel of 
the natural surjection $R'_{n+1}\to R'_n$. 
Note That $K_n$ is a 
sheaf of $R_0$-modules on $T^{\mathbb{N}^{\rm op}}$ 
essentially isomorphic to $R_0$. 
If $R_0\otimes^L_{R_\bullet}C$ is essentially zero, 
then 
$R'_n\otimes^L_{R_\bullet}C$ is also essentially zero for each $n$ by $2$ and by induction on $n$. Then, for each $n\geq0$,  
there exists $m\geq n$ such that the transition map 
$i^\ast_m{ H}^i(R'_n\otimes_{R_\bullet}^LC)\to i^\ast_n{ H}^i(R'_n\otimes^L_{R_\bullet}C)=
i^\ast_n{ H}^i(C)$ is zero. Hence $i^\ast_m{ H}^i(C)\to
i^\ast_m{ H}^i(R'_n\otimes^L_{R_\bullet}C)\to i^\ast_n{ H}^i(C)$ is zero, which shows the first assertion.

Suppose now that $R_0\otimes^L_{R_\bullet}C$ is acyclic. Then $R'_n\otimes^L_{R_\bullet}C$ is also acyclic by induction on $n$ 
 as $K_n\otimes^L_{R_\bullet}C\cong K_n\otimes_{R_0}^LR_0\otimes^L_{R_\bullet}C\cong 0$. 
 Therefore, $i^\ast_n{ H}^i(C)\cong
i^\ast_n{ H}^i(R'_n\otimes^L_{R_\bullet}C)$ is zero. 
\qed}

We recall the notion 
of (normalized) $R$-complexes, following \cite{adic}. 
\begin{df}\label{normalizedR}
\begin{enumerate}
\item We say that a complex $K\in D^b(T^{\mathbb{N}^{\rm op}},R_\bullet)$ is an 
{\rm $R$-complex} if $L\pi^\ast  R_0\otimes^L_{R_\bullet}K$ is 
essentially constant. 
\item We say that a complex $K\in D^b(T^{\mathbb{N}^{\rm op}},R_\bullet)$ is a 
{\rm normalized $R$-complex} if, for each $n\in\mathbb{N}$, the canonical map 
$i_{n+1}^\ast K\otimes^L_{R_{n+1}}R_n\to i^\ast_nK$ is an isomorphism. 
\end{enumerate}
\end{df}
In the following lemma, we collect key results in order to construct the derived category of $\mathfrak{m}$-adic sheaves. 
\begin{lm}\label{normR}
Let $K\in D^b(T^{\mathbb{N}^{\rm op}},R_\bullet)$ be a complex. 
\begin{enumerate}
\item The canonical map $L\pi^\ast R_0\otimes^L_{R_\bullet}K\to 
\pi^{-1} R_0\otimes^L_{R_\bullet}K$ is an essential isomorphism. 
\item If $K$ is a normalized $R$-complex, then it is an $R$-complex. 
\item If $K$ is an $R$-complex, then $L\pi^\ast R\pi_\ast K$ is bounded and the canonical map 
$L\pi^\ast R\pi_\ast K\to K$ is an essential isomorphism. 
\item The following are equivalent. 
\begin{enumerate}
\item The complex $K$ is a normalized $R$-complex. 
\item 
The canonical morphism $L\pi^\ast R\pi_\ast K\to K$ is an isomorphism. 
\item There is a complex $M\in D(T,R)$ such that 
$L\pi^\ast M\cong K$. 
\end{enumerate}
If these equivalent conditions hold, then there exists a bounded complex $M\in D^b(T,R)$ such that $L\pi^\ast M\cong K$. 
\end{enumerate}
\end{lm}
\proof{
$1$. This follows from Lemma \ref{eses1}.1 and 2. 

$2$. Suppose that $K$ is a normalized $R$-complex. By $1$, the natural map 
$L\pi^\ast R_0\otimes^L_{R_\bullet}K\to 
\pi^{-1} R_0\otimes^L_{R_\bullet}K$ is an essential isomorphism. As $K$ is normalized, the target is isomorphic to $\pi^{-1}i_0^{-1}K$, hence the assertion. 

$3$. First we show that $R_0\otimes^L_RR\pi_\ast K$ is bounded, which implies the boundedness of $L\pi^\ast R\pi_\ast K$. 
Since $R_0$ has finite projective dimension as an $R$-module, we have $R_0\otimes^L_RR\pi_\ast K\cong R\pi_\ast
(L\pi^\ast R_0\otimes^L_{R_\bullet}K)$. Since 
$L\pi^\ast R_0\otimes^L_{R_\bullet}K$ is bounded and essentially constant, 
the assertion follows from Lemma \ref{eses}.3. 

 We show that 
$L\pi^\ast R\pi_\ast K\to K$ is an essential isomorphism. 
By Lemma \ref{eses1}.3, it suffices to show that $R_0\otimes^L_{R_\bullet}
L\pi^\ast R\pi_\ast K\to R_0\otimes^L_{R_\bullet}K$ is an essential 
isomorphism. 
By 1, it is equivalent to showing  that 
\begin{equation*}
L\pi^\ast R_0\otimes^L_{R_\bullet}
L\pi^\ast R\pi_\ast K\to L\pi^\ast R_0\otimes^L_{R_\bullet}K
\end{equation*}
 is an essential 
isomorphism. 
The former complex is isomorphic to 
$L\pi^\ast R\pi_\ast(L\pi^\ast R_0\otimes^L_{R_\bullet}K)$ and 
we have a commutative diagram 
\begin{equation*}
\xymatrix{
L\pi^\ast R\pi_\ast(L\pi^\ast R_0\otimes^L_{R_\bullet}K)\ar[r]&
L\pi^\ast R_0\otimes^L_{R_\bullet}K\\
\pi^{-1}R\pi_\ast(L\pi^\ast R_0\otimes^L_{R_\bullet}K)
\ar[u]\ar[ru]
}
\end{equation*}
 in $D(T^{\mathbb{N}^{{\rm op}}},\mathbb{Z})$, where 
the vertical arrow is induced from $\pi^{-1}R\to R_\bullet$ and 
the slant one is the counit of adjunction. 
Since $L\pi^\ast R_0\otimes^L_{R_\bullet}K$ is essentially constant, 
the slant one is an essential isomorphism 
and $R\pi_\ast(L\pi^\ast R_0\otimes^L_{R_\bullet}K)$ is bounded 
by Lemma \ref{eses}.3. 
Since the cohomology sheaves of $ R\pi_\ast(L\pi^\ast R_0\otimes^L_{R_\bullet}K)
\cong R_0\otimes^L_RR\pi_\ast K$ 
are sheaves of $R_0$-modules, the vertical one is an essential isomorphism 
by Lemma \ref{eses1}.1. The assertion follows.

$4$. We show $(a)\Rightarrow(b)$. 
Let $K$ be a normalized $R$-complex. Then it is an $R$-complex by $2$ and $L\pi^\ast R\pi_\ast K$ is bounded by $3$. 
Applying Lemma \ref{eses1}.3 to the mapping cone of 
$L\pi^\ast R\pi_\ast K\to K$, we know that it suffices to show that 
$R_0\otimes^L_{R_\bullet}L\pi^\ast R\pi_\ast K\to
R_0\otimes^L_{R_\bullet}K$ is an isomorphism. 
We have 
\begin{equation*}
R_0\otimes^L_{R_\bullet}L\pi^\ast R\pi_\ast K\cong R_0\otimes^L_R\pi^{-1}R\pi_\ast K\cong\pi^{-1}R\pi_\ast(
L\pi^\ast R_0\otimes^L_{R_\bullet}K)\cong
\pi^{-1}R\pi_\ast(
\pi^{-1} R_0\otimes^L_{R_\bullet}K). 
\end{equation*}
Here the last map is an isomorphism since 
$L\pi^\ast R_0\otimes^L_{R_\bullet}K\to 
\pi^{-1} R_0\otimes^L_{R_\bullet}K$ is an essential isomorphism by $1$ and $R\pi_\ast$ kills essentially zero complexes by Lemma \ref{eses}.2.; we also use the fact that $
\pi^{-1} R_0\otimes^L_{R_\bullet}K$ is bounded, which follows from the definition of normalized $R$-complex. As every cohomology sheaf of $\pi^{-1} R_0\otimes^L_{R_\bullet}K$ is of the form $\pi^{-1}N$ for some sheaf $N$ on $T$, the counit $\pi^{-1}R\pi_\ast(
\pi^{-1} R_0\otimes^L_{R_\bullet}K)\to
\pi^{-1}R_0\otimes^L_{R_\bullet}K$ of the adjunction is an isomorphism by Lemma \ref{eses}.1.

$(b)\Rightarrow(c)$ is obvious. We show $(c)\Rightarrow(a)$. If $M$ is bounded, then this  implication and the last assertion are clear. We show that there exist morphisms of complexes of 
sheaves of $R$-modules on $T$ 
\begin{equation}\label{ab}
M\xrightarrow{\alpha} M'\xleftarrow{\beta} M''
\end{equation}
such that $M''$ is bounded and the cohomology sheaves of the mapping cones of $\alpha,\beta$ are uniquely divisible by a uniformizer $\varpi\in R$. 
For $N\in D(T,R)$ with uniquely divisible cohomology sheaves, the pull-back $L\pi^\ast N$ is acyclic. 
Hence $L\pi^\ast M$ is quasi-isomorphic to $L\pi^\ast M''$, which implies that 
we may assume that $M$ is bounded by replacing $M$ by $M''$. 

We construct (\ref{ab}). Since 
$R_0\otimes_R^LM\cong i_0^\ast 
L\pi^\ast M$ is bounded, 
 the cohomology sheaves $H^i(M)$ are 
uniquely divisible for any $i\in\mathbb{Z}$ whose absolute value is large enough. Indeed, as 
$R\xrightarrow{\varpi}R$ is a flat resolution of $R_0$, we have a 
distinguished triangle 
\begin{equation*}
M\xrightarrow{\varpi}M\to R_0\otimes_R^LM\to. 
\end{equation*}
Let $n$ be a positive integer such that $H^i(R_0\otimes_R^LM)$ is zero when $|i|\geq n$. Then 
the multiplication-by-$\varpi$ map 
$H^i(M)\to H^i(M)$ is an isomorphism when $|i|\geq n+1$. 
Therefore, considering truncations of $M$,  we find 
(\ref{ab}) with the desired properties. 
\qed
}

Let $\mathcal{A}$ and $D_{\rm norm}(T^{\mathbb{N}^{\rm op}},R_\bullet)$ be 
the full subcategories of $D^b(T^{\mathbb{N}^{\rm op}},R_\bullet)$ consisting of $R$-complexes and normalized $R$-complexes respectively. 
We also let 
$\mathcal{B}$ be the full subcategory of $D^b(T^{\mathbb{N}^{\rm op}},R_\bullet)$ consisting of complexes which are essentially zero when regarded as objects in $D(T^{\mathbb{N}^{\rm op}},\mathbb{Z})$.

Since $\mathcal{B}$ is a thick triangulated subcategory of 
$D^b(T^{\mathbb{N}^{\rm op}},R_\bullet)$, the quotient $D^b(T^{\mathbb{N}^{\rm op}},R_\bullet)/
\mathcal{B}$ is a triangulated category. By Lemma \ref{eses1}.2, we have an inclusion $\mathcal{B}\subset\mathcal{A}$. 
Let $D^b(T-R)$ be the quotient category $\mathcal{A}/\mathcal{B}$. 
Since the subcategory of essentially constant complexes is stable under 
extensions (which follows from Lemma \ref{eses}.3) and the shift functor, $D^b(T-R)$ is a triangulated subcategory of $
D^b(T^{\mathbb{N}^{\rm op}},R_\bullet)/
\mathcal{B}$. 

\begin{lm}\label{hatab}
\begin{enumerate}
\item For $K\in\cal A$, put $\hat{K}:=L\pi^\ast R\pi_\ast 
K$. Then this is a normalized $R$-complex and is acyclic for 
$K\in\cal B$. 
\item For $K\in\cal A$, $R_n\otimes^L_RR\pi_\ast K$ belongs to $D^b(T,R_n)$ and is acyclic when $K\in\cal B$. 
\end{enumerate}
\end{lm}
\proof{
$1$. The first assertion follows from Lemma \ref{normR}.3.~and 4. The second one follows from Lemma \ref{eses}.2. 

$2$. This is a special case of $1$, since $R_n\otimes^L_RR\pi_\ast K\cong i_n^{-1}\hat{K}$. 
\qed}

\begin{df}\label{alpha}
\begin{enumerate}
\item 
By Lemma \ref{hatab}.1, the functor ${\cal A}\to 
D(T^{\mathbb{N}^{\rm op}},R_\bullet)$ 
sending $K\mapsto\hat{K}:=L\pi^\ast R\pi_\ast 
K$ induces a functor 
$ D^b(T-R)\to D_{\rm norm}(
T^{\mathbb{N}^{\rm op}},R_\bullet)$, for which we write $\Phi$. 
\item By Lemma \ref{hatab}.2, the assignment $K\mapsto R_n\otimes^L_RR\pi_\ast K$ induces a functor $ D^b(T-R)\to D^b(T,R_n)$. We write $R_n
\otimes^L_RK$ for the image of $K\in D^b(T-R)$. 
\end{enumerate}
\end{df}
By Lemma \ref{normR}.2, we can define the functor $D_{\rm norm}(
T^{\mathbb{N}^{\rm op}},R_\bullet)\to D^b(T-R)$ to be the composition of 
$D_{\rm norm}(T^{\mathbb{N}^{\rm op}},R_\bullet)\to\mathcal{A}
\to \mathcal{A}/\mathcal{B}=D^b(T-R)$. 
\begin{lm}\label{essequ}
The functor $D_{\rm norm}(T^{\mathbb{N}^{\rm op}},R_\bullet)\to D^b(T-R)$ 
is an equivalence of categories with a quasi-inverse $\Phi$. 
\end{lm}
\proof{
We show that the compositions of the two functors are isomorphic to the 
identity functors. 
For a normalized $R$-complex $K$, we know that $\Phi(K)=L\pi^\ast R\pi_\ast K\to K$ 
is an isomorphism by Lemma \ref{normR}.4. 

Let $K$ be an $R$-complex. The map $L\pi^\ast R\pi_\ast K\to K$ is an essential isomorphism by Lemma \ref{normR}.3. Therefore, it is an isomorphism in 
$D^b(T-R)$. 
\qed
}

We impose a finiteness condition on (normalized) $R$-complexes. 
From now on, we assume that $T$ is a noetherian topos. 
Let $D^b_c(T,R_0)$ be the full subcategory of $D^b(T,R_0)$ 
consisting of bounded complexes whose cohomology sheaves are constructible. 
\begin{df}\label{consadic}
\begin{enumerate}
\item We denote by $D_{c,{\rm norm}}(T^{\mathbb{N}^{\rm op}},R_\bullet)$ 
the full subcategory of $D_{\rm norm}(T^{\mathbb{N}^{\rm op}},R_\bullet)$ 
consisting of $K\in D_{\rm norm}(T^{\mathbb{N}^{\rm op}},R_\bullet)$ such that 
$i^\ast_0K\in D^b_c(T,R_0)$. 
\item We denote by $D^b_c(T,R)$ the full subcategory of $D^b(T-R)$ 
consisting of $K\in D^b(T-R)$ such that $R_0\otimes^L _RK\in D^b_c(T,R_0)$ (the definition of $R_0\otimes^L _RK$ is given in Definition \ref{alpha}.2). 
An element of $D^b_c(T,R)$ is called  {\rm a constructible complex of $R$-sheaves}. 
\end{enumerate}
\end{df}

\begin{lm}
 The functor in Lemma \ref{essequ} induces an equivalence 
$D_{c,{\rm norm}}(T^{\mathbb{N}^{\rm op}},R_\bullet)\cong
D^b_c(T,R)$. 
\end{lm}
\proof{
It follows from Lemma \ref{essequ} and an 
isomorphism $i_0^\ast K\cong R_0\otimes^L_RR\pi_\ast K$ for $K\in 
D_{\rm norm}(T^{\mathbb{N}^{\rm op}},R_\bullet)$. 
\qed
}

Let $\ell$ be a prime number and let $\Ql$ be an algebraic closure of $\mathbb{Q}_\ell$. For a finite subextension $E/\mathbb{Q}_\ell$ in $\Ql$, we write 
$\mathcal{O}_E$ for the ring of integers in $E$. 
Define the category 
$D^b_c(T,\Zl)$ to be the $2$-colimit 
$\varinjlim_ED^b_c(T,\mathcal{O}_E)$ where 
$E$ runs through the finite subextensions of $\Ql/\mathbb{Q}_\ell$.

\section*{Acknowledgement}
The author would like to express sincere gratitude to his advisor Professor Takeshi Saito for 
many  valuable suggestions and for his warm and patient encouragement. 
He is also grateful to Quentin Guignard for suggesting that 
the results be generalized to the case over a general perfect field. The author would like to thank the anonymous referees for their careful reading and helpful comments. 
This work was supported in part by the Program for Leading Graduate Schools, MEXT, Japan, by 
JSPS KAKENHI Grant Numbers 19J11213 and  25KJ0122, and by 
RIKEN Special Postdoctoral Researcher Program.

\end{document}